\documentclass[reqno]
{amsart}

\usepackage{url}
\usepackage{latexsym}
\usepackage{lscape} 
\usepackage{pdflscape}
\usepackage{float}
\usepackage{amssymb}
\usepackage{amscd}

\usepackage{tikz,tikz-cd}
\usetikzlibrary{arrows,positioning,automata,shadows,fit,shapes}

\usepackage{faktor}

\usepackage[pdftex,
unicode=true,
colorlinks=false
]{hyperref}

\usepackage{subfiles} 

\newtheorem{thm}[equation]{Theorem}
\newtheorem{lem}[equation]{Lemma}
\newtheorem{cor}[equation]{Corollary}

\newtheorem{prop}[equation]{Proposition}

\theoremstyle{remark}
\newtheorem{rem}[equation]{Remark}

\newtheorem{ex}[equation]{Example}

\theoremstyle{definition}
\newtheorem{Def}[equation]{Definition}

\numberwithin{equation}{section}

\newcommand{\rank}{\operatorname{rank}}
\newcommand{\End}{\operatorname{End}}
\newcommand{\Hom}{\operatorname{Hom}}

\newcommand{\Ad}{\operatorname{Ad}}
\newcommand{\ad}{\operatorname{ad}}

\newcommand{\vspan}{\operatorname{span}}

\newcommand{\fa}{{\mathfrak a}}             

\newcommand{\fg}{{\mathfrak g}}

\newcommand{\fh}{{\mathfrak h}}
\newcommand{\fk}{{\mathfrak k}}
\newcommand{\fl}{{\mathfrak l}}

\newcommand{\fn}{{\mathfrak n}}

\newcommand{\fp}{{\mathfrak p}}
\newcommand{\fo}{{\mathfrak o}}
\newcommand{\fq}{{\mathfrak q}}
\newcommand{\ft}{{\mathfrak t}}

\newcommand{\cal}{\mathcal}

\newcommand{\be}{\begin{equation}}
\newcommand{\beu}{\begin{equation*}}

\newcommand{\diag}{\operatorname{diag}}
\newcommand{\bbar}{\,\big|\,}

\newcommand{\pr}{\operatorname{pr}}

\usepackage{graphicx,color,cancel}

\newcommand\eps{{\varepsilon}}

\newcommand{\twedge}{{\textstyle{\bigwedge}}}

\newcommand{\Kt}{\widetilde{K}}

\newcommand\gr{\operatorname{gr}}

\newcommand\im{\operatorname{im}}

\newcommand\id{\operatorname{id}}
\newcommand\Id{\operatorname{id}}
\newcommand\Spin{\operatorname{Spin}}

\def\frT{\mathfrak{T}}

\newcommand{\pf}{\begin{proof}}
\newcommand{\epf}{\end{proof}}
\newcommand{\eq}{\begin{equation}}
\newcommand{\eeq}{\end{equation}}
\newcommand{\eqn}{\begin{equation*}}
\newcommand{\eeqn}{\end{equation*}}

\newcommand{\fra}{\mathfrak{a}}

\newcommand{\frg}{\mathfrak{g}}
\newcommand{\frh}{\mathfrak{h}}
\newcommand{\frk}{\mathfrak{k}}

\newcommand{\frn}{\mathfrak{n}}

\newcommand{\frp}{\mathfrak{p}}

\newcommand{\frt}{\mathfrak{t}}

\newcommand{\frsl}{\mathfrak{sl}}

\newcommand{\frso}{\mathfrak{so}}

\newcommand{\bbC}{\mathbb{C}}

\newcommand{\bbZ}{\mathbb{Z}}

\newcommand{\even}{\operatorname{even}}
\newcommand{\odd}{\operatorname{odd}}

\newcommand{\smat}{\left(\begin{smallmatrix}}
\newcommand{\esmat}{\end{smallmatrix}\right)}
\newcommand{\pmat}{\begin{pmatrix}}
\newcommand{\epmat}{\end{pmatrix}}

\newcommand{\fsl}{\mathfrak{sl}}
\newcommand{\fso}{\mathfrak{so}}
\newcommand{\fgl}{\mathfrak{gl}}

\newcommand{\fsp}{\mathfrak{sp}}
\newcommand{\fe}{\mathfrak{e}}
\newcommand{\ff}{\mathfrak{f}}

\newcommand{\tB}{\widetilde{B}}

\DeclareMathOperator{\Cl}{Cl}

\DeclareMathOperator{\HC}{hc}
\DeclareMathOperator{\dd}{d}
\DeclareMathOperator{\ev}{ev}

\newcommand{\bas}{\operatorname{bas}}
\newcommand{\hor}{\operatorname{hor}}

\newcommand{\sft}{\mathsf{t}}

\DeclareMathOperator{\Der}{Der}
\DeclareMathOperator{\Res}{Res}
\DeclareMathOperator{\sym}{sym}

\newcommand{\rel}{\mathrm{rel}}

\newcommand{\cA}{\mathcal{A}}
\newcommand{\cB}{\mathcal{B}}

\newcommand{\cW}{\mathcal{W}}

\newcommand{\what}{\widehat}
\newcommand{\wtilde}{\widetilde}
\newcommand{\cF}{\mathcal{F}}

\newcommand{\ul}{\underline}
\newcommand{\tW}{\widetilde{W}}

\title[Clifford algebra analogue of the Cartan theorem]{
  Clifford algebra analogue of Cartan's theorem for symmetric pairs}
\author{Kieran Calvert}
\address[K.~Calvert]{Department of Mathematics and Statistics, Lancaster University, Bailrigg, Lancaster, UK}
\email{kieran.calvert@lancaster.ac.uk}

\author{Karmen Grizelj}
\address[K.~Grizelj]{Department of Mathematics, Faculty of Science, University of Zagreb, Bijeni\v{c}ka cesta 30, 10000 Zagreb, Croatia}
\email{karmen.grizelj@math.hr}

\author{Andrey Krutov}
\address[A.~Krutov]{Mathematical Institute of Charles University, Sokolovsk\'a 83, 186 75 Prague, Czech Republic} 
\email{andrey.krutov@matfyz.cuni.cz}

\author{Pavle Pand\v{z}i\'c}
\address[P.~Pand\v{z}i\'c]{Department of Mathematics, Faculty of Science, University of Zagreb, Bijeni\v{c}ka cesta 30, 10000 Zagreb, Croatia}
\email{pandzic@math.hr}

\thanks{K.~Grizelj and P.~Pand\v{z}i\'{c} were supported by the QuantiXLie Center of Excellence, a project
  co-financed by the Croatian Government and European Union through the European Regional Development Fund -
  the Competitiveness and Cohesion Operational Programme (grant PK.1.1.10.0004) 
  P.~Pand\v{z}i\'{c} was supported by the Croatian Science Foundation (HRZZ), grant no. IP-2025-02-6514, and by the European Union – NextGenerationEU through the National Recovery and Resilience Plan 2021-2026. Institutional grant of University of Zagreb Faculty of Science (IK IA 1.1.3. Impact4Math).  
  A.~Krutov was supported by the GA\v{C}R projects~24-10887S and~26-20231L.
  This article is based upon work from COST Action CaLISTA CA21109 supported by COST (European
  Cooperation in Science and Technology, \url{www.cost.eu}).
}



\begin{document}

\begin{abstract}
  We extend Kostant's results about $\fg$-invariants in the Clifford algebra~$\Cl(\fg)$ of a~complex semisimple
  Lie algebra~$\fg$ to the relative case of $\fk$-invariants in the Clifford algebra $\Cl(\fp)$, where
  $(\fg,\fk)$ is a~classical symmetric pair and $\fp$ is the $(-1)$-eigenspace of the corresponding involution.
  In this setup we prove 
  the Cartan theorem for Clifford algebras,
  a~relative transgression theorem,
  the Harish--\/Chandra isomorphism for~$\Cl(\fp)$, and a~relative version of Kostant's Clifford algebra
  conjecture. 
\end{abstract}

\maketitle

\section{Introduction}

Let $G$ be a~compact simple connected Lie group and $K$ a~closed subgroup (hence a Lie subgroup).
We do not assume $K$ is connected, but we will be assuming that~$K$ is a~symmetric subgroup, i.e., there is an~involution~$\Theta$ of~$G$ such that
\[
G^\Theta_0\subseteq K\subseteq G^\Theta,
\]
where $G^\Theta$ denotes the subgroup of $G$ consisting of the points fixed under $\Theta$, and $G^\Theta_0$
denotes the connected component of $G^\Theta$.
Typical examples of compact symmetric spaces $G/K$ are the classical ones, such as
$SU(p+q)/S(U(p)\times U(q))$, $SO(p+q)/S(O(p)\times O(q))$, $Sp(p+q)/Sp(p)\times Sp(q)$.
Among these, $K=S(O(p)\times O(q))$ is not connected.
Note that in all these examples $K = G^\Theta$ for an appropriate~$\Theta$.

By well known classical results of \'E.~Cartan and de Rham, the de Rham cohomology (with complex coefficients)
$H(G/K)$ of~$G/K$ as above can be described as follows.
Let $B$ be a~nondegenerate invariant symmetric bilinear form on the complexified Lie algebra~$\frg$ of~$G$;
such a~form exists since $G$ is reductive, so one can take the Killing form on the semisimple part of~$\frg$
and extend by any nondegenerate form on the center of~$\frg$.
Let $\frk\subseteq\frg$ be the complexified Lie algebra of~$K$ and let~$\frp$ be the orthogonal complement of~$\frk$ with respect to~$B$, so that
\[
  \frg=\frk\oplus\frp,
\]
and this decomposition is compatible with the adjoint $K$-action. 
The space~$\frp$ can also be described as the $(-1)$-eigenspace of $\theta = \dd\Theta$.
Denote by $(\twedge\frp^\ast)^K$ the algebra of $K$-invariants in $\twedge\frp^\ast$.
Then
\[
H(G/K)\cong(\twedge\frp^\ast)^K
\]
as graded algebras.

As mentioned above, this result is old and very well known, but it is not easy to find a~proof in the
literature. One place where one can find it is~\cite[\S 4.1]{CNP}.

The structure of the algebra $(\twedge\frp^\ast)^\frk$ is well understood through the work of Borel, H.~Cartan,
Hopf, Koszul, Samelson and others, see~\cite{Cartan1951trans,Cartan1951,Borel1953,Michel1962,Rashevskii1969}.
This is achieved using a~series of quasi-isomorphisms and it is quite involved and
not completely explicit; moreover, there is no treatment of $(\twedge\frp^\ast)^K$ for disconnected~$K$.
See~\cite[12.3 Theorem~1]{Onishchik1994} and~\cite[X.2 Theorem~III]{CCCvol3}.

On the other hand, one can replace the algebra $\twedge\frp$ (which is isomorphic to $\twedge\frp^{\ast}$ via~$B$) by its filtered deformation (or  quantisation), the Clifford algebra $\Cl(\frp)$ with respect to $B$. This is an associative algebra with unit, generated by $\frp$, with relations
\[
XY+YX=2B(X,Y),\qquad X,Y\in\frp.
\]
Clearly, $\Cl(\frp)$ is a quotient of the tensor algebra $T(\frp)$ by a nonhomogeneous ideal, so it inherits a filtration, and the corresponding graded algebra is $\twedge\frp$. 

One might think that the filtered algebras $\Cl(\frp)$ and $\Cl(\frp)^K$ are more complicated than their associated graded algebras $\twedge\frp$ and $(\twedge\frp)^K$, but in many aspects the filtered versions are in fact simpler. To illustrate this, let us temporarily assume that $\frg$ and $\frk$ have equal rank, and let $\frt$ be a common Cartan subalgebra. Then $\dim\frp$ is even, and the Clifford algebra $\Cl(\frp)$ has only one simple module, the spin module $S$; moreover, any $\Cl(\frp)$-module is a direct sum of copies of $S$, and $\Cl(\frp)\cong\End S$. To construct $S$, write
\[
\frp=\frp^+\oplus\frp^-,
\]
where $\frp^+$ and $\frp^-$ are maximal isotropic subspaces of $\frp$ in duality under $B$. Then one can take $S=\twedge\frp^+$, with elements of $\frp^+$ acting by wedging, and elements of $\frp^-$ by contracting. The spin module $S$ carries an action of the spin double cover $\Kt$ of $K$; the corresponding action of $\frk$ is described
by the Lie algebra map $\alpha_\frp:\frk\to \Cl(\frp)$,
\[
\alpha_\frp(X)=\frac{1}{4}\sum_i [X,b_i]d_i,\qquad X\in\frk,
\]
where $b_i$ is any basis of $\frp$ and $d_i$ is the dual basis with respect to $B$.
As a $\frk$-module, $S$ is multiplicity free and decomposes as
\[
S=\bigoplus_{w\in W^1} E_{w\rho_\frg-\rho_\frk}.
\]
Here $\rho_\frg$ and $\rho_\frk$ are the half sums of compatible positive roots of $(\frg,\frt)$ respectively $(\frk,\frt)$. Furthermore, 
\[
W^1=\{w\in W_\frg\bbar w\rho_\frg\text{ is }\frk\text{-dominant}\},
\]
where $W_\frg$ denotes the Weyl group of $(\frg,\frt)$. Finally, 
 $E_{w\rho_\frg-\rho_\frk}$ denotes the irreducible finite-dimensional $\frk$-module with highest weight $w\rho_\frg-\rho_\frk$.

Now $\Cl(\frp)\cong\End S$ implies 
\[
\Cl(\frp)^\frk=\End_\frk S,
\]
and by Schur's Lemma this is the algebra of projections $\pr_w: S\to E_{w\rho_\frg-\rho_\frk}$. This is a very simple commutative algebra, isomorphic to $\bbC^{|W^1|}$ with coordinate-wise multiplication. (The isomorphism is given by identifying $\pr_w$ with $(0,\dots,0,1,0,\dots,0)$, with $1$ in the place corresponding to $w$.)
This was explained to the fourth-named author by Kostant \cite{kostantemail}; see~\cite{LRpaper} and~\cite{CNP} for more details. We can now extend our map $\alpha_\frp$ to the universal enveloping algebra $U(\frk)$ of $\frk$, restrict it to the center $Z(\frk)$ of $U(\frk)$, and identify $Z(\frk)$ with $\bbC[\frt^*]^{W_\frk}$ by the Harish-Chandra isomorphism (here $W_\frk$ is the Weyl group of $(\frk,\frt)$). In this way we get a map 
\[
\alpha_\frp:\bbC[\frt^*]^{W_\frk}\to \Cl(\frp)^\frk,\qquad \alpha_\frp(P)=\sum_{w\in W^1} P(w\rho)\pr_w,
\]
which is easily seen to be onto (in the equal rank case). Moreover, the kernel of $\alpha_\frp$ is generated by the $W_\frg$-invariants in $\bbC[\frt^*]$ that evaluate to 0 at $\rho$, and thus one gets to understand $\Cl(\frp)^\frk$ by generators and relations. One can then pass to associated graded algebras and obtain results about the structure of $(\twedge\frp)^\frk$. This approach is much more direct and explicit than the classical approach described earlier, and has an additional advantage of simultaneously understanding the filtered deformation $\Cl(\frp)^\frk$ of $(\twedge\frp)^\frk$. Finally, it is not difficult to extend these considerations to the algebras $\Cl(\frp)^K$ and $(\twedge\frp)^K$ for disconnected $K$; see \cite{CNP} and \cite{LRpaper}.

If $\frg$ and $\frk$ do not have equal rank, $\Cl(\frp)^\frk$ is typically much bigger than the algebra  $\Pr(S)$ of projections  of the spin module; moreover, it is not abelian any more. In fact, it is the tensor product of the algebra of projections and a Clifford algebra isomorphic to $\Cl(\fra)$, where $\fra$ is the split part of the fundamental Cartan subalgebra $\frh=\frt\oplus\fra$ of $\frg$.
The situation is thus analogous to the case of $(\twedge\frp)^\frk$ studied by Borel, H.~Cartan, Hopf, Koszul and Samelson. The Clifford algebra setting was first studied by Kostant \cite{KostantRho} who proved that in the ``absolute" case when $\frk$ and $\frp$ are both isomorphic to a complex Lie algebra $\frg$ (and what we denoted by $\frg$ earlier is now $\frg\oplus\frg$), $\Cl(\frg)^\frg$ is isomorphic to the Clifford algebra $\Cl(\frh)$, where $\frh$ is a Cartan subalgebra of $\frg$, and this copy of $\Cl(\frh)$ is realized inside $\Cl(\frg)$ as a Clifford algebra over the primitive elements $P_\wedge(\frg)$ which correspond to $\frh\subset\Cl(\frh)$, but have various high degrees. This is analogous to the classical Hopf-Koszul-Samelson theorem about $(\twedge\frg)^\frg$.

The purpose of this article is to prove results analogous to Kostant's in the ``relative" case. We show that for a symmetric pair $(G,K)$ as above, the algebra $\Cl(\frp)^\frk$ is isomorphic to a tensor product of the algebra of projections of the spin module with a Clifford algebra isomorphic to $\Cl(\fra)$. The subalgebra $\Cl(\fra)$ is realized inside $\Cl(\frp)$ as the Clifford algebra over the ``primitives" $P_\wedge(\frp)$ corresponding to $\fra$ but lying in various high degrees. Explicitly $P_\wedge(\frp)$ is the projection of $P_\wedge(\frg)$ to $(\twedge\frp)^\frk$. We denote by $P_{\Cl}(\frp)$ the image of $P_\wedge(\frp)$ under the Chevalley map (skew-symmetrization) $q:\twedge\frp\to\Cl(\frp)$. 
We define a form $\tilde{B}$ on the vector space $(\twedge\frp)^\frk\cong \Cl(\frp)^\frk$ by setting  $\tilde{B}(a,b)$ to be the 0th component of $\iota_a b$, where $\iota_a$ is the contraction of $\twedge\frp$ by $a$.  We
can now state our main theorem: 

\begin{thm}\label{thm:main}
  Let $(G,K)$ be a~compact symmetric pair such that $G$ is simple and connected.  
  Assume that $(\fg,\fk)$ is different from $(\fe(6),\fsp(8))$.
 
  (a) With the above notation, the inclusion $P_{\Cl}(\frp) \hookrightarrow \Cl(\fp)^K$ extends to a~filtered algebra homormorpism 
  \[ \Cl(P_{\Cl}(\frp), \tilde{B}) \to \Cl(\fp)^K,\]
  which is an~isomorphism in the primary cases, i.e., when the spin module $S$ contains only one $\frk$-type.
  For the ``almost primary cases" $(G,K)=(SU(2n),SO(2n))$, the statement remains true if $K$ is replaced by~$O(2n)$.
  
  (b) There is an algebra isomorphism 
  \[     \Cl(P_{\Cl}(\frp),\tilde{B}) \otimes \Pr(S)   \cong \Cl(\fp)^\frk.\]
  
  (c) There is a~filtered algebra isomorphism 
  \[ \Pr(S) \cong \bbC[\ft^*]^{W_\frk} / I_\rho,\]
  where $I_\rho$ is the ideal of $\bbC[\ft^\ast]^{W_\fk}$ generated by $a\in \bbC[\ft^\ast]^{W_\fg}$ such that
  $a(\rho) = 0$.
  
  (d) Passing to the associated graded algebras we recover the well known Cartan (or Cartan-Borel) Theorem for symmetric spaces 
  \[\twedge P_\wedge(\frp) \otimes \bbC[\ft^*]^{W_\frk}/I_+   \cong (\twedge\fp)^\frk,\]
  where $I_+ = \gr I_\rho$ is the ideal of $\bbC[\ft^\ast]^{W_{\fk}}$ generated by elements
  of~$\bbC[\ft^\ast]^{W_\fg}$ with zero constant term.
\end{thm}

For the equal rank case, Theorem~\ref{thm:main}(a) is trivial and (b)
was proved by Kostant (see above).
For the rest, see~\cite[Section 2]{LRpaper} and~\cite[Section 3.4]{CNP}.
For $(\frg,\frk)=(\frso(2p+2q+2,\bbC),\frso(2p+1,\bbC)\oplus\frso(2q+1,\bbC))$ Theorem 1.1 is proved
in~\cite[Section 3.5]{CNP}; this case is called almost equal rank because $\dim\fra=1$.
In this paper we prove this result in the remaining cases, the primary and almost primary cases. 

One of Kostant's main tools in the absolute case was the Transgression Theorem, which allows one to construct primitive invariants in $\twedge\frg$ starting from the better known primitive invariants in the symmetric algebra $S(\frg)$. Analogues of this theorem play also an important role in our approach. 
Kostant also used the Hodge decomposition and the corresponding duality for Lie algebra (co)homology.
This is not available in the relative case, so we developed new tools based on transgression and
Harish--\/Chandra maps.
This also gives an~alternative way to prove Kostant's results.

We emphasise that although our main interests are (relative) quantum and commutative Weil algebras, the proofs of our main results (Theorem~\ref{thm:main} and, in particular, Theorem~\ref{thm:IotaTClP})
fundamentally utilise the (relative) noncommutative Weil algebra. The reason why noncommutative Weil algebras are
crucial in our approach is the fact that they are universal algebras in the category of locally free $G$-differential
algebras, see~\cite{AlekseevMeinrenken2005}.

The noncommutative Weil algebra of a~Lie algebra~$\fg$ is defined as 
\[
\tW(\frg)=T(\frg^\ast[-2]\oplus\frg^\ast[-1]),
\]
see~\cite{AlekseevMeinrenken2005} and Section~\ref{subsec:noncomm Weil}.
We revisit the transgression construction and show that it factorises via~$\tW(\fg)$.
This allows us to show that the image of the transgression map consists of ``primitives''
  a~la~\cite{Rashevskii1969} and obtain an alternative proof of Kostant's results which is more suitable for our generalisations.
  
In a similar manner we define transgression maps for the relative Weil algebras 
\[
W(\frg,K)=W(\frg)_{K-\text{basic}}=(S(\frg^\ast)\otimes\twedge\frp^\ast)^K
\]
and prove the relative
  version of the transgression theorem, Theorem \ref{thm:relTransgression}. 
  It says that we can produce primitive $K$-invariants in $\twedge\frp^\ast$ starting from those primitive $\frg$-invariants in $S(\frg^\ast)$ that project trivially to $S(\frk^\ast)^K$.

We also study other related absolute and relative Weil algebras, as well as corresponding transgression maps. For example, we consider the quantum Weil algebra
\[
\cW(\frg)=U(\frg)\otimes\Cl(\frg)
\]
and its relative version
\[
\cW(\frg,K)=(U(\frg)\otimes\Cl(\frp))^K.
\]
To relate the results about transgression and Theorem \ref{thm:main}, we use certain Harish-Chandra type
projections of the Clifford algebra. In the absolute case, a Harish-Chandra map $\HC$ from $\Cl(\frg)^\frg$ to
$\Cl(\frh)$ was studied by Kostant, Bazlov and others; see \cite{MeinrenkenBook}. In the relative cases, we
define analogues of this map  attached to elements of $W^1$. We show that each
$\HC_w:\Cl(\frp)^\frk\to\Cl(\fra)$ is a surjective algebra map, and that~$\HC_{\id}$
is an~isomorphism in the primary cases, equivalently, $W^1 = \{\id\}$; see Propositions \ref{prop:HCsurjective} and \ref{prop:HCisIsoInPrimary}.

In the absolute cases, it was conjectured by Kostant and proved by Bazlov~\cite{BazlovPhD}, Joseph~\cite{Joseph2012KostantClifford} and Alekseev-Moreau~\cite{AlekseevMoreau2012}, that $\HC$ is an isomorphism of filtered algebras, if we introduce a filtration of $\frh$ using the action of a principal $\frsl(2)$ subalgebra of the dual Lie algebra $\frg^{\vee}$. We show that an analogous result holds in primary and almost primary cases, see Theorem \ref{thm:relKostCliffAlgConj}. 

Proposition \ref{prop:HCsurjective}, Proposition \ref{prop:HCisIsoInPrimary} and Theorem \ref{thm:relKostCliffAlgConj} can be summarized as follows:

\begin{thm} \label{thm:RelKostConj}
Let $(G,K)$ be a symmetric pair as in Theorem \ref{thm:main}. Let $\fra$ be filtered using the action of the principal $\frsl(2)$ subalgebra of $\frg^\vee$.
Then for each $w\in W^1$, $\HC_w:\Cl(\frp)^\frk\to\Cl(\fra)$ is a surjective filtered algebra homomorphism. In the primary case, $\HC_1$ is an isomorphism. 
\end{thm}

The paper is organized as follows. After reviewing some preliminaries in Section~2, we discuss various Weil
algebras and transgression maps in Section~3 and~4. Section~5 is devoted to the proof of Relative
Transgression Theorem~\ref{thm:relTransgression} and in Section~6 we study the Harish-Chandra maps. In Section~7 we prove the main
theorem, Theorem~\ref{thm:main}. Finally, in Section~8 we discuss the relative version of Kostant's Clifford
algebra conjecture.

\subsection*{Acknowledgements}
We thank Anton Alekseev and Kyo Nishiyama for stimulating discussions and continued interest in our work.
We are also grateful to the referee for careful reading and many useful suggestions that led to significant
improvement of the manuscript.

\section{Preliminaries}
In this section we survey some  definitions and facts about Clifford algebras and symmetric pairs. 

\subsection{Contractions and bilinear forms}
Let $V$ be a vector space with a nondegenerated symmetric bilinear form $B$. 

The transpose anti-automorphism $\bullet^T$ of $\twedge V$ is defined by sending
$a = v_1 \wedge \ldots \wedge v_k\in \twedge V$ to  
$a^T = v_k \wedge \ldots \wedge v_1$. Since~$\twedge V$ is super-commutative, $\bullet^T$ is multiplication by~$(-1)^{k(k-1)/2}$ on each graded component $\twedge^k V$. 

\begin{Def}\label{def:contractions}
  For $f \in V^\ast$, we define a~derivation
  $\iota_f: \twedge V \to \twedge V$ by
  \[
    \iota_f (v_1 \wedge \ldots \wedge v_k) = \sum_{i=1}^k(-1)^{i-1}\langle f, v_i\rangle v_1 \wedge \ldots \wedge v_{i-1} \wedge v_{i+1}\wedge\ldots \wedge v_k.
  \]
  The map $\iota$ extends to an~algebra homomorphism $\twedge V^\ast \to \End(\twedge V)$ and is called
  contraction.

  Using~$B$ we can identify $V^\ast$ with~$V$ via $B^\sharp\colon V\to V^\ast$, $B^\sharp(v) = B(v, -)$ and
  define contractions by elements of~$V$. Explicitly,
  $\iota_v = \iota_{B^\sharp(v)}$ for $v\in V$.
\end{Def}
  
\begin{Def}\label{def:biforms}   
  For $a \in \twedge V$ let $(a)_{[k]}\in\twedge^k V$ be the $k^{th}$ degree part of~$a$.
  We extend the bilinear form $B$ to $\twedge V$ in two ways;
  \[ B(a,b) = (\iota_{a^T}b)_{[0]}, \quad \tB(a,b) = (\iota_ab)_{[0]}.\]
\end{Def}

\subsection{Clifford algebras}\label{sec:CliffAlg}
In this section we recall the definition and main properties of Clifford algebras. Our main reference
is~\cite[Chapter~2]{MeinrenkenBook}; see also~\cite[\S2]{HuangPandzicBook}.

\begin{Def}
  For a~vector space~$V$ with symmetric nondegenerate inner product~$B$, define the Clifford algebra 
  \[ \Cl(V, B) = T(V) / \langle u\otimes v + v\otimes  u - 2 B( u, v ) \rangle .\]
\end{Def}
\noindent We often drop the bilinear form in notation and write $\Cl(V,B)$ as $\Cl(V)$.

Depending on the dimension of~$V$, $\Cl(V)$ is either an~endomorpism algebra or two copies of an endomorphism algebra of spin modules $S$ and $S_1,S_2$ respectively. 
\[
  \Cl(V) \cong
  \begin{cases} \End(S) & \text{if $\dim V$ is even},\\
    \End(S_1) \oplus \End(S_2) & \text{if $\dim V$ is odd}.
  \end{cases}
\] 

The associated graded algebra of the Clifford algebra~$\Cl(V)$ is the exterior algebra $\twedge V$.
These are isomorphic as vector spaces.
The isomorphisms are given by the symbol map $\mathrm{symb}:\Cl(V) \to\twedge V$
and the quantization map $q:\twedge V \to \Cl(V)$; see~\cite[Section 2.2.5]{MeinrenkenBook}.

For $v\in V$ we define an odd derivation $\iota_v\colon \Cl(\fg)\to\Cl(\fg)$ by
\[
  \iota_v (v_1\ldots v_k) = \sum_{i=1}^{k}(-1)^{i-1} B(v,v_i) v_1\ldots v_{i-1} v_{i+1}\ldots v_k,
\]
The quantisation map intertwines contraction operators in $\twedge V$ and $\Cl(V)$:
for all $\omega \in \twedge V$, $v\in V$ we have that $\iota_v(q(\omega)) = q(\iota_v\omega)$.

The filtration on~$\Cl(V)$ is compatible with the $\mathbb{Z}$-grading of~$\twedge V$ in the following manner
\begin{equation}\label{eq:Clfiltration}
  \Cl^{(k)}(V) = q\left( \bigoplus_{m=0}^k \twedge^m V \right)
  \qquad\text{for $k=0,\ldots,\dim(V)$}.
\end{equation}
The degree~$2$ space $q(\twedge^2 V)\subseteq \Cl^{(2)}(V)$ is closed under the Lie bracket $[X,Y] = XY - YX$.
The Lie algebra $q(\twedge^2 V)$ is isomorphic to $ \fo (V)$.

Let $L$ be a compact Lie group and $\fl$ be the corresponding Lie algebra (in examples we need, $L$ will be~$G$
or~$K$).
Now assume that $V$ is an~$\fl$-module which admits a~non-degenerate invariant symmetric bilinear form~$B$.
Since $B$ is symmetric and invariant, the action of~$\fl$ on~$V$ defines a~Lie algebra homomorphism  $\fl \to \fo(V)$.
Hence we get an~$\fl$-equivariant linear map
\[
  \begin{tikzcd}
    \lambda_V \colon \fl \arrow[r, ""] & \fo(V) \arrow[r, "\simeq"] & \twedge^2 V.
  \end{tikzcd}  
\]
Composing it with the standard embedding of the Lie algebra~$\fo(V)\cong \twedge^2V$ into the Clifford algebra~$\Cl(V)$, we get a~Lie algebra map 
\begin{equation} \label{eq:def of alpha}
  \begin{tikzcd}
    \alpha_V \colon \fl \arrow[r, ""] & \fo(V) \arrow[r, "\simeq"] & \twedge^2 V \arrow[r, "q"]  & \Cl^{(2)}(V).
  \end{tikzcd}  
\end{equation}
Explicitly, if $e_i$ is a~basis of $V$ with dual basis $f_i$ with respect to~$B$, then
\[
  \lambda_V(X) = \frac{1}{4}\sum_{i}^{(V)}[X,e_i]\wedge f_i,\qquad
  \alpha_V(X) =  \frac{1}{4}\sum_{i}^{(V)} [X,e_i]f_i, \quad X \in \fl,
\]
where $\sum^{(V)}_i$ denotes the summation over a~basis of~$V$ indexed by~$i$,
and $[X,e_i]$ denotes the action of $X$ on~$e_i\in V$ (in all examples it will be the adjoint action, hence
the notation).
By the universal properties, $\lambda_V$ and $\alpha_V$ extend to algebra maps
\[
  \lambda_V\colon S(\fl)\to \twedge V,\qquad
  \alpha_V \colon U(\fl) \to \Cl(V),
\]
(in particular, $\alpha_V$ defines a~quantum moment map in the sense of~\cite{Lu1993}).
Therefore, the spin module $S$ (or $S_1$ and $S_2$) of~$\Cl(V)$ can be viewed as an $\fl$-module.
We define $\alpha_Z = \alpha_V|_{Z(U(\fl))}$ to be the restriction of $\alpha_V$ to the centre of $U(\fl)$.

In a similar manner one defines the corresponding $\tilde{L}$-module structure on $S$ (or $S_1$, $S_2$),
where $\tilde{L}$ is the double cover of~$L$ defined by the following pullback diagram:
\[
  \begin{tikzcd}
    \tilde{L} \arrow[d] \arrow[r] & \Spin(V) \arrow[d] \\
    L \arrow[r] & SO(V, B) 
  \end{tikzcd}
\]
See~\cite[\S 3.2.1]{HuangPandzicBook} for details.

\subsection{Symmetric pairs of Lie algebras}\label{ssec:SymPairs}

Let $(G,K)$ be a~symmetric pair of compact Lie groups as in the introduction, with $\Theta$ being the
corresponding involution of~$G$. Then the complexified Lie algebra~$\fg$ of~$G$ has an~involution~$\theta =
\dd\Theta$ such that $\fg^\theta = \fk$ is the complexified Lie algebra of~$K$.
The pair $(\fg,\fk)$ is called a~\emph{symmetric pair} (of Lie algebras).
Then $\fg$ has a~nondegenerate symmetric bilinear form~$B$ (e.g. the Killing form if~$\fg$ is
semisimple) and this form is invariant under~$\theta$.
The $(-1)$-eigenspace of~$\theta$ is denoted by~$\fp$. It is clear that $\fp$ is a~$\fk$-module.
Let $\ft$ denote a~Cartan subalgebra of~$\fk$ and let $\fh = \ft \oplus \fa$ be a~Cartan subalgebra of~$\fg$
with $\fa\subset\fp$.

\begin{Def}
The Clifford algebra $\Cl(\fp)$ is defined with 
respect to a (fixed) nondegenerate $\frg$-invariant form~$B$ on~$\fg$ restricted to~$\fp$.
We denote by~$S$ (or~$S_1$ and~$S_2$ if $\dim \fp$ is odd) the irreducible (spin) module(s) for $\Cl(\fp)$.
\end{Def}

\subsection{Types of symmetric pairs}

\begin{Def}\label{def:typeOfPairs}
  A symmetric pair~$(\fg,\fk)$ is called
  \begin{itemize}
  \item \emph{primary} if $S$ is $\fk$-primary, i.e., it contains only one $\fk$-type;
  \item \emph{almost primary} if $S$ has two $\fk$-types and $\dim\fa>1$;
  \item \emph{of equal rank} if $\dim\fa=0$;
  \item \emph{of almost equal rank} if $\dim\fa=1$.
  \end{itemize}  
\end{Def}

The classification of compact symmetric spaces~$G/K$ for simple $\fg = \mathrm{Lie}(G)$ is given, for example,
in~\cite{Helgason1962}.
It turns out that all corresponding symmetric pairs of Lie algebras fall into one of the four categories of
Definition~\ref{def:typeOfPairs} except for the pair~$(\fe(6),\fsp(8))$.

In this paper we will consider mostly primary symmetric pairs, such as
\[
  (\fsl(2n+1),\fo(2n+1),\ (\fsl(2n),\fsp(2n)),\ (\fe(6),\ff(4)),
\]
and almost primary symmetric pairs
\[
  (\fsl(2n),\fo(2n)).
\]

\subsection{Algebra structure of $\Cl(\fp)^\fk$ as endomorphisms of spin modules}

Recall that, if $\dim \frp$ is even, the Clifford algebra is isomorphic to the endomorphisms of the spinor module and if $\dim \fp$ is odd, equal to the direct sum of two such spaces. 

\[
  \Cl(\fp) = \begin{cases} \End(S) & \dim \fp \text{ even,} \\ \End(S_1) \oplus \End(S_2) & \dim \fp \text{ odd}.
             \end{cases}
\]

\begin{lem}
    The algebra $\Cl(\fp)^\fk$ is equal to the $\fk$-endomorphisms of $S$ (respectively the direct sum of $\fk$-endomorphisms of~$S_1$ and~$S_2$).
\end{lem}

\begin{proof}
This follows from the definition of the action of $\frk$ on $\End(V)$ for a $\frk$-module~$V$.
\end{proof}

\begin{lem}[\cite{Parthasarathy} or \S2.3.6 in~\cite{HuangPandzicBook}]
  The highest weights of the spin module $S$ (respectively $S_1$ and $S_2$) are given by
  \[ w \rho_\fg-\rho_\fk \quad \text{ for } w \in W^1.\]
\end{lem}

Let $\fp = \fa \oplus \fq$ where $\fq = \fq^+\oplus \fq^- = \fn^+ \cap \fp \oplus \fn^- \cap \fp$;
here $\frn^\pm$ are defined using our fixed choice of positive roots.
Define $E_w$ to be the irreducible finite-dimensional $\fk$-module with highest weight $w \rho_\fg-\rho_\fk$.
Let $\fa^+,\fra^-$ be a~dual pair of maximal isotropic subspaces of~$\fa$ with respect to~$B$. If $\fra$ or
equivalently $\frp$ is odd dimensional, let $\fra_0$ be the $1$-dimensional space orthogonal to $\fra^+\oplus \fra^-$,
spanned by $z$ of norm~$1$.
Then for $\dim\fp$ even, $S$ is isomorphic to $\twedge \fa^+ \otimes \twedge \fq^+$.
For $\dim\fp$ odd, $S_1$, $S_2$ are both isomorphic to $ \twedge \fa^+ \otimes \twedge \fq^+$,
with $z\in\fa_0$ acting on $S_1^{\even}$, $S_2^{\odd}$ by~$1$ and $S_1^{\odd}$, $S_2^{\even}$ by~$-1$.

\begin{lem}\cite[\S2.3.6]{HuangPandzicBook}\label{lem::spindecomp}
  As $\fk$-modules, 
  \[
    S \cong \bigoplus_{w \in W^1} \twedge{\fa^+}\otimes E_w,
  \quad
    S_1  \cong S_2 \cong \bigoplus_{w \in W^1} \twedge \fa^+\otimes E_w,
  \]
  with trivial $\fk$-action on $\twedge \fa^+$.    
\end{lem}

\begin{Def}\label{def:PrSigma}
  Define an~algebra $\Pr(S)$ isomorphic to $\mathbb{C}^{|W^1|}$ with basis of idempotents labeled by
  $\pr_{w}$, $w\in W^1$, such that  
  \[
    \pr_{w_1} \pr_{w_2} = \delta_{w_1,w_2} \pr_{w_1}.
  \]    
\end{Def}

\begin{lem}\label{lem:ClifDecomp} 
  The $\fk$-invariants in $\Cl(\fp)$ are isomorphic as an algebra to 
  \[
    \Cl(\fp)^\fk \cong
    \begin{cases} \End(\twedge\fa^+) \otimes \Pr(S) & \dim \fp \text{ even,} \\
      \left(\End(\twedge\fa^+) \oplus \End(\twedge\fa^+)\right) \otimes \Pr(S) & \dim \fp \text{ odd.}\end{cases} 
  \]
  In both cases
  \[
    \Cl(\fp)^\fk \cong \Cl(\fa) \otimes \Pr(S).
  \]
  Moreover, $\Pr(S) = \im\alpha_Z$.
\end{lem}

\begin{proof}
  The algebra $\Cl(\fp)^\fk$ is the $\fk$-endomorphisms of~$S$ or~$S_1$ and~$S_2$.
  For each $\fk$-type~$E_w$ there is a~projection~$\pr_w$ onto the~$E_w$ isotypic component.
  Furthermore, this isotypic component is a~direct sum of $\dim\twedge\fa^+$ copies of~$E_w$,
  hence the space of $\fk$-endomorphisms of this isotypic component is isomorphic to~$\End(\twedge\fa^+)$.
  An~identical argument can be made for~$S_1$ or~$S_2$.
  To finish, one notes that~$\fq$ is always even dimensional as it contains both positive and negative root
  spaces.
  Hence $\dim \fp$ has the same parity as $\dim \fa$.
  Therefore, if $\dim \fp$ is even, then $\Cl(\fa) = \End(\twedge\fa^+)$ and
  if $\dim\fp$ is odd, then $\Cl(\fa) = \End(\twedge\fa^+) \oplus \End (\twedge\fa^+)$. 

  The fact that $\Pr(S) = \im\alpha_Z$ was proved in~\cite{LRpaper}.
\end{proof}

It is worth noting this is not a~filtered algebra isomorphism if we use the usual filtration on~$\Cl(\fa)$.
A~description of the filtration on~$\Cl(\fa)$ which is compatible with the above isomorphism will be given
in Section~\ref{sec:KostantCliffordConjecture}.

\section{Transgressions in Weil algebras}
The goal of the next several sections is to prove that the contractions of a~primitive element belongs to
the image of the quantum moment maps $\alpha_\fg$, respectively $\alpha_\fp$ in the relative case
of~$\Cl(\fp)^K$.
To show this, we prove that the subspaces of the primitive invariants coincides with the images of the
transgression maps defined for various Weil algebras.
We prove that all transgression maps factor through the noncommutative Weil algebra.
Then using the universal property of the noncommutative Weil algebra we show that contractions of the
image of the transgression map belong to the image of the noncommutative analogue of~$\alpha_\fg$,
respectively of~$\alpha_\fp$, which we denote by $\Phi_{\Cl(\fg)}$, respectively $\Phi^{\mathrm{hor}}_{\Cl(\fg)}$.
The original goal now follows from the analogous statement for the noncommutative Weil algebra, respectively
its relative version, which we prove in Theorem~\ref{thm:iotaTtilde}, respectively
Proposition~\ref{prop:contr_trans_formula}.

In this section we define various Weil algebras and construct the corresponding transgression maps.
The results of this section will allow us to reprove Kostant's  description of $\Cl(\frg)^\frg$ as the Clifford
algebra on primitives. Furthermore, we are able to extend our results to the relative case which will be
illuminated in section~\ref{sec:Cliff Alg analogue}. Theorem~\ref{thm:iotaTtilde} shows that the image of the
transgression map consists of primitives as defined in~\cite{Rashevskii1969}.
It turns out that the transgression maps in various Weil algebras all factor through the noncommutative Weil
algebra.
The universal property of the noncommutative Weil algebra will be crucial for defining various derivations and
understanding their properties.

\subsection{$G$-differential algebras}\label{ssec:GdiffAlg}
H.~Cartan introduced the notion of $\fg$-differential algebras as a~generalisation
of algebras of differential forms on manifolds with $\fg$-action,  \cite{Cartan1951trans,Cartan1951}.
Later $\fg$-differential algebras appeared in the study of equivariant cohomology~\cite{GuilleminSternbergBook,AlekseevMeinrenken2000,AlekseevMeinrenken2005duke},
in Chern--\/Weil theory~\cite{AlekseevMeinrenken2005,MeinrenkenBook}, and in relation to (algebraic) Dirac
operators and Vogan's conjecture~\cite{HP1,HuangPandzicBook,AlekseevMeinrenken2005}.
Here we review the definition and main properties of $G$-differential algebras; see~\cite[\S2-\S4]{GuilleminSternbergBook}
and~\cite[\S6]{MeinrenkenBook} for more details.

Let $\fg$ be a~Lie algebra and~$\xi$ be a~generator of the Grassmann algebra $\twedge[\xi]$, i.e., $\xi^2 = 0$.
Let $\dd$ denote the differential on~$\twedge[\xi]$ given by $\dd = \frac{\partial}{\partial\xi} \in \Der\twedge[\xi]$.
Consider a~semisimple $\mathbb{Z}$-graded Lie superalgebra 
\[
  \what{\fg} = \twedge[\xi] \otimes \fg \oplus \mathbb{C} \dd
  = \what\fg_{-1} \oplus \what\fg_0 \oplus \what\fg_{1},
\]
where $\what{\fg}_{-1} = \vspan(\iota_x\mid x\in\fg)$, $\what{\fg}_0 = \vspan(L_x \mid x\in\fg)$, and $\what{\fg}_1 = \vspan(\dd)$.
In the following we refer to $\iota_x$ as contractions and $L_x$ as Lie derivatives.
The explicit brackets in~$\what{\fg}$ are as follows
\begin{align*}
  [L_{x}, L_{y}] &= L_{[x,y]},&
  [L_{x}, \iota_{y}] &= \iota_{[x, y]}, &                                
  [\iota_{x},\dd] &= L_{x},\\
  [\dd,\dd] &= 0, &
  [L_{x}, \dd] &= 0, &
  [\iota_{x},\iota_{y}] &= 0,
\end{align*}
for all $x, y \in \fg$. The equality $[\iota_x,\dd] = L_x$ is called Cartan's magic formula.

A~(graded, filtered) \emph{$\fg$-differential space} is a~(graded, filtered) super vector space~$V$, together
with a~structure of $\what\fg$-representation $\rho\colon\what\fg\to\fgl(V)$
(compatible with gradation/filtration).
A~$\fg$-differential algebra is
a~super-algebra~$\cA$, together with a~structure of a~$\fg$-differential space, such that $\rho$ takes values in~$\Der(\cA)$.

Similarly, let $G$ be a~compact Lie group and denote by~$\fg$ its Lie algebra.
A~\emph{$G$-differential algebra} is a~superalgebra $\mathcal{A}$ with a~representation~$\rho$ of~$G$
by automorphisms of~$\mathcal{A}$ and a~structure of $\fg$-differential algebra on~$\cA$ which satisfies
\begin{align*}
  \frac{d}{dt} \rho (\exp t x) \Big|_{t=0} &= L_{x}, &
  \rho(g) L_{x} \rho(g^{-1}) &= L_{\Ad_g x}, \\
  \rho(g) \iota_{x} \rho(g^{-1}) &= \iota_{\Ad_g x}, &
  \rho(g) \dd \rho(g^{-1}) &= \dd,
\end{align*}
for all $g \in G$ and $x \in \fg$, where $\Ad$ denotes the adjoint action of~$G$ on~$\fg$.

Let $\mathcal{A}$ be a~$G$-differential algebra.
The $G$-horizontal subalgebra of~$\mathcal{A}$, denoted by~$\mathcal{A}_{G-\hor}$, is defined by
\[
  \cA_{G-\hor} = \left\{ a\in\cA \mid \iota_x a = 0\quad\text{for all $x\in\fg$}\right\}.
\]
This could also be denoted by $\cA_{\fg-\hor}$.
Let us emphasise that $\cA_{G-\hor}$ is not a differential subalgebra.

The $G$-basic subalgebra of~$\mathcal{A}$, denoted
by~$\mathcal{A}_{G-\bas}$ is the algebra of $G$-invariants in~$\mathcal{A}$ which are horizontal,
i.e. annihilated by all the contractions by the elements of $\fg$. In other words:
\[
  \mathcal{A}_{G-\bas} =  \mathcal{A}^{G}\cap \mathcal{A}_{G-\hor}.
\]
Similarly, the $\fg$-basic subalgebra of~$\cA$ is defined by
\[
  \mathcal{A}_{\fg-\bas} =  \mathcal{A}^{\fg}\cap \mathcal{A}_{\fg-\hor}.
\]
Note that $\cA_{G-\bas} \subseteq \cA_{\fg-\bas}$.

A \emph{connection} on a~(graded, filtered) $\fg$-differential algebra~$\cA$ 
is an~odd linear map $\vartheta\colon \fg^\ast \to \cA$ (of degree~1)
such that $L_x\vartheta(\mu) = \vartheta(\ad^*_x\mu)$ 
and $\iota_x\vartheta(\mu) = \langle \mu, x \rangle$ for all
$\mu\in\fg^\ast$ and $x\in\fg$. 
A~$\fg$-differential algebra admitting a~connection is called \emph{locally free}.
If $e_i$ is a~basis of $\fg$ then a~connection can be represented by the element $\vartheta = \sum_{i}^{(\fg)}
\vartheta^i\otimes e_i \in \cA \otimes \fg$.

The \emph{curvature of a~connection}~$\vartheta$ is an~even map $F^\vartheta\colon \fg^\ast \to \cA$ of degree~2, defined by
\begin{equation}\label{eq:Curvature}
  F^\vartheta  = \dd\vartheta + \frac12 [\vartheta,\vartheta].
\end{equation}
The curvature map can be represented by the element
$F^\vartheta = \sum_{i}^{(\fg)} (F^\vartheta)^i\otimes e_i \in \cA\otimes\fg$, where
\[
  (F^\vartheta)^i = \dd\vartheta^i + \frac12 \sum_{a,b}^{(\fg)}c^i_{a,b}\vartheta^a\vartheta^b,
\]
where $c^i_{a,b}$ are the structure constants of~$\fg$ in the basis~$e_a$,
i.e., $[e_a,e_b]=\sum_i c^i_{a,b} e_i$.
It is easy to see that the curvature map is $\fg$-equivariant and takes values in~$\cA_{\fg-\hor}$.

\begin{ex}\label{ex:dWedge}
  Consider the case of $\cA=\twedge\fg^\ast$.
  Define contractions~$\iota_x$~of $\twedge\fg^\ast$ by an~element $x\in\fg$ by
Definition~\ref{def:contractions}.
  Let $L_x$ denote the extension of the coadjoint action~$\ad_x^\ast$ to~$\twedge\fg^\ast$ for $x\in\fg$.
  Let $\dd_\wedge\colon \fg^\ast\to\twedge^2\fg^\ast$ be the dual map to the bracket $[-,-]\colon \twedge^2\fg\to\fg$
  extended to a~degree~1 odd derivation of~$\twedge\fg^\ast$; explicitly
  \[
    \dd_\wedge = \frac12 \sum_{a}^{(\fg)} f_a\circ L_{e_a},
  \]
  where $e_a$ is a~basis of~$\fg$ and $f_a$ is the corresponding dual basis of $\fg^\ast$ considered as the
  left multiplication operator.
  In a~basis $f_i$ this takes the following form
  \[
    \dd_\wedge(f_i) = - \frac12\sum_{a,b}^{(\fg)} c^i_{a,b} f_a\wedge f_b.
  \]  
  It is easy to see that $\twedge\fg^\ast$ is a~$\fg$-differential algebra.
  The $\fg$-differential algebra~$\twedge\fg^\ast$ is locally free with connection~$\vartheta_{\wedge}$
  given by~$\vartheta_{\wedge}(\mu) = \mu$ for all $\mu\in\fg^\ast$. It is easy to see that the curvature of
  the connection~$\vartheta_{\wedge}$ is zero. See~\cite[\S6.7]{MeinrenkenBook} for further details.
\end{ex}

\begin{ex}\label{ex:ClisGDiff}
  Assume that $\fg$ admits a non-degenerate symmetric invariant bilinear form~$B$ so we can identify $\fg^\ast \cong \fg$
  via the map $B^\sharp \colon \fg\to\fg^\ast$ and its inverse $(B^{\sharp})^{-1}$.
  Consider the corresponding Cartan 3-tensor
  \begin{equation}\label{eq:Cartan3tensor}
    \phi =  \frac13 \sum_i^{(\fg)} \lambda_\fg(e_i)\wedge f_i \in (\twedge\fg)^\fg.
  \end{equation}
  The Clifford algebra $\Cl(\fg)$ is a~filtered $\fg$-differential algebra with differential given by
  $\dd_{\Cl} = [q(\phi), - ]$, the Lie derivatives are induced by the adjoint action, and the contractions are
  defined in~\S\ref{sec:CliffAlg}.
  Then $\Cl(\fg)$ admits a~flat connection~$\vartheta_{\Cl}$ given by $\vartheta_{\Cl}(\mu) =
  (B^\sharp)^{-1}(\mu)$ for all $\mu\in\fg^\ast$.
\end{ex}

\subsection{The commutative Weil algebra}\label{subsec comm weil}
We recall the definition of the (commutative) Weil algebra of~$\fg$ and its properties
following~\cite[\S6]{MeinrenkenBook}. 
Let
\[
  W(\fg) = S(\fg^\ast)\otimes \twedge \fg^\ast.
\]
For $\mu\in\fg^\ast$, let $\mu$ denote generators of~$1\otimes\twedge\fg^\ast$ and $\what{\mu}$ denote the
generators of~$S(\fg^\ast)\otimes 1$. For $x\in\fg$, let $L_x$ be given by the coadjoint 
action of~$x$ on~$\fg^\ast$ and let $\iota_x$ be defined as
\[
  \iota_x := \id\otimes \iota_x,
\]
where $\iota_x$ in the second factor is as in Example \ref{ex:dWedge}.
Let $\dd_{\mathrm{K}}$ denote the Koszul differential relative to the generators $\mu$ and $\what\mu$, that is,
$\dd_{\mathrm{K}}\mu = \what\mu$, $\dd_{\mathrm{K}}\what\mu = 0$. Let $\dd_{\mathrm{CE}}$ denote the
Chevalley--\/Eilenberg differential
given by 
\[
  \dd_{\mathrm{CE}} = \sum_{a}^{(\fg)} L_{e_a}\otimes f_a + \id \otimes \dd_\wedge,    
\]
where $e_a$ is a~basis of~$\fg$ and $f_a$ is the corresponding dual basis of~$\fg^\ast$ considered as the left
multiplication operator.
Then a~differential on~$W(\fg)$ is given by 
\[
  \dd_W = \dd_{\mathrm{CE}} + 2\dd_{\mathrm{K}}.
\]
The graded $G$-differential algebra~$W(\fg)$ is called the \emph{commutative Weil algebra} of~$\fg$.

\begin{rem}
  The factor 2 in the above definition of~$\dd_W$ could be replaced by any non-zero complex number, in
  particular, this number could be 1. We however use the factor~2 following~\cite{MeinrenkenBook} to ensure
  that the differential~$\dd_W$ is equal to a~Poisson bracket with a~certain element,
  in other words it is a~Hamiltonian operator.
  See~\cite[Remark~6.8 on p.~151 and \S7.2.1]{MeinrenkenBook}.
  Another reason to use the factor~2 is to ensure that the quantisation map~$q_W$ is given by antisymmetrization
  rather than by some more complicated formula,  see Subsection~\ref{ssec:cWg}.
\end{rem}

We have that
\[
  \dd_W(f_i) = 2\what{f}_i + \dd_{\wedge}(f_i) = 2\what{f}_i - \frac12\sum_{a,b}^{(\fg)} c^i_{a,b} f_a\wedge f_b.
\]
If $\fg$ admits a~nondegenerate invariant symmetric bilinear form~$B$, then
\begin{equation}\label{eq:dWmu}
  \dd_W(\mu) = 2\what{\mu} + 2\lambda_\fg(\mu)\qquad\text{for $\mu\in\fg^\ast$.}
\end{equation}

The Weil algebra is locally free, with connection $\vartheta_{W}\colon \fg^\ast \to W(\fg)$ given by
$\mu\mapsto \mu$ for $\mu\in\fg^\ast$. The commutative Weil algebra is universal among (super)commutative
$\fg$-differential algebras: for any commutative (graded, filtered) $\fg$-differential algebra~$\cA$ with
connection $\vartheta_{\cA}$, there is a~unique morphism of (graded, filtered) $\fg$-differential algebras
$c_{\cA}\colon W(\fg) \to \cA$ such that $c_{\cA}\circ \vartheta_W = \vartheta_{\cA}$.
The homomorphism $c_{\cA}$ is called \emph{characteristic homomorphism}.
As it was noticed in Example~\ref{ex:dWedge}, the exterior algebra~$\twedge\fg^\ast$ is (super)commutative
and  locally free. Let $\pi_{\wedge\fg^\ast} \colon W(\fg)\to\twedge\fg^\ast$ denote the corresponding
characteristic homomorphism. On the generators $\mu$ and $\what{\mu}$ it is given by
\[
  \pi_{\wedge\fg^\ast}(\mu) = \mu,\qquad \pi_{\wedge\fg^\ast}(\what{\mu}) = 0.
\]

The Weil algebra~$W(\fg)$ and its $\fg$-invariant part~$W(\fg)^\fg$ are acyclic differential algebras; 
for example, see~\cite[Proposition~6.9 on p.~149]{MeinrenkenBook}.
Therefore, we can define the \emph{transgression} map in the Weil algebra~$W(\fg)$ as follows. Since $\dd_W$
vanishes on~$W(\fg)_{\fg-\bas} = S(\fg^\ast)^\fg$, any element $p\in S^+(\fg^\ast)^\fg$ is a~cocycle, hence
is a coboundary. A~\emph{cochain of transgression} for $p\in S^+(\fg^\ast)^\fg$ is an~odd element
$C_p\in W(\fg)^\fg$ such that $\dd_W(C_p)=p$.
For the symmetric and exterior algebra, using a~$+$ in the superscript means we are talking about elements of
strictly positive degree.

\begin{Def}
For $p\in S^+(\fg^\ast)^\fg$ with a~cochain of transgression~$C_p$, the \emph{transgression} of~$p$ is defined
by~$\mathsf{t}_{\wedge\fg^\ast}(p) = \pi_{\wedge\fg^*}(C_p)$.
\end{Def}

\begin{Def}\label{def:PwedgeG}
Let $P_\wedge(\fg)\subset(\twedge\fg^\ast)^\fg$ be the space of primitive invariants
which can be defined as the~orthogonal complement to $((\twedge^+\fg^\ast)^\fg)^2$  in
$(\twedge^+\fg^\ast)^\fg$.
\end{Def}
In particular this agrees with the definition using the Hopf algebra structure on~$(\twedge\fg^\ast)^\fg$,
see~\cite[Theorem~10.1 and Proposition~10.11]{MeinrenkenBook}.

\begin{thm}[Transgression theorem]\label{thm:transgression}
Let $\fg$ be a~complex reductive Lie algebra. Then
  \begin{enumerate}
  \item For $p\in S^+(\fg^\ast)^\fg$ the transgression map is independent of the choice of cochain.\label{item:tWellDef}
  \item The transgression map~$\mathsf{t}_{\wedge\fg^\ast}$ satisfies\label{item:Tiso}
    \[
      \ker \mathsf{t}_{\wedge\fg^\ast} = (S^+(\fg^\ast)^\fg)^2,\qquad
      \im \mathsf{t}_{\wedge\fg^\ast} = P_\wedge(\fg).
    \]
  \end{enumerate}
\end{thm}

The construction of the transgression map is due to H.~Cartan~\cite{Cartan1951trans}. The transgression theorem was
also explicitly stated by Chevalley in~\cite{ChevalleyICM}; see also~\cite{Leray1951}.
For more details see~\cite[\S10.7]{MeinrenkenBook}.

\subsection{The noncommutative Weil algebra}\label{subsec:noncomm Weil}
Although our main interest in this paper is the study of commutative and quantum Weil algebras and their
relative analogues, we also need the more universal noncommutative Weil algebras.
The main technical reason for introducing them is the existence of certain homotopies which make it possible
to do computations that we can not do directly in the setting of quantum Weil algebras, where such homotopies
are not available.

For a~$\mathbb{Z}$-graded vector space~$E = \bigoplus_kE^k$ and~$n \in \bbZ$ define a graded space with a degree shift~$E[n]$ so that~$(E[n])^k = E^{n+k}$. In other words, if~$x \in E$ is of degree~$k$, then its degree in~$E[n]$ is~$k-n$. In particular, if $E$ is just a vector space we consider it as graded in degree 0, so $E[n]$ is the same space viewed in degree $-n\in\bbZ$.

Following~\cite{AlekseevMeinrenken2005}, define the \emph{noncommutative Weil} algebra of~$\fg$ as the tensor algebra
\[
  \widetilde{W}(\fg) = T(\fg^\ast[-2] \oplus \fg^\ast[-1]).
\]			
For~$\mu \in \fg^\ast$ denote by~$\mu$ the corresponding element of~$\fg^\ast[-1]$ and by~$\bar{\mu}$ the
corresponding element of~$\fg^\ast[-2]$. The algebra $\tW(\frg)$ has a natural $\bbZ$-grading, with generators $\mu$ being of degree 1 and generators $\bar\mu$ of degree 2. In particular, we can view $\tW(\frg)$ as a $\bbZ$-graded superalgebra.

The algebra $\tW(\frg)$ has a natural structure of a graded $\fg$-differential algebra defined as follows; note that by the universal property of the tensor algebra we can define algebra homomorphisms or derivations freely on the generators $\mu$ and $\bar\mu$. 

The differential of $\tW(\frg)$ is the degree 1 derivation $d_{\tW}$ defined on generators by 
\eq
\label{d tW}
\dd_{\tW}\mu = \bar{\mu}; \qquad \dd_{\tW}\bar{\mu}=0.
\eeq
Contraction by an element~$x \in \fg$ is the odd derivation defined on generators by
\eq\label{iota tW}
\iota_x \mu = \langle x, \mu \rangle;\qquad \iota_x \bar{\mu} = \ad_x^\ast\mu.
\eeq
The Lie derivatives $L_x$ for $x\in\fg$ are the even derivations defined on generators by
\eq\label{L tW}
L_x\mu=\ad_x^\ast\mu;\qquad L_x\bar\mu=\ad_x^\ast\bar\mu=\overline{\ad_x^\ast\mu}.
\eeq

The $\fg$-differential algebra $\wtilde{W}(\fg)$ is locally free, with connection $\vartheta_{\tW}\colon \fg^\ast \to \wtilde{W}(\fg)$
given by 
\eq\label{connection tW}
\vartheta_{\tW}(\mu) = \mu,\qquad \mu \in \fg^\ast.
\eeq
Let $e_i$ be a~basis of $\fg$ and $f_i$ be the corresponding dual basis of~$\fg^*\cong\fg^\ast[-1]$.
Since $\vartheta_{\tW}$ is the identity map from $\fg^\ast$ to $\fg^\ast\cong\fg^\ast[-1]\subset
\tW(\frg)$, it can be written as
\[
\vartheta_{\tW} = \sum_{i}^{(\fg)} f_i \otimes e_i\in \tW\otimes\frg
\]
under the identification $\Hom(\frg^*,\tW(\frg))=\tW(\frg)\otimes\frg$.

The curvature $F^{\vartheta_{\tW}} = \sum^{(\fg)}_{i} (F^{\vartheta})^i\otimes e_i \in \widetilde{W}(\fg) \otimes
\fg$ of the connection~$\vartheta_{\tW}$ is given by
\[
  (F^{\vartheta})^i  = \dd_{\tW} f_i + \frac12 \sum_{a,b}^{(\fg)} c^i_{a,b} f_a \otimes f_b =\bar f_i - \delta(f_i)
\]
where  $c^{i}_{a,b}=\langle f_i, [e_a,e_b]\rangle$ are the structure constants of~$\fg$ corresponding to the basis $e_i$, i.e., 
\[
[e_a,e_b]=\sum_i c^i_{a,b} e_i,
\]
and $\delta:\fg^\ast[-1] \to \wtilde{W}(\fg)$ is the linear map defined by
\begin{equation}\label{eq:delta}
  \delta(f_i) = - \frac{1}{2} \sum_{a,b}^{(\fg)} c^i_{a,b} f_a \otimes f_b.
\end{equation}
It is easy to see that $(F^{\vartheta})^i\in \wtilde{W}(\fg)_{\fg-\hor}$.
The map $\delta$ is analogous to $\dd_{\wedge}$ and satisfies (partial) Cartan calculus, see
Lemma~\ref{lem:TCartan} below.
It also satisfies 
\begin{equation}\label{eq:iotaiotadelta}
  \iota_x\iota_y\delta(\mu) = \langle \mu, [x,y] \rangle,\qquad \mu\in\fg^\ast[-1],\ x,y\in\fg,
\end{equation}
i.e., it is dual to the Lie bracket of $\frg$.

We define the hat variables, or curvature variables, by
\[
  \what{f}_i := (F^{\vartheta})^i 
  = \bar{f}_i -\delta(f_i). 
\]
Note that $\what f_i\in\wtilde{W}^2(\fg)_{\fg-\hor}$.

\begin{lem}
  Let~$W'$ be the tensor algebra of the vector space spanned by all~$\mu, \widehat{\mu}$ for $\mu \in \fg^\ast$ and
  define $\chi \colon W' \to \widetilde{W}(\fg)$ by
  $\chi(\mu) = \mu$, $\chi(\widehat{\mu}) = \bar{\mu} - \delta(\mu)$.
  Then~$\chi$ is an algebra isomorphism. 
\end{lem}
\begin{proof}
  By the universal property of tensor algebra~$\chi$ is an algebra morphism. It is obviously injective. We can
  construct the inverse homomorphism by putting $\chi^{-1} (\mu) = \mu$, $\chi^{-1}(\bar{\mu}) = \widehat{\mu} +  \delta(\mu)$
  and extending to $\widetilde{W}(\fg)$ by the universal property.  
\end{proof}

Therefore, $\mu$ and $\what{\mu}$ form another set of generators of~$\wtilde{W}(\fg)$. The main advantage of
this change of generators is that
\[
  \iota_x\mu = \langle \mu, x\rangle,\qquad \iota_x\what{\mu} = 0,\quad\text{for $x\in\fg$}.
\]
\begin{Def}\label{def: delta}
We extend $\delta$ of~\eqref{eq:delta} to a~degree 1 odd derivation of~$\tW(\fg)$ by setting 
$\delta(\what{f}_i)=0$. Explicitly, we have that  
\[
   \delta(f_i) = - \frac12\sum_{a,b}^{(\fg)} c^i_{a,b} f_a\otimes f_b,\qquad
  \delta(\what{f}_i) = 0.
\]
\end{Def}

Consider a $\mathbb{Z}^2$-grading on~$\tW(\fg)$
\begin{equation}\label{eq:tWZ2grading}
  \tW(\fg) = \bigoplus_{i,j\in \mathbb{N}_0} \tW^{i,j}(\fg)
\end{equation}
given by $\deg(\what{\mu})=(1,0)$ and $\deg(\mu)=(0,1)$.
The natural $\bbZ$-grading and this $\bbZ^2$-grading are related by $\tW^{i,j}(\fg) \subseteq \tW^{2i+j}(\fg)$.
Let $\tW^{+,0}(\fg)$ denote the augmentation ideal in the subalgebra $\tW^{\bullet,0}(\fg)$ of~$\tW(\fg)$ generated by~$\what{\mu}$.

The following lemma shows that $\delta$, $\iota_x$ and~$L_x$ define an~analogue of Cartan's calculus on~$\tW^{0,\bullet}(\fg)$. It will be useful in proving the main result of this subsection, Theorem \ref{thm:iotaTtilde}, but it can also serve as motivation for our definition of $\delta$.

\begin{lem}\label{lem:TCartan}
  For $x\in\fg$ we have that on $\tW^{0,\bullet}(\fg)$
  \[
   [\delta,\iota_x] = L_x,
  \]
  and on~$\tW(\fg)$ we have that
    \[
    [\delta,L_x]=0,\qquad [\iota_x,\delta^2] = 0.
  \]

\end{lem}
\begin{proof}
  1) Let $e_i\in\fg$ and $f_i\in\fg^*$ be dual bases. 
  We have
  \[
    [\iota_{e_k},\delta](f_i)
    = \iota_{e_k}\left(-\frac12 \sum_{a,b}c^i_{a,b}f_a\otimes f_b\right) 
    = - \sum_a c^i_{k,a}f_a = L_{e_k}f_i.
  \]
  The last equality follows from
  \[
    \langle L_{e_k}f_i, e_r \rangle =-\langle f_i, [e_k,e_r]\rangle =-\sum_j\langle f_i,c_{k,r}^j e_j\rangle
    = -c_{k,r}^i= - \left\langle \sum_a c^i_{k,a}f_a, e_r\right\rangle.
  \]
  The claim follows.

  2) Since for $\mu\in\fg^\ast$ the element~$\delta(\mu)$ is an anti-symmetric tensor and the action of~$\fg$ preserves anti-symmetric tensors,
  it is enough show that 
  $\iota_x\circ\iota_y\circ L_z\circ \delta (\mu) = \iota_x\circ\iota_y\circ\delta\circ L_z(\mu)$
  for $x,y,z\in\fg$.
  We have
  \begin{align*}
    \iota_x\circ\iota_y\circ L_z\circ\delta(\mu)
    = {}& \iota_x\circ(L_z\circ\iota_y - \iota_{[z,y]})\circ\delta(\mu)\\
    {}={}&(L_z\circ\iota_x\circ\iota_y - \iota_{[z,x]}\circ\iota_y - \iota_x\circ\iota_{[z,y]})\circ\delta(\mu)\\
    {}={}&L_z\circ\iota_x\circ\iota_y\circ\delta(\mu) - \iota_{[z,x]}\circ\iota_y\circ\delta(\mu)
           - \iota_x\circ\iota_{[z,y]}\circ\delta(\mu).\\
    \intertext{Since $\iota_x\circ\iota_y\circ\delta(\mu)\in \mathbb{C}$, the first term is zero, and using~\eqref{eq:iotaiotadelta}
    we get }
    {}={}& - \langle \mu, [[z,x],y] + [x,[z,y]] \rangle = - \langle \mu, [z,[x,y]] \rangle\\
    {}={}& \langle L_z\mu, [x,y] \rangle =  \iota_x\circ\iota_y\circ\delta\circ L_z(\mu).
  \end{align*}
  For $\what{\mu}$ we have that
  \[
    [\delta,L_x](\what{\mu})  = \delta(L_x(\what{\mu})) - 0 = 0.
  \]

  3) By 1) we have that $[\delta,\iota_x] = L_x$, thus by the Leibniz rule 
  \[
    [\iota_x,\delta^2] 
    =[\iota_x,\delta]\circ\delta - \delta\circ[\iota_x,\delta] = [L_x,\delta]
  \]
  which equals $0$ by 2).
\end{proof}

As shown in~\cite[\S3.3]{AlekseevMeinrenken2005}, the noncommutative Weil algebra has  the following universal
property. Given a (graded, filtered) $\fg$-differential algebra~$\cA$ with a connection~$\vartheta_{\cA}$ (so that $\cA$ is locally free), 
there is
a~unique morphism of (graded, filtered) $\fg$-differential algebras 
\eq\label {pi A}
\wtilde{\pi}_{\cA}\colon \tW(\fg)\to \cA,\quad\text{such that}\quad \wtilde{\pi}_{\cA}\circ\vartheta_{\tW(\fg)} = \vartheta_{\cA}.
\eeq
The morphism $\wtilde{\pi}_{\cA}$ is called the \emph{characteristic homomorphism}.

Explicitly, if $\vartheta_{\cA} = \sum_{i}^{(\fg)} \vartheta_{\cA}^i\otimes e_i \in \cA\otimes\fg$, then
\eq\label {pi A expl}
  \wtilde{\pi}_{\cA}(f_i)  = \vartheta_{\cA}^i,\quad
  \wtilde{\pi}_{\cA}(\bar{f}_i)  = \wtilde{\pi}_{\cA}(\dd_{\tW}f_i) = \dd_{\cA}\vartheta_{\cA}^i,\quad
  \wtilde{\pi}_{\cA}(\what{f}_i) = \dd_{\cA}\vartheta_{\cA}^i + \frac12 \sum_{a,b}^{(\fg)}c^i_{a,b} \vartheta_{\cA}^a\vartheta_{\cA}^b.
\eeq

We will be especially interested in the case when the connection $\vartheta_{\cA}$ of the $\fg$-differential
algebra $\cA$ is flat, i.e., its curvature $F^{\vartheta_{\cA}}$ defined by~\eqref{eq:Curvature} vanishes.
Then the corresponding characteristic homomorphism $\wtilde{\pi}_{\cA}: \widetilde{W}(\fg) \to \cA$ has the following properties:

\begin{lem}\label{lem:Aflat}
  Let $\cA$ be a~locally free algebra with a flat connection and $\wtilde{\pi}_{\cA}\colon\tW(\fg)\to \cA$ be the corresponding characteristic
  homomorphism. Then
  \begin{enumerate}
  \item $\wtilde{\pi}_{\cA}(\what{\mu})=0$ for $\mu\in\fg^\ast$,\label{item:Aflat1}  
  \item $\wtilde{\pi}_{\cA}\circ \delta^2 = 0$.\label{item:Aflat3}
  \end{enumerate}
\end{lem}
\begin{proof} 
  (\ref{item:Aflat1}) Follows from the definition of $\what{\mu}$ as curvature components and the fact that
  the connection on~$\cA$ is flat. 
  Namely, since $\wtilde{\pi}_{\cA}$ is the characteristic homomorphism, it is a~morphism of $\fg$-differential
  algebras that agrees with connections: $\wtilde{\pi}_{\cA}(\vartheta_{\tW}^i )= \vartheta_{\cA}^i$.
  Hence, the curvature of the connection on~$\tW(\fg)$ maps to the
  curvature of the connection on~$\cA$:
  \begin{align*}
    \wtilde{\pi}_{\cA}( (F^{\vartheta_{\tW}})^i)
    = {}&\wtilde{\pi}_{\cA}\left( \dd_{\tW} \vartheta_{\tW}^i + \frac12 \sum_{a,b} c^i_{a,b} \vartheta_{\tW}^a\otimes\vartheta_{\tW}^b \right)\\
    {} ={}& \dd_{\cA}(\wtilde{\pi}_{\cA}(\vartheta_{\tW}^i)) + \frac12 \sum_{a,b} c^i_{a,b} \wtilde{\pi}_{\cA}(\vartheta_{\tW}^a)\wtilde{\pi}_{\cA}(\vartheta_{\tW}^b)\\
   {} ={}& \dd_{\cA}\vartheta_{\cA}^i + \frac12\sum_{a,b} c^i_{a,b} \vartheta_{\cA}^a\vartheta_{\cA}^b
    = (F^{\vartheta_{\cA}})^i.
  \end{align*}    

  (\ref{item:Aflat3}) Since $\delta$ is an odd derivation and the map $\wtilde{\pi}_\cA$ is an algebra morphism, then $\wtilde{\pi}_\cA \circ \delta$ satisfies the property, for $a,b \in \tW(\frg)$, 
  \[
    \wtilde{\pi}_\cA \circ \delta(a \otimes b) = (\wtilde{\pi}_\cA\circ \delta(a))( \wtilde{\pi}_\cA(b)) + (-1)^{|a|} (\wtilde{\pi}_\cA(a))( \wtilde{\pi}_\cA\circ \delta(b)).
  \]

  In other words, $\wtilde{\pi}_\cA \circ \delta$ is a derivation of $\tW(\frg)$ with values in $\cA$.
  It is simple to check that $\dd_\cA \circ \wtilde{\pi}_\cA$  is also a derivation with values in $\cA$. We claim $\wtilde{\pi}_\cA \circ \delta = \dd_\cA \circ \wtilde{\pi}_\cA$, it is sufficient to prove this on generators.
  Note that both are clearly zero on $\what{\mu}$. For~$\mu$
  \[
    0 = \wtilde{\pi}_{\cA}(\what{\mu})  = \wtilde{\pi}_{\cA}(\bar\mu-\delta(\mu))=
    \wtilde{\pi}_{\cA}(\dd_{\tW}\mu - \delta(\mu))  = \dd_{\cA}\circ\wtilde{\pi}_{\cA}(\mu) - \wtilde{\pi}_{\cA} \circ\delta(\mu).
  \]
  We conclude $ \wtilde{\pi}_\cA \circ \delta=\dd_\cA \circ \wtilde{\pi}_\cA$.
  Using this identity twice, and $\dd_\cA^2=0$, we find 
  \[
  \wtilde{\pi}_{\cA}\circ\delta^2  = \dd_{\cA}\circ\wtilde{\pi}_{\cA}\circ\delta = \dd_{\cA}^2\circ\wtilde{\pi}_{\cA} = 0.
  \]
  This completes the proof.
\end{proof}

We would like to define a transgression map for $\tW(\frg)$, taking values in an arbitrary $\frg$-differential algebra $\cA$ with a flat connection. To prove the results we require, in particular Theorem \ref{thm:iotaTtilde}, we will need to restrict to a certain subspace that we call the transgression space, defined below.

To construct the transgression map we need a homotopy operator. To define this homotopy operator, we first consider $\mathbb{C}[t,dt]$,  the graded commutative differential algebra with an even generator~$t$ of degree~0 and 
an odd generator~$dt$ of degree~1, and with a single relation $(dt)^2=0$. The differential on
$\mathbb{C}[t,dt]$ is defined by $\dd(t) = dt$, $\dd(dt) = 0$.
One can think of $\mathbb{C}[t,dt]$ as an~algebraic counterpart to differential forms on the unit interval. 

Define an ``integration operator'' $J\colon \mathbb{C}[t,dt] \to \mathbb{C}$ by
\begin{equation}\label{eq:J}
  J\left( \sum_{k\geq0} a_kt^k + \sum_{l\geq0} b_lt^ldt\right)  = \sum_{l\geq0} \frac{b_l}{l+1} \qquad
  \text{for $a_k,b_l\in\mathbb{C}$}.
\end{equation}
In other words, if $P$ and $Q$ are polynomials in $t$, then
\[
J(P+Qdt)=\int_0^1 Qdt.
\]
Let $\eps\colon\tW(\fg)\to\mathbb{C}$ be the augmentation map and $i\colon\mathbb{C}\to\tW(\fg)$ be the unit map.
The standard homotopy operator~$h$ between~$\id$ and $i \circ \eps$ on~$\tW(\fg)$ is defined by
\begin{equation}\label{eq:h}
  h = (J\otimes\id)\circ \phi_h,
\end{equation}
where~$\phi_h$ is the differential algebra morphism from
$\tW(\fg)$ to $\bbC[dt,t] \otimes  \tW(\frg)$ defined on generators by 
\eq\label{def phih}
\phi_h(\mu) = t \mu,\qquad
\phi_h(\bar{\mu}) = dt \mu + t \bar{\mu}.
\eeq

\begin{lem}\label{lem htpy}
Let $h$ be the homotopy map defined above, and let $i$ respectively $\eps$ be the unit and the counit map of $\tW(\frg)$. Then 
\[
h\circ\dd_{\tW}+\dd_{\tW}\circ h=\id-i\circ\eps.
\]
 Moreover, $h$ commutes with $L_x$ for all $x\in\fg$ and $h$ preserves $\tW(\frg)^\frg$.
\end{lem}
\pf This is analogous to the commutative case, see \cite[Section 6.4]{MeinrenkenBook}. 
\epf

To define the transgression space, we first recall the symmetrization map
\[
  \sym_W\colon W(\fg)\to\tW(\fg)
\]
described in \cite[Proposition~6.13]{MeinrenkenBook}. This map is obtained by symmetrizing the (even) variables $\bar\mu$ and skew-symmetrizing the (odd) variables $\mu$. Explicitly, $\sym_W$ is the linear map defined on monomials $v_1\ldots v_k$ in $W(\frg)$ such that each of $v_1,\dots,v_k$ is either some $\mu$ or some $\bar\mu$, by
\[
  \sym_W(v_1\ldots v_k) =  \frac{1}{k!}\sum_{\sigma\in\mathbf{S}_k}(-1)^{N_\sigma(v_1,\ldots,v_k)} \phi(v_{\sigma^{-1}(1)})\ldots\phi(v_{\sigma^{-1}(k)}),
\]
where $\mathbf{S}_k$ is the symmetric group, and $N_\sigma(v_1,\ldots,v_k)$ is the number of pairs $i<j$ such that
$v_i$, $v_j$ are odd elements and 
$\sigma^{-1}(i) > \sigma^{-1}(j)$.

It is proved in \cite[Proposition~6.13]{MeinrenkenBook} that $\sym_W$ is an~injective homomorphism of $\fg$-differential
spaces. (It is not an algebra homomorphism, but it is an algebra morphism up to homotopy; see \cite[Corollary~3.4]{AlekseevMeinrenken2005}.)
From this and the above definition, it is clear that
\begin{align*}
  & \sym_W(f_i) =  f_i,\qquad \sym_W(\dd_Wf_i) = \dd_{\tW}f_i=\bar f_i, \\
  & \sym_W(\dd_\wedge(f_i)) = \delta(f_i).
\end{align*}
From~\eqref{eq:dWmu}, in~$W(\fg)$ we have that
\[
  \what{f}_i = \frac12\dd_{W}(\mu) - \lambda_\fg(f_i) =  \frac12(\dd_{W}\mu - \dd_\wedge f_i)
\]
and in~$\tW(\fg)$ we have that
\[
  \what{f}_i =  \dd_{\tW}f_i  -\delta(f_i).
\]
Since $\sym_W$ is a~morphism of $\fg$-differential spaces we get that
\[
  \sym_W(\what{f}_i) = \frac12 \sym_W(\dd_Wf_i - \dd_\wedge(f_i)) = \frac12(\dd_{\tW}f_i-\delta(f_i))
  = \frac12\what{f}_i.
\]
The factor $1/2$ appears due to difference in the definition of differentials in~$W(\fg)$ and~$\tW(\fg)$, for
more details see~\cite[Remark~6.8 and \S7.2.1]{MeinrenkenBook}.

If $p\in S(\fg^\ast)^\fg\subset W(\fg)$, then $\dd_{W}p = 0$, hence $\dd_{\tW}(\sym_W(p))=0$.
Moreover, since $\tW(\fg)$ is acyclic by Lemma~\ref{lem htpy}. (Here acyclicity means that $H^0(\tW(\fg)) = \bbC$
and $H^i(\tW(\fg)) = 0$ for $i>0$.)
We also have that~$\sym_W(p)\in\im\dd_{\tW}$ if $p\in S^+(\fg^\ast)^\fg$.
Recall, the plus in the superscript signifies strictly positive degree.

\begin{Def}\label{def:trangr}
 The \emph{transgression space} for~$\tW(\fg)$ is
  \[
    \wtilde{\mathfrak{T}}(\fg) = \{ \what{\sym}(p)  \mid p\in S^+(\fg^\ast)^\fg\}
    \subset \tW^{+,0}(\fg)^\fg,
  \]
  where $\what{\sym}\colon S(\fg^\ast) \to \tW^{\bullet,0}(\fg)$ is symmetrization in the hat
  variables (which is clearly $\fg$-equivariant).
\end{Def}
We note that $\wtilde{\mathfrak{T}}(\fg)$ is a~graded subspace of~$\tW(\fg)$.
We now define the transgression map for the non-commutative Weil algebra as follows.

\begin{Def}\label{def trngr tW}
Let $\cA$ be a locally free $\frg$-differential algebra with a flat connection and let $\wtilde{\pi}_\cA\colon \tW(\fg)\to \cA$ be the corresponding characteristic homomorphism defined in \eqref{pi A} and~$h$ the homotopy operator defined in~\eqref{eq:h}. 
The (universal)
\emph{transgression} map of $\tW(\frg)$ with values in~$\cA^\frg$ is the linear map
\[
  \wtilde{\mathsf{t}}_{\cA} = \wtilde{\pi}_\cA \circ h\colon \wtilde{\mathfrak{T}}(\fg)\to \cA^\fg.
\]
(The image of $\wtilde{\mathsf{t}}_\cA$ is contained in $\cA^\frg$, since $\wtilde{\pi}_\cA$ and $h$ commute with $L_x$ and $\wtilde{\mathfrak{T}}(\fg) \subset \tW(\frg)^\frg$.)

\end{Def}

\begin{rem}
   We could not formulate a set of conditions to define the transgression map for~$\tW(\fg)$ with values
  in~$\cA^\fg$ using cochains of transgression.
  Instead we use the homotopy~$h$.
\end{rem}

\subsection{On introduced derivations of~$\tW(\fg)$}
We have already defined multiple derivations throughout the next two sections we will introduce many more derivations of $\tW(\frg)$. For reference we include the following table of derivations as well as several Lemmas describing their properties and relations between them.

For each derivation $D$ of~$\tW(\fg)$ we list its parity $p(D)\in\bbZ/2\bbZ = \{\bar0, \bar1\}$ and the action of~$D$ on the generators $\mu$,
$\what{\mu}$, $\bar{\mu} = \what{\mu} + \delta(\mu)$.
\begin{equation}\label{tableder} 
\begin{minipage}[l]{11.25cm}
\begin{center}
\begin{tabular}{c|c|c|c|c}
  $D$ & $p(D$) & $D(\mu)$ & $D(\what{\mu})$ & $D(\bar{\mu})$ \\
  \hline 
  $\delta$ & $\bar1$ & $-\frac{1}{2}\sum\limits_{a,b}\langle \mu , [e_a, e_b] \rangle f_a \otimes f_b$  & $0$ & $\delta^2(\mu)$ \\
  $\dd_{\tW}$ & $\bar1$ & $\bar{\mu}$ & $-\dd_{\tW} \delta(\mu)$ & $0$ \\
  $\dd_{\what{K}}$& $\bar1$ & $\what{\mu}$ &  $0$ & $\dd_{\what{K}}\delta(\mu)$ \\
  $L_x$ & $\bar0$ & $\ad^\ast_x \mu$ & $\what{\ad^\ast_x \mu}$ & $ \overline{\ad^*_x\mu} $ \\
  $\iota_x$ & $\bar1$ & $\langle x, \mu\rangle$ & 0 & $\ad^\ast_x\mu$ \\ 
  $\iota_{\what{x}}$ & $\bar0$ & 0 & $\langle x, \mu\rangle$  & $\langle x, \mu\rangle$  \\
  $g$ &  $\bar1$ &  $0$ & $\mu$ & $\mu$ \\
  $\eta$ & $\bar0$ & $0$ & $\delta(\mu)$ & $\delta(\mu)$ \\
  $\zeta = [\eta, \dd_{\what{K}}]$  & $\bar1$ & $ \delta(\mu)$ & $-\dd_{\what{K}}\delta(\mu)$ & $-\dd_{\what{K}}\delta(\mu) +\zeta\delta(\mu)$ \\
\end{tabular}
\end{center}
\end{minipage}
\end{equation}
\bigskip

Recall that the super commutator of any two derivations is a derivation, in particular, $\delta^2$ is a derivation.

As remarked after \eqref{eq:delta} the derivation $\delta \in \Der
\tW(\frg)$ is the counterpart to $\dd_\wedge \in \Der W(\frg)$, although $\delta$ fails to be a differential. If one replaces $\delta$ by $\dd_\wedge$ and the algebra $\tW(\frg)$ by $W(\frg)$ in all of the definitions in Table \eqref{tableder} one recovers counterparts of these derivations for the commutative Weil algebra $W(\frg)$.
The counterpart of the derivation~$\zeta \in\Der\tW(\fg)$
is the Chevalley--Eilenberg differential~$\dd_{\mathrm{CE}}\in\Der W(\fg)$.
The counterpart of~$g\in\Der\tW(\fg)$ is the de Rham differential of~$W(\fg)$.

The following two lemmas give some relations between various derivations of~$\tW(\fg)$ that we have
defined. They will also be useful in the proof of Theorem~\ref{thm:iotaTtilde}.

Let $\varphi:\tW(\frg)\to\tW(\frg)$ be the linear map acting by 0 on the constants, and by $\frac{1}{i+j}$ on
the subspace $\tW^{i,j}(\frg)\subset\tW(\frg)$ consisting of elements of degree $i$ with respect to the hat
variables and of degree $j$ with respect to the $\mu$ variables, with $i+j>0$.

\begin{Def}
Let $I_{\delta^2}$ be the two-sided ideal of~$\tW(\fg)$ generated by the image of~$\delta^2$.
\end{Def}
\begin{lem}\label{lem:derIdelta2}
  The derivations of~$\tW(\fg)$,  $\iota_x$, $\iota_{\what{x}}$, $\eta$ and $g$ preserve $I_{\delta^2}$.
  Moreover,
  we have the following identities
  \begin{subequations}
  \begin{gather}
    [\iota_x ,\delta^2] = 0,\tag{\ref{lem:derIdelta2}.a}\\
    [\iota_{\what{x}}, \delta^2]=0,\tag{\ref{lem:derIdelta2}.b}\\
    [g,\delta^2](\mu) = 0,\quad
    [g,\delta^2](\what{\mu}) = - \delta^2(\mu)\tag{\ref{lem:derIdelta2}.c},\\
    [\eta,\delta^2](\mu) = 0,\quad
    [\eta,\delta^2](\what{\mu}) = - \delta^3(\mu).\tag{\ref{lem:derIdelta2}.d}
  \end{gather}
\end{subequations}
Also the linear maps~$\varphi$ and~$h$ preserve $I_{\delta^2}$.
\end{lem}
\begin{proof}
  First we will prove the identities.
  
  \noindent (a) This identity is proved in Lemma~\ref{lem:TCartan}. 
  
  \noindent (b) For $\mu\in\fg^\ast$, we have that
    \[
    [\iota_{\what{x}},\delta^2](\mu) = \iota_{\what{x}}(\delta^2(\mu)) - \delta^2(0) = 0,\quad 
    [\iota_{\what{x}},\delta^2](\what{\mu}) = \iota_{\what{x}}(0) - \delta^2(\iota_{\what{x}}\what{\mu})= 0.
  \]
  \noindent (c) This is obvious.
  
  \noindent (d) For $\eta$ we have that
  \[
    [\eta,\delta^2](\mu)  = \eta(\delta^2(\mu)) - 0 = 0,\quad
    [\eta,\delta^2](\what{\mu}) = \eta(0) - \delta^2(\delta(\mu))  = - \delta^3(\mu).
  \]
 
  Now we prove that the derivations preserve~$I_{\delta^2}$.
  We note that any derivation which commutes with~$\delta^2$ also preserves $I_{\delta^2}$.
  It follows now from (a)
  and the identities above that  $\iota_x$ and $\iota_{\what{x}}$ preserve~$I_{\delta^2}$

  For $g$, we have that $g\circ\delta^2 (\mu) = \delta^2\circ g(\mu)$.  

  Let us note that elements of the form $a\otimes\delta^2(\mu)\otimes c$, where
  $a,c\in\tW(\fg)$, $\mu\in\fg^\ast$, span~$I_{\delta^2}$ (recall that $\delta^2(\what{\mu}) = 0$).
  Then
  \begin{equation*}
    g(a\otimes \delta^2(\mu)\otimes c) = g(a)\otimes\delta^2(\mu)\otimes c
    \pm a\otimes(\delta^2\circ g (\mu)) \otimes c 
    \pm a\otimes\delta^2(\mu)\otimes g(c)
  \end{equation*}
  which belongs to~$I_{\delta^2}$. Thus, $g$ preserves~$I_{\delta^2}$. 
  
  Using a~similar argument with (d), we find that $\eta$ preserves $I_{\delta^2}$.  
  The same applies to $h$, because $h = (J\otimes\id)\circ\phi_h$
  and $\phi_h(\delta^2(\mu)) = t^3\delta^2(\mu)$.
  Since $I_{\delta^2}$ is a graded subspace, it is also preserved by $\varphi$.
\end{proof}

Let $\dd_{\what{K}} \in
\Der^1\tW(\fg)$ be the Koszul differential for $\mu$ and $\what{\mu}$ generators defined by (see Table~\ref{tableder})
\eq\label{def dKhat}
\dd_{\what{K}}(\mu) = \what{\mu},\qquad \dd_{\what{K}}(\what{\mu})=0.
\eeq
Recall from the Table~\ref{tableder}, the odd degree~1 derivation  $\zeta := [\eta,\dd_{\what{K}}]\in\Der^1\tW(\fg)$. We note that $\zeta$ is
an analogue of the Chevalley--Eilenberg differential $\dd_{\mathrm{CE}}\in\Der W(\fg)$ for the noncommutative
Weil algebra~$\tW(\fg)$. However, $\zeta^2\neq 0$. The following lemma describes the interaction of
the derivations $\zeta$, $\eta$ and $\delta^2$. 
\begin{lem}\label{lem:M}  
  For $X\in\tW(\fg)$ and $p\in\wtilde{\frT}(\fg)$ we have
  \begin{enumerate}
  \item $[\eta,\zeta](X)  \in I_{\delta^2}$,
  \item $\dd_{\tW} = \dd_{\what{K}} + \zeta - [\delta^2,g]$, hence
    $(\dd_{\tW} - \dd_{\what{K}} - \zeta)(X) \in I_{\delta^2}$.
  \item $\dd_{\tW}(p)\in I_{\delta^2}$ and $\zeta\circ\eta^k(p) \in I_{\delta^2}$, for all $k\geq 0$.
  \end{enumerate}
\end{lem}

\begin{proof}
(a)  For $\mu\in\fg^\ast$ we have
  \[
    \zeta(\mu) = [\eta,\dd_{\what{K}}](\mu)= \eta(\what\mu)-0= \delta(\mu),
 \]
 which implies
 \[
    [\eta,\zeta](\mu) = \eta\circ\delta(\mu) - \zeta(0) = 0.
 \]
Furthermore, we have 
\begin{equation}\label{eq:MmuHat}
  \zeta(\what{\mu}) = [\eta,\dd_{\what{K}}](\what\mu)=0-\dd_{\what{K}}\circ\eta(\what\mu)=
  - \dd_{\what{K}}\circ\delta(\mu),
\end{equation}
which implies
\[
  [\eta,\zeta](\what{\mu}) = \eta\circ \zeta(\what\mu) - \zeta\circ\eta(\what\mu)=
  -\eta\circ\dd_{\what{K}}\circ\delta(\mu) - \zeta\circ\delta(\mu).
\]
We rewrite this using $\eta\circ d_{\what K} = \zeta - d_{\what K}\circ\eta$ and $\eta\circ\delta(\mu)=0$; we obtain
\[
  [\eta,\zeta](\what{\mu}) = -2\zeta\circ\delta(\mu).
\]
Finally, since~$\zeta$ and~$\delta$ are both odd derivations of $\tW^{0,\bullet}(\fg)$  and since they agree on the generators $\mu$ of 
$\tW^{0,\bullet}(\fg)$, they agree on all of $\tW^{0,\bullet}(\fg)$, which implies $\zeta\circ\delta(\mu)=\delta^2(\mu)$.
So we see that $[\eta,\zeta](\what{\mu}) =-2\delta^2(\mu)\in I_{\delta^2}$.

(b) First we see that
\begin{align*}
  (\dd_{\tW}-\dd_{\what{K}})(f_i) = {} & \what{f}_i + \delta(f_i) - \what{f}_i = \delta(f_i) = \zeta(f_i),
  \end{align*}
  so 
  \eq\label{display 1}
  (\dd_{\tW}-\dd_{\what{K}}-\zeta  + [\delta^2,g])(f_i)=0\in I_{\delta^2}.
  \eeq
  Furthermore, 
  \begin{align*}
  (\dd_{\tW}-\dd_{\what{K}})(\what{f}_i) = {} & \dd_{\tW}(\bar{f}_i - \delta(f_i))  = \dd_{\tW}(- \delta(f_i))\\
  {}={}& - \frac12 \sum_{a,b}c^i_{ab} \left(\dd_{\tW}(f_a)\otimes f_b - f_a\otimes \dd_{\tW}(f_b)\right)\\
  {}={}& - \frac12 \sum_{a,b}c^i_{ab} \left((\what{f}_a + \delta(f_a))\otimes f_b - f_a\otimes(\what{f}_b + \delta(f_b))\right)\\
  \text{(By~\eqref{eq:MmuHat})}{}={}& \zeta(\what{f}_i) - \frac12 \sum_{a,b}c^i_{ab} \left( \delta(f_a)\otimes f_b - f_a\otimes \delta(f_b)\right)\\
    {}={}& \zeta(\what{f}_i) - \delta^2(f_i) = \zeta(\what{f}_i) - [\delta^2,g](\what{f}_i).
\end{align*}
So $(\dd_{\tW}-\dd_{\what{K}}-\zeta)(\what{f_i})=- \delta^2(f_i)\in I_{\delta^2}$, and together with \eqref{display 1} this implies~(b).

(c) We use induction on~$k$.
Assume that $k=0$.
Since $p\in\wtilde{\frT}(\fg)$ we have that $\dd_{\what{K}}(p) =0$.
Using (b), we have that $(\dd_{\tW}-\dd_{\what{K}} - \zeta)(p)\in I_{\delta^2}$,
since $[\delta^2,g](p) = \delta^2\circ g (p) \in I_{\delta^2}$, thus  $\dd_{\tW}(p) \in I_{\delta^2}$.
 It follows that
\[
  \zeta(p) = (\dd_{\tW} - \dd_{\what{K}} + [\delta^2,g])(p)\in I_{\delta^2}.
\]
This is the base of our induction.

Assume now that $\zeta\circ\eta^k(p) \in I_{\delta^2}$. We have that
\[
  \zeta\circ\eta^{k+1}(p) = \eta\circ \zeta\circ\eta^k(p) + [\zeta,\eta](\eta^k(p))
\]
The first term belongs to~$I_{\delta^2}$ by the induction hypothesis, since $\eta$ preserves~$I_{\delta^2}$ by Lemma \ref{lem:derIdelta2}.
The second term belongs to~$I_{\delta^2}$ by~(a). Thus $\zeta \circ \eta^{k+1}(p)$ is in the ideal $I_{\delta^2}$.
\end{proof}

\subsection{Contractions and transgressions}

Let $\Phi\colon\tW(\fg)\to\tW(\fg)$ be the $\fg$-equivariant algebra map defined on generators by
\begin{equation}\label{eq:Phi}
  \Phi(\mu) = \mu,
  \qquad \Phi(\what{\mu}) = \delta(\mu).
\end{equation}

Furthermore, let~$\cA$ be a~locally free $\frg$-differential algebra with a~flat connection and let $\wtilde{\pi}_\cA\colon \tW(\fg)\to \cA$ be the corresponding characteristic homomorphism defined by~\eqref{pi A}.
Define the $\frg$-equivariant algebra homomorphism 
\eq\label{eq:PhiA}
\Phi_{\cA} = \wtilde{\pi}_{\cA} \circ \Phi: \tW(\fg) \to \cA.
\eeq

The main result of this subsection is the following theorem; it is a generalisation of Theorem~73 in~\cite{KostantRho}, and its proof is a generalisation
of~\cite[Proposition~6.19 on p.~160]{MeinrenkenBook}. This result will be used in Subsection~\ref{ssec:cWg},
and  Section \ref{sec:trans in rel weil} to prove analogous (relative) versions,
Theorem~\ref{thm:iotaxImAlphaG} and Theorem~\ref{thm:IotaTClP}, which will in turn be essential for the proof
of our main result, Theorem~\ref{thm:main}.

\begin{thm}\label{thm:iotaTtilde}
Let $\cA$ be a~locally free $\frg$-differential algebra with a~flat connection and let~$\Phi_{\cA}$ be the algebra homomorphism defined by~\eqref{eq:PhiA}. Furthermore, let $\wtilde{\mathsf{t}}_{\cA}$ be the transgression map defined in Definition~\ref{def trngr tW}, and let $\iota_x,\iota_{\widehat{x}}$ be the contractions defined in Table~\ref{tableder}.
    
    If $x \in \fg$ and $p \in \wtilde{\mathfrak{T}}^{m+1}(\fg)$, then
  \[
    \iota_x \wtilde{\mathsf{t}}_{\cA}(p) = \frac{(m!)^2}{(2m)!} \cdot\Phi_{\cA}(\iota_{\widehat{x}}p).
  \]
\end{thm}

\pf

Our first task is to rewrite the left side of the claimed equality using Lemma \ref{lem htpy}. 
Recall from Table \ref{tableder} the odd derivation $g$ of $\tW(\frg)$, defined on generators by
\[
  g(\what\mu)=\mu,\qquad g(\mu)=0.
\]
The following proposition allows us to write $\wtilde{\mathsf{t}}_{\cA}$ as the composition of~$\Phi_{\cA}$
and~$g$ which will make it easier to apply $\iota_x$.
\begin{prop}\label{prop:trformula}  
  For any $p \in \wtilde{\mathfrak{T}}^{m+1}(\fg)$,
\[
\wtilde{\mathsf{t}}_{\cA}(p) = \frac{(m!)^2}{(2m+1)!} \Phi_{\cA} (g(p)).
\]
\end{prop}

\begin{proof}
Since $p$ is in $\wtilde{\mathfrak{T}}^{m+1}(\fg)\subseteq \tW^{+,0}(\fg)^\fg$, we can write it as
\[
p=\sum \what\mu_1\otimes\dots\otimes\what\mu_{m+1}.
\]
  It is easy to see that $\wtilde{\mathsf{t}}_{\cA} = \wtilde{\pi}_{\cA}\circ h = (J \otimes \Id) \circ (\Id \otimes \wtilde{\pi}_{\cA}) \circ \phi_h$.
  Note that
  \[
    \phi_h(\widehat\mu_i)=\phi_h(\bar\mu_i-\Phi(\widehat\mu_i))=dt\mu_i+t\bar\mu_i-t^2\Phi(\widehat\mu_i)=dt\mu_i+t\widehat\mu_i+(t-t^2)\Phi(\widehat\mu_i),
  \]
  using the fact that $\Phi(\widehat\mu_i)=\delta(\mu_i)$, and that it is in  $\tW^{0,2}(\frg)$, so $\phi_h$
  acts on $\Phi(\what\mu_i)$ by multiplying by $t^2$.  
This implies
  \begin{multline*}
    (J \otimes \Id) \circ (\Id \otimes \wtilde{\pi}_{\cA}) \circ \phi_h(p)  = {}\\
    {} = \sum(J \otimes \Id) \circ (\Id \otimes \wtilde{\pi}_{\cA}) \Bigl( ( dt \mu_1  + t \widehat{\mu}_1 + (t - t^2) \Phi(\what{\mu}_1)) \cdot \ldots\\
    \ldots  \cdot (dt \mu_{m+1}  + t \widehat{\mu}_{m+1} + (t - t^2) \Phi(\what{\mu}_{m+1})) \Bigr).
  \end{multline*}
  Applying $\wtilde{\pi}_{\cA}$ 
and remembering that $\wtilde{\pi}_{\cA}(\what\mu_i)=0$ by Lemma \ref{lem:Aflat}, we find 
  \begin{multline*}
    (J \otimes \Id) \circ (\Id \otimes \wtilde{\pi}_{\cA}) \circ \phi_h(p) = {}\\
    {} = \sum \left( J \otimes \Id \right)\Bigl( ( dt \wtilde{\pi}_{\cA}(\mu_1) + (t - t^2) \Phi_{\cA}(\what{\mu}_1))\cdot  \ldots\\
      \ldots \cdot (dt \wtilde{\pi}_{\cA}(\mu_{m+1})  + (t - t^2) \Phi_{\cA}(\what{\mu}_{m+1})) \Bigr).
  \end{multline*}
  which has linear $dt$ term 
  \[
    (J\otimes\id)\sum (t-t^2)^m dt  \sum_{i} \Phi_{\cA}(\widehat{\mu}_1 \otimes \cdots \otimes \widehat{\mu}_{i-1})
    \wtilde{\pi}_{\cA}(\mu_i) \Phi_{\cA}(\widehat{\mu}_{i+1} \otimes \cdots \otimes \widehat{\mu}_{m+1}).
  \]
    The lemma follows by noting that $J \left( (t-t^2)^m dt \right) = \frac{(m!)^2}{(2m+1)!}$.
\end{proof}

Getting back to the proof of Theorem \ref{thm:iotaTtilde}, we now see that the left side of the claimed equality is equal to
\[
\frac{(m!)^2}{(2m+1)!} \iota_x\Phi_\cA(g(p))=\frac{(m!)^2}{(2m+1)!} \wtilde{\pi}_{\cA} (\iota_x\Phi(g(p))).
\]
On the other hand, the right side of the equality we need to prove is
$\frac{(m!)^2}{(2m)!}\Phi_\cA(\iota_{\what x}p),
$
so the equality to prove becomes
\eq\label{toprove}
\wtilde{\pi}_{\cA} (\iota_x\Phi(g(p)))=(2m+1)\Phi_{\cA} (\iota_{\what x}p).
\eeq
Now we rewrite the right side of \eqref{toprove} as follows.
Recall from Table~\eqref{tableder} that the even derivation $\eta: \wtilde{W}(\fg) \to \wtilde{W}(\fg)$ acts
on generators by $\eta(\mu) = 0$ and $\eta(\what{\mu}) = \delta(\mu)$.
We claim that for any $q \in \tW^{m,0}(\fg)$ we have
  \begin{equation}\label{eq:piAq}
    \Phi_{\cA}(q) = \wtilde{\pi}_{\cA} \left(\frac{\eta^m}{m!}(q) \right).
  \end{equation}
In particular, this will hold for $q=\iota_{\what x} p\in \tW^{m,0}(\fg)$.

  To prove \eqref{eq:piAq}, we can assume $q=\what\mu_1\otimes\dots\otimes\what\mu_m$ and write
  \[
    \Phi_\cA(q)=\wtilde{\pi}_{\cA}(\Phi(\what\mu_1)\otimes\dots\otimes\Phi(\what\mu_m))=\wtilde{\pi}_{\cA}(\delta(\mu_1)\otimes\dots\otimes\delta(\mu_m))=\wtilde{\pi}_{\cA}(\eta(\what\mu_1)\otimes\dots\otimes\eta(\what\mu_m)).
  \]
  On the other hand,
  since $\eta^2(\what\mu_i)=0$,
  \[
    \wtilde{\pi}_{\cA} \left(\frac{\eta^m}{m!}(q)\right)= \wtilde{\pi}_{\cA}(\eta(\what\mu_1)\otimes\dots\otimes\eta(\what\mu_m)).
  \]
  This proves~\eqref{eq:piAq}.
In particular, \eqref{eq:piAq} for $q=\iota_{\what x}p$ implies that the right side of \eqref{toprove} is equal to
\[
  \frac{2m+1}{m!} \wtilde{\pi}_{\cA} ( \eta^m(\iota_{\what{x}}p))
\]
Since the derivation $[\iota_{\what{x}},\eta]$ is zero,
which can easily be checked on generators, the last expression is equal to
\eq\label{rhs}
\frac{2m+1}{m!} \wtilde{\pi}_{\cA} (\iota_{\what{x}} \eta^m(p)).
\eeq

To pass between $\iota_{\what{x}}$ and $\iota_x$, we need the following lemma. 

\begin{lem}
\label{final lem}
Let~$p = \sum\what{\mu}_1 \otimes \dots \otimes \what{\mu}_{m+1} \in \wtilde{\mathfrak{T}}^{m+1}(\fg)$.
  Then
  \[
 \wtilde{\pi}_{\cA}(\iota_{\what{x}}\eta^m(p))= \wtilde{\pi}_{\cA}(\iota_x\circ g\circ\varphi(\eta^m(p))).
 \]
\end{lem}

Assuming Lemma \ref{final lem}, we now finish the proof of Theorem \ref{thm:iotaTtilde}. Recall that we need to prove \eqref{toprove}, and that we have rewritten the right side of \eqref{toprove} in \eqref{rhs}. Using Lemma \ref{final lem}, \eqref{rhs} becomes 
\[
\frac{2m+1}{m!}\wtilde{\pi}_{\cA}(\iota_x\circ g\circ\varphi(\eta^m(p))), 
\]
which is equal to $\frac{1}{m!}\wtilde{\pi}_{\cA}(\iota_x\circ g(\eta^m(p)))$ since $\eta^m(p)$ has bidegree $(1,2m)$ so $\varphi$ acts on it by $\frac{1}{2m+1}$.

It now suffices to prove that
\eq\label{toprove 2}
g(\eta^m(p))=m!\Phi(g(p));
\eeq
this will show that the above expression for the right side of~\eqref{toprove} is equal to the left side of~\eqref{toprove}.
The proof of~\eqref{toprove 2} is straightforward, using the fact that~$\eta$ sends the~$\mu$ variables to~$0$. 
This finishes the proof of Theorem~\ref{thm:iotaTtilde} modulo 
Lemma~\ref{final lem}.
\epf

\begin{proof}[Proof of Lemma \ref{final lem}.]
By Lemma \ref{lem:Aflat}(\ref{item:Aflat3}) $\wtilde\pi_{\cA}$ annihilates the ideal $I_{\delta^2}$ of $\tW(\frg)$.
So it suffices to prove
\begin{equation}\label{toprove 3}
  \iota_{\what{x}}\eta^m(p) - (\iota_x\circ g\circ\varphi)(\eta^m(p))\in I_{\delta^2}.
\end{equation}

Recall from Table~\eqref{tableder} that
$g$ sends $\what\mu$ to $\mu$ and $\mu$ to $0$,
while $\\d_{\what K}$ sends $\mu$ to $\what\mu$ and $\what\mu$ to $0$,
so it follows that $[\dd_{\what{K}},g]$ acts as multiplication by~$i+j$ on~$\wtilde{W}^{i,j}(\fg)$.
The map $\varphi$ is the inverse of $[\dd_{\what{K}},g]$ on~$\tW^+(\fg)$. Moreover, $g \circ \varphi$ is
a~homotopy operator for~$\dd_{\what{K}}$, i.e., if $i$ and $\varepsilon$ are the unit and counit maps
of~$\tW(\fg)$ as in Lemma~\ref{lem htpy}. Then
  \eq\label{dhat htpy}
  [\dd_{\what{K}},g \circ \varphi] = \id - i\circ\varepsilon. 
  \eeq
It is clear that
  \[
    \iota_{\what{x}}\circ \eta^m(p)\in \tW^{0,+}(\fg),\quad
    \iota_x\circ g\circ \varphi\circ\eta^m(p) \in \tW^{0,+}(\fg).
  \]  
  Since $\ker\dd_{\what{K}}\cap \tW^{0,+}(\fg) = 0$, equation \eqref{toprove 3} (and hence the lemma) will follow if we show that
\eq\label{toprove 4}
\left(\iota_{\what{x}}- \iota_x\circ g\circ\varphi\right)(\eta^m(p)) \in\ker \dd_{\what{K}} +  I_{\delta^2}.
\eeq
We claim that the operator  $\iota_{\what{x}} - \iota_x\circ g\circ \varphi$ preserves $\ker\dd_{\what{K}}+I_{\delta^2}$.

To see this, we first show that 
\eq\label{comm dkhat iotahat}
[\dd_{\what{K}},\iota_{\what{x}}] = -\iota_x.
\eeq
Indeed, if $\mu \in \fg^*$ and~$x \in \fg$, then
  \begin{align*}
    [\dd_{\what{K}},\iota_{\what{x}}](\mu) ={}&{} - \iota_{\what{x}} \dd_{\what{K}}(\mu) = - \iota_{\what{x}}\what{\mu} =- \langle \mu, x \rangle =- \iota_x \mu,\\
    [\dd_{\what{K}},\iota_{\what{x}}](\what{\mu}) ={}&{} 0 = -\iota_x \what{\mu},
  \end{align*}
  so the derivations $[\dd_{\what{K}},\iota_{\what{x}}]$ and $ -\iota_x$ agree on generators and hence are equal. Using \eqref{comm dkhat iotahat}, \eqref{dhat htpy},  and the obvious fact $[\dd_{\what{K}}, \iota_x] = 0$, we get
  \begin{align*}
    [\dd_{\what{K}},\iota_{\what{x}} - \iota_x \circ g \circ \varphi] &{}= - \iota_x + \iota_x \circ [\dd_{\what{K}},g \circ \varphi] \\
   {} &{} = -\iota_x \circ i\circ \varepsilon = 0.
  \end{align*}
  So $\iota_{\what{x}} - \iota_x\circ g\circ \varphi$ preserves $\ker\dd_{\what{K}}$.
  Lemma~\ref{lem:derIdelta2} shows that $\iota_{\what{x}} - \iota_x\circ g\circ \varphi$ also preserves $I_{\delta^2}$.
  This finishes the proof of the claim that $\iota_{\what{x}}-\iota_x\circ g\circ \varphi$ preserves
  $\ker\dd_{\what{K}}+ I_{\delta^2}$.

We now see that \eqref{toprove 4} (and hence the lemma) will follow if we prove
\eq\label{toprove 5}
\eta^m(p)\in \ker\dd_{\what{K}}+ I_{\delta^2}
\eeq
We show next that
\begin{equation}\label{eq:dhatEtaMp}
  \dd_{\what{K}}\eta^k(p) \in I_{\delta^2}\qquad\text{for any $k\geq 0$}
\end{equation}
by using induction on~$k$.

We note that $p\in\wtilde{\mathfrak{T}}(\fg)\subset \ker\dd_{\what{K}}$, therefore
\[ \dd_{\what{K}}(\eta^0(p)) = 0,\]
thus the base case is established.

Assume  now that $\dd_{\what{K}}(\eta^k(p))\in I_{\delta^2}$, then we have (since $[\dd_{\what{K}},\eta]=\zeta$, see Table~\eqref{tableder})
\[
  \dd_{\what{K}}(\eta^{k+1}(p)) =
  \eta(\dd_{\what{K}}(\eta^k(p))) - \zeta(\eta^k(p)).
\]
The first term belongs to~$I_{\delta^2}$ by the induction hypothesis and the fact that $\eta$
preserves~$I_{\delta^2}$ (see Lemma~\ref{lem:derIdelta2}).
The second term belongs to $I_{\delta^2}$ by Lemma~\ref{lem:M}(c).
This proves~\eqref{eq:dhatEtaMp} for any~$k$, in particular for $k=m$.

Note that $i\circ\eps(\eta^m(p)) =0$ since $\eta^m(p)\in\tW^{1,2m}(\fg)$.
Using~\eqref{dhat htpy} we get
\[
  \eta^m(p) = \dd_{\what{K}}\circ g \circ\varphi (\eta^m(p)) + g\circ\varphi( \dd_{\what{K}}(\eta^m(p))).
\]
The first term belongs to~$\ker\dd_{\what{K}}$. Since $g$ and $\varphi$ preserve~$I_{\delta^2}$,
\eqref{eq:dhatEtaMp} implies that the second term belongs to~$I_{\delta^2}$. Hence $\eta^m(p) \in \ker \dd_{\what{K}} + I_{\delta^2}$, thus proving ~\eqref{toprove 5} and Lemma~\ref{final lem}.
\end{proof}

\subsection{The quantum Weil algebra}\label{ssec:cWg}
In this subsection we define the transgression map $\ul{\mathsf{t}}_{\Cl(\fg)}$ for the quantum Weil
algebra~$\cW(\fg) = U(\fg)\otimes \Cl(\fg)$ with values in the Clifford algebra~$\Cl(\fg)$ (Definition~\ref{def:uTClg}).
The main result of this subsection is to show that $\iota_x \ul{\mathsf{t}}_{\Cl(\fg)} (U^+(\fg)^\fg)$ is
contained in the image of~$\alpha_\fg$. We have by Theorem~\ref{thm:iotaTtilde} that $\iota_x
\wtilde{\mathsf{t}}_{\Cl(\fg)}(\tilde{\frT}(\fg))$ is contained in the image of~$\Phi_{\Cl(\fg)}$.
Using the universal property of~$\tW(\fg)$ we show that the images of transgressions in $\tW(\fg)$ and
$\cW(\fg)$ with values in~$\Cl(\fg)$ coincide (Lemma~\ref{lem:imTsame}), as do the images of quantum moment maps
$\Phi_{\Cl(\fg)}$ and~$\alpha_\fg$ (Lemma~\ref{lem:imalphaisimphi}).
Thus the main result of this subsection (Theorem~\ref{thm:iotaxImAlphaG}) is a~direct consequence of the
universal property of~$\tW(\fg)$ and~Theorem~\ref{thm:iotaTtilde}.

As before, suppose that~$\fg$ admits a~nondegenerate invariant symmetric bilinear form~$B$, so we can identify
$\fg^\ast \cong \fg$.
Recall from Example~\ref{ex:ClisGDiff} that $\Cl(\fg)$ is a $\fg$-differential algebra with the flat connection $\vartheta_{\Cl(\fg)}$ given by
$\vartheta_{\Cl(\fg)}(\mu)=\mu$. Let $\wtilde{\pi}_{\Cl(\fg)}\colon \tW(\fg)\to \Cl(\fg)$ be the corresponding
characteristic homomorphism. For $\mu\in\fg^\ast$ we get that
\[
\wtilde{\pi}_{\Cl(\fg)}(\mu) = \wtilde{\pi}_{\Cl(\fg)}\circ \vartheta_{\tW}(\mu)=\vartheta_{\Cl(\fg)}(\mu)= \mu
\]
and it follows from Lemma~\ref{lem:Aflat} that $\wtilde{\pi}_{\Cl(\fg)}(\what{\mu})=0$.

Define the quantum Weil algebra as
\[
\cW(\fg) = U(\fg) \otimes \Cl(\fg),
\]
and let $\what{e}_i$ denote generators of~$U(\fg)\otimes 1$ and $e_i$ denote generators of $1\otimes\Cl(\fg)$.
Consider the cubic Dirac element~\cite{GoetteDirac,Kostant1999cubic}
\[
  \mathcal{D} = \sum_i^{(\fg)} \what{e}_if_i + q(\phi) \in \cW(\fg),
\]
where $\phi$ is the Cartan 3-tensor defined by~\eqref{eq:Cartan3tensor} and we consider $f_i\in\fg^\ast$, dual to $e_i$, as
elements of~$\fg$ via the above identification $\fg^\ast\cong\fg$.

The quantum Weil algebra~$\cW(\fg)$ is a filtered $\fg$-differential algebra with the differential given by
$\dd_{\cW} = [\mathcal{D}, - ]$, with the Lie derivatives induced by the adjoint action, and with the contractions defined by $\iota_x = \id\otimes\iota_x$, $x\in\fg$.
The filtration is induced by the usual filtrations on $U(\frg)$ and on $\Cl(\frg)$: $\deg\what{e}_i=2$, $\deg e_i=1$. 
We have that
\begin{equation}\label{eq:duWmu}
  \dd_{\cW}\mu = 2\what\mu + 2\alpha_\fg(\mu),
\end{equation}
recalling $\alpha_\fg$ is defined in Equation \ref{eq:def of alpha}.

\begin{rem}
  Let us compare the explicit formulas for differentials in~$\tW(\fg)$, $W(\fg)$, and~$\cW(\fg)$.
  For $\mu\in\fg^\ast$ we have
  \begin{align*}
    \text{in $\tW(\fg)$:}& & & \dd_{\tW}(\mu) = \what{\mu} + \delta(\mu),\\
    \text{in $W(\fg)$:}& & & \dd_{W}(\mu) = 2\what{\mu} + 2\lambda_\fg(\mu),\\
    \text{in $\cW(\fg)$:}& & &\dd_{\cW}(\mu) = 2\what{\mu} + 2\alpha_\fg(\mu).
  \end{align*}
  The factor~$2$ that appears when we pass to~$W(\fg)$ and~$\cW(\fg)$ is related to the Poisson structure
  on~$W(\fg)$, see the discussion in~\cite[\S7]{MeinrenkenBook}.
\end{rem}

One can define~$\ul{\pi}_{\Cl(\fg)}: \cW(\fg) \to \Cl(\fg) $ as the projection along~$U^+(\fg) \otimes 1$,
where $U^+(\fg)$ is the augmentation ideal~of $U(\fg)$. 
This projection is a~$\fg$-differential algebra morphism.
From \cite[Section 7.3]{MeinrenkenBook} we have an isomorphism of~$\fg$-differential spaces~$q_{W} : W(\fg)
\to \cW(\fg)$ such that 
\eq\label{q and pi}
\ul{\pi}_{\Cl(\frg)} \circ q_W =q \circ  \pi_{\wedge(\frg)},
\eeq
and its restriction to $S(\fg^\ast)\subset W(\fg)$, generated by~$\what{e}_i$,
is the Duflo map $S(\fg^\ast) \cong S(\frg)\to U(\fg)$, see~\cite{Dulfo1977}.
(Recall that $q\colon\twedge\fg \to \Cl(\fg)$ is the quantisation map, see~\S\ref{sec:CliffAlg}.)
Since $q_W$ is an isomorphism of~$\fg$-differential spaces, we have that $H(\cW(\fg)) = \mathbb{C}$.

Using~$\ul{\pi}_{\Cl(\fg)}$ we can define another transgression.
\begin{Def}\label{def:uTClg}
For~$x \in U^+(\fg)^{\fg}$ denote by~$C_x \in \cW(\fg)^{\fg}$ a~cochain of transgression for~$x$, i.e., $\dd_{\cW}C_x = x$.
(Such a~cochain of transgression for~$x$ exists since there is a~cochain of transgression for~$q^{-1}_W(x)$.)
Then the transgression of~$x$ is given by
\[
  \ul{\mathsf{t}}_{\Cl(\fg)}(x) = \ul{\pi}_{\Cl(\fg)}(C_x).
\]
\end{Def}
Let us show that $\ul{\mathsf{t}}_{\Cl(\frg)}$ does not depend on a~choice of a~cochain.
Namely, suppose $C_x$ and $C_x'$ are two cochains of transgression for $x$, then $q_W^{-1}(C_x)$ and $q_W^{-1}(C_x')$
are cochains of transgression for $y := q_W^{-1}(x) \in W(\frg)$.
Since the transgression map $\mathsf{t}_{\wedge(\frg)}$ does not depend on a~choice of a~cochain
(Theorem~\ref{thm:transgression}\ref{item:tWellDef}),
we have that $\pi_{\wedge\frg}(q_W^{-1}(C_x - C_x')) = 0$. Thus we conclude that
\[
  \ul{\pi}_{\Cl(\frg)}(C_x - C_x') =q\circ \pi_{\wedge\frg} \circ q_W^{-1}(C_x - C_x') = 0,
\]
which shows that $\ul{\mathsf{t}}_{\Cl(\frg)}$ is independent of the choice of cochain hence well-defined.

\begin{ex}
  The Dirac element~$\mathcal{D}$ is a~cochain of transgression for the quadratic Casimir element
  $\mathrm{Cas}_\fg= \sum_i\what{e}_i\what{f}_i\in U(\fg)^\fg$. The image of~$\mathrm{Cas}_\fg$ under the
  transgression map $\ul{\mathsf{t}}_{\Cl(\fg)}$ is the quantisation of the Cartan 3-tensor~$q(\phi)$.
\end{ex}

The quantum Weil algebra~$\cW(\fg)$ also has a connection given by~$\vartheta_{\cW}(\mu) = \mu$, but this
connection is not flat.
By the universal property of the noncommutative Weil algebra~$\tW(\fg)$, there exists a
unique~$\fg$-differential algebra morphism, the characteristic homomorphism,~$\ul{c} : \tW(\fg) \to \cW(\fg)$
such that~$\ul{c} \circ \vartheta_{\tW} = \vartheta_{\cW}$. Before we calculate~$\ul{c}$ on generators let us
remark the following.
Identifying~$\fg^\ast$ and~$\fg$ via bilinear form~$B$, for~$f_i \in \fg^\ast\cong\fg$, we have
\[
\langle L_{f_i}e_a, e_b \rangle =-\langle L_{e_a}f_i,e_b \rangle =  \langle f_i, [e_a,e_b] \rangle = c_{a,b}^i,
\]
so~$L_{f_i}e_a = \sum_b c_{a,b}^i f_b$ and
\begin{equation}\label{eq:alphafi}
\alpha_{\fg}(f_i) = \frac{1}{4} \sum_{a}L_{f_i}e_a \cdot f_a = \frac{1}{4} \sum_{a,b} c_{a,b}^i f_b \cdot f_a 
=- \frac{1}{4} \sum_{a,b} c_{b,a}^i f_b \cdot f_a.
\end{equation}
Recalling~$\dd_{\cW}(\mu)  = 2 \what{\mu} + 2 \alpha_{\fg}(\mu)$ we get
\begin{subequations}\label{eq:ulc}
\begin{align}
  \ul{c}(f_i) &{}= \ul{c}(\vartheta_{\tW}(f_i)) = \vartheta_{\cW}(f_i) = f_i, \\
  \ul{c}(\what{f}_i) &= \ul{c} \left(
  \dd_{\tW} f_i + \frac{1}{2} \sum_{a,b} c_{a,b}^i \, f_a \otimes f_b \right) = \dd_{\cW}f_i + \frac{1}{2} \sum_{a,b} c_{a,b}^i \, f_a \cdot f_b \notag\\
                &{}= \dd_{\cW}f_i -2 \alpha_{\fg}(f_i) = 2(\what{f_i} + \alpha_{\fg}(f_i)) - 2 \alpha_{\fg} (f_i) = 2 \what{f}_i.
\end{align}
\end{subequations}

The main result of this subsection is the following theorem:

\begin{thm}\label{thm:iotaxImAlphaG}
  Let $\ul{\mathsf{t}}_{\Cl(\fg)}:U^+(\frg)^\frg\to \Cl(\frg)^\frg$ be the transgression map defined in Definition~\ref{def:uTClg}.
  Then for any $x\in\fg$, $p\in U^+(\fg)^\fg$, we have that $\iota_x\ul{\mathsf{t}}_{\Cl(\fg)}(p) \in \im \alpha_\fg$.  
\end{thm}

The proof uses the following lemma whose proof is postponed to the end of this subsection:

\begin{lem}\label{lem:imTsame}
  For $p\in S^{+}(\fg^\ast)^\fg$,
  we have~$\wtilde{\mathsf{t}}_{\Cl(\fg)}(\what{\sym}(p)) = \ul{\mathsf{t}}_{\Cl(\fg)}(\ul{c}(\sym_{W}(p)))$.
\end{lem}
\begin{proof}[Proof of Theorem~\ref{thm:iotaxImAlphaG}]
By Theorem~7.2 on p.~174 in~\cite{MeinrenkenBook}, the map~$q_W\colon W(\fg)\to\cW(\fg)$ factors through the
noncommutative Weil algebra:
\begin{equation}\label{eq:qWfactors}
  \begin{tikzcd}
    q_W\colon W(\fg) \arrow[r, "\sym_W"] & \tW(\fg)\arrow[r,"\ul{c}"] & \cW(\fg).
  \end{tikzcd}
\end{equation}
Therefore, Lemma~\ref{lem:imTsame} shows that the transgression map in~$\cW(\fg)$ factors through~$\tW(\fg)$ and
the diagram
\[
  \begin{tikzcd}
    S^+(\fg^\ast)^\fg \arrow[r,"\what{\sym}"]\arrow[d,"q_W"]  & \wtilde{\mathfrak{T}}(\fg)\arrow[d,"\wtilde{\mathsf{t}}_{\Cl(\fg)}"] \\
    U^+(\fg)^\fg \arrow[r,"\ul{\mathsf{t}}_{\Cl(\fg)}"] & \Cl(\fg)^\fg
  \end{tikzcd}
\]
commutes. 

We next claim that $\im \wtilde{\mathsf{t}}_{\Cl(\fg)} = \im \ul{\mathsf{t}}_{\Cl(\fg)}$. To conclude this, we must show that $\ul{c}$ restricts to a~surjection of the transgression spaces. The transgression space for $\cW(\frg)$ is $U^+(\frg)^\frg$, so must prove that $\ul{c}(\frT(\frg)) = U^+(\frg)^\frg$.

  Since the restriction of~$q_W$ to~$S(\fg^\ast)\cong S(\frg)$ is the Duflo map, then $U^+(\frg)^\frg = q_W(S^+(\frg^\ast)^\frg)$. The transgression space $\frT(\frg)$ equals $\sym_W(S^+(\frg)^\frg)$, 
  and from~\ref{eq:qWfactors} we get
  \[U^+(\frg)^\frg = q_W(S^+(\frg^\ast)^\frg) = \ul{c} (\frT(\frg)).\]
  Therefore $\im \underline{\mathsf{t}}_{\Cl(\frg)} = \im\wtilde{\mathsf{t}}_{\Cl(\frg)}$.   Furthermore, Theorem~\ref{thm:iotaTtilde} shows that $\iota_x (\im \wtilde{\mathsf{t}}_{\Cl(\frg)}) \subset \Phi_{\Cl(\frg)}\left( \tW^{\bullet,0}(\fg) \right)$,
  hence
  \[
    \iota_x (\im \underline{\mathsf{t}}_{\Cl(\frg)}) \subset  \Phi_{\Cl(\frg)}\left( \tW^{\bullet,0}(\fg) \right).
  \]
  Theorem~\ref{thm:iotaxImAlphaG}
  now follows from the following Lemma.
  \end{proof}

\begin{lem}\label{lem:imalphaisimphi}  
  We have that
  \(
    \Phi_{\Cl(\fg)}\left( \tW^{\bullet,0}(\fg) \right) = \alpha_\fg( U(\fg)).
  \)
\end{lem}
\begin{proof} 
  Since both $\Phi_{\Cl(\fg)}$ and $\alpha_{\frg}$ are algebra morphisms
  and the map $\ul{c}$ is surjective,
  it suffices to show that~$\Phi_{\Cl(\fg)}(\what{\mu})$ and $\alpha_{\frg}(\ul{c}(\what{\mu}))$ are equal for~$\mu \in \fg^{\ast}$.
  Take a generator~$f_i$ as before: we have, using~\eqref{eq:Phi} and~\eqref{eq:PhiA},
\begin{align*}
  \Phi_{\Cl(\fg)}(\what{f}_i)&{}= \wtilde{\pi}_{\Cl(\fg)}(\delta(f_i)) = \wtilde{\pi}_{\Cl(\fg)} \left( - \frac{1}{2} \sum_{a,b} c_{a,b}^i \, f_a \otimes f_b \right) = - \frac{1}{2}\sum_{a,b} c_{a,b}^i \, f_a \, f_b. \\
\intertext{From~\eqref{eq:alphafi} we get that}
  \alpha_{\fg}(\ul{c}(\what{f}_i))
 &{}= \alpha_{\fg}(2\what f_i)= - \frac{1}{2} \sum_{a,b} c_{b,a}^i\,  f_b\, f_a = \Phi_{\Cl(\fg)}(\what f_i).
\end{align*}
This proves the claim.
\end{proof}

\medskip

\noindent{\it Proof of Lemma \ref{lem:imTsame}.}
  For~$p \in S^+(\fg^\ast)^\fg$ denote~$\wtilde{p}:=\what{\sym}(p)\in\wtilde{\mathfrak{T}}(\fg)$,
  then $\ul{p}:= \ul{c}(\wtilde{p})\in U^+(\fg)^\fg$.
   We claim that~$\ul{c}(h(\wtilde{p}))$ is a cochain of transgression for~$\ul{p}$. Since~$\ul{c}$
   intertwines the differentials, $\ul{c}(h(\dd_{\tW}(\wtilde{p}))) = 0$ (Lemma~\ref{lem:M}(c), Lemma~\ref{lem:derIdelta2} and Lemma~\ref{lem:Aflat}(\ref{item:Aflat3})), and
   $\dd_{\tW} \circ h + h \circ \dd_{\tW} = \id -i \circ \varepsilon$ (Lemma~\ref{lem htpy}), we find
  \[
    \dd_{\cW}\circ \mathop{\ul{c}} (h(\wtilde{p})) = \ul{c}\circ \dd_{\wtilde{W}}(h(\wtilde{p})) =
    \ul{c}(\wtilde{p}) = \ul{p}.
  \]
  Furthermore, we have~$\wtilde{\pi}_{\Cl(\fg)} = \ul{\pi}_{\Cl(\fg)} \circ \ul{c}$;
  using the definition of the characteristic homormorphism~$\wtilde{\pi}_{\Cl(\fg)}$ and Lemma~\ref{lem:Aflat}, this can be checked on generators:
  \begin{align*}
      &\wtilde{\pi}_{\Cl(\fg)}(\mu) = \mu,  & &\ul{\pi}_{\Cl(\fg)} \circ \ul{c}(\mu) = \ul{\pi}_{\Cl(\fg)}(\mu) = \mu, \\
      &\wtilde{\pi}_{\Cl(\fg)}(\what{\mu}) = 0,  & &\ul{\pi}_{\Cl(\fg)}\circ \ul{c}(\what{\mu}) = \ul{\pi}_{\Cl(\fg)}(2 \what{\mu}) = 0.
  \end{align*}
  Finally, since $\ul{c}(\sym_W(p)) = \ul{c} (\what{\sym}(p))$, we get that
  \[
    \ul{\mathsf{t}}_{\Cl(\fg)}(\ul{c}(\sym_{W}(p))) 
    = \ul{\pi}_{\Cl(\fg)}\circ \ul{c} (h(\wtilde{p})) = \wtilde{\pi}_{\Cl(\fg)} \circ h(\wtilde{p}) = \wtilde{\mathsf{t}}_{\Cl(\fg)}(\what{\sym}(p)).    
  \]
  This proves the claim.
\qed

\subsection{Relating transgressions in the commutative and quantum Weil algebras}
In this subsection we show that the images of the transgression maps in the commutative and quantum Weil
algebras are related by the Chevalley map~$q$.

\begin{lem}\label{lem:qintprims}
  For $p\in S^+(\fg^\ast)^\fg$, we have that $q\circ \mathsf{t}_{\wedge\fg}(p) = \ul{\mathsf{t}}_{\Cl(\fg)}(q_{W}(p))$. 

\end{lem}
\begin{proof}
  Since~$q_W$ intertwines the differentials we have: if~$C_p$ is a cochain of transgression for~$p$,
  then~$q_W(C_p)$ is a cochain of transgression for~$q_W(p)\in U^+(\fg)^\fg$. The claim now follows from \eqref{q and pi}.
\end{proof}

\begin{cor}
  We have that $\im \ul{\mathsf{t}}_{\Cl(\fg)} = q(P_\wedge(\fg)) $, where $P_\wedge(\fg)$ denotes the space of
  primitive invariants in~$\twedge\fg^\ast$, see Definition~\ref{def:PwedgeG}.
\end{cor}
\begin{proof}
  Recall from Theorem~\ref{thm:transgression}(\ref{item:Tiso}) that $P_{\wedge}(\frg)$ is equal to $\im
  \mathsf{t}_{\wedge(\frg)}$, hence the Corollary follows from  Lemma~\ref{lem:qintprims}.
\end{proof}

\section{Transgressions in relative Weil algebras}\label{sec:trans in rel weil}

The results of this section will allow us to extend Kostant's results~\cite{KostantRho} about the structure
of~$\Cl(\fg)^\fg$ to the relative case of~$\Cl(\fp)^K$.
Analogous to the ``absolute'' case, the transgression maps in various relative Weil algebras all factor through
the relative noncommutative Weil algebra.

\subsection{Relative noncommutative Weil algebra}\label{sec:tWrel}

Let $(\fg,\fk)$ be a symmetric pair of complex Lie algebras and $K$ be a~compact Lie group such that

\begin{enumerate}
\item $\mathrm{Lie}(K) =  \fk$,

\item the adjoint action of~$K$ on~$\fk$ extends to an~action on~$\fg$ which preserves~$B$,
  so $\fp$ is a~$K$-submodule.
\end{enumerate}
For such a~triple $\tW(\fg)$ is a~$K$-differential algebra.

In particular, we are interested in the following cases of 
$(\fg,\fk,K)$:
$(\fsl(2n+1),\fo(2n+1),SO(2n+1))$,
$(\fsl(2n), \fo(2n), O(2n))$,
$(\fsl(2n), \fsp(2n), Sp(2n))$,
$(\fe(6), \ff(4), F(4))$, where $F(4)$ is the compact connected simply connected group of type~$\ff(4)$.

Let $G$ be a~compact connected Lie group and $K$ be its symmetric subgroup,
then $(\mathrm{Lie}(G), \mathrm{Lie}(K), K)$ is an~example of such a~triple.
The only example that we consider which is not of this kind is $(\fsl(2n), \fo(2n), O(2n))$.

\begin{Def}
Define the \emph{relative noncommutative Weil algebra} $\tW(\fg,K)$ to be
the $K$-basic subalgebra of~$\tW(\fg)$, in other words the subalgebra of $\tW(\frg)^K$ which is annhilated by
every contraction by an element of $\frk$.
\end{Def}

In what follows we denote the derivation~$\delta$ of~$\tW(\fg)$ by~$\delta_\fg$.
Since $[\fk,\fp] \subset \fp$, 
for any $\mu\in\fk^\ast$ the element $\delta_\fg(\mu)\in\tW^{0,2}(\fg)$ given by~\eqref{eq:delta} decomposes as
\begin{equation}\label{eq:DeltaGasDeltaKP}
  \delta_\fg(\mu) = \delta_\fk(\mu) + \delta_\fp(\mu),
\end{equation}
where $\delta_\fk(\mu)\in T^2(\fk^\ast[-1])\subset\tW^{0,2}(\fk)\subset\tW^{0,2}(\fg)$ and $\delta_\fp(\mu) \in T^2(\fp^\ast[-1])\subset \tW^{0,2}(\fg)$.
It is easy to see that this decomposition of~$\delta_\fg(\mu)$ is $K$-equivariant.
Explicitly, for the basis vector $f_i\in\fk^\ast$
\begin{equation}\label{eq:deltaP}
  \delta_\fp(f_i) =  - \frac12 \sum_{a,b}^{(\fp)} c^i_{a,b} f_a\otimes f_b.
\end{equation}
We denote by $\Res^\fg_\fp$ the restriction of linear functionals from~$\fg$ to~$\fp$. For $\mu\in\fg^\ast$ we set
\begin{equation}\label{eq:PhiP}
  \Phi_\fp(\mu) = \Res^\fg_\fp\mu,\qquad
  \Phi_\fp(\what{\mu}) =
  \begin{cases}
    \delta_\fp(\mu) & \text{if $\mu\in\fk^\ast$}, \\
    0 &  \text{if $\mu\in\fp^\ast$},
  \end{cases}
\end{equation}
and extend it to an~algebra morphism
$\Phi_{\fp}\colon\tW(\fg)\to\tW(\fg)$. Clearly,
\begin{equation}\label{eq:PhiPeqResPhi}
  \Phi_{\fp}  = \Res^\fg_\fp\circ \Phi,
\end{equation}
where the algebra morphism $\Phi\colon \tW(\fg)\to\tW(\fg)$ is defined by~\eqref{eq:Phi}.
Moreover, $\Phi_\fp$ is $K$-equivariant and preserves the $\fk$-horizontal subspace of~$\tW(\fg)$, hence
$\Phi_\fp$ preserves~$\tW(\fg,K)$. In what follows we restrict~$\Phi_\fp$ to an~algebra morphism $\Phi_\fp\colon \tW(\fg,K)\to\tW(\fg,K)$.

Note that $\tW(\fg)$ is a~locally free $K$-differential algebra with connection given by~$\mu\mapsto \mu$ for
$\mu\in\fk^\ast$. Let $j\colon \tW(\fk)\to \tW(\fg)$ denote the corresponding characteristic homomorphism of
$K$-differential algebras.
For $\mu\in\fk^\ast$ we have
\begin{align}
  j(\mu) = {}& \mu,\label{eq:j}\\
  j(\what{\mu}) =
  {}& j(\dd_{\tW}(\mu) - \delta_{\fk}(\mu))
      = \dd_{\tW}(\mu) - \delta_{\fk}(\mu)
      = \what{\mu} + \delta_{\fg}(\mu) - \delta_{\fk}(\mu)
      = \what{\mu} + \delta_{\fp}(\mu).\notag
\end{align}
Note that $\mu$ and $\what{\mu}$ denote elements of both $\tW(\fk)$ and~$\tW(\fg)$; likewise $\dd_{\tW}$ denotes the
differential of both algebras.

We identify $\frg^*=\frk^*\oplus\frp^*$ in the usual way, and denote $\mu=\mu_\frk+\mu_\frp$ accordingly.
Define the projection $\pr_{\fk}:\tW(\fg) \to \tW(\fk)$ to be the algebra morphism given on the generators by
\begin{equation}\label{eq:DefPrK}
  \pr_{\fk}(\mu) = \mu_{\fk},\qquad \pr_{\fk}(\bar{\mu}) = \overline{\mu_{\fk}}.
\end{equation}
Then~$\pr_{\fk}$ is a~$K$-differential algebra morphism which satisfies
\[
\pr_{\fk}(\what{\mu}) = \pr_{\fk}(\what{\mu_{\fk}} + \what{\mu_{\fp}}) = \pr_{\fk}(\overline{\mu_{\fk}} - \delta_{\fg}(\mu_{\fk})+ \overline{\mu_{\fp}} - \delta_{\fg}(\mu_{\fp})) = \overline{\mu_{\fk}} - \delta_{\fk}(\mu_{\fk})=\what{\mu_{\fk}},
\]
i.e., $\pr_\frk$ takes $\what\mu\in\tW(\frg)$ to~$\what{\mu_{\fk}}\in\tW(\fk)$.

Let~$\psi_1,\psi_2 : E \to E'$ be morphisms of~$K$-differential spaces~$(E,\dd)$ and~$(E',\dd
')$. A~\emph{$K$-homotopy} between~$\psi_1$ and~$\psi_2$ is a morphism of~$K$-differential spaces~$\psi: E \to \bbC
[t, dt] \otimes E'$ such that
\[
\psi_1 = (\ev (0) \otimes \id) \circ \psi, \quad \psi_2 = (\ev (1) \otimes \id) \circ \psi,
\]
where $\ev(a)\colon \mathbb{C}[t,dt]\to\mathbb{C}$ is defined by $\ev(a)(P+Qdt)=P(a)$.
If such a $\psi$ exists, we say that~$\psi_1$ and~$\psi_2$ are~\emph{$K$-homotopic}.
Similarly, a~\emph{$K$-homotopy operator} between~$\psi_1$ and~$\psi_2$ is a morphism of super spaces~$h: E[1] \to E'$ such that
\begin{align*}
&{}h \circ \iota_{x} + \iota_{x} \circ h = 0,\qquad \forall x \in \fk, \\
&{}h \circ d + d' \circ h = \psi_2 - \psi_1.
\end{align*}
The maps~$\psi_1$ and~$\psi_2$ are~$K$-homotopic if and only there exists a~$K$-homotopy operator between them. We need the second part of this statement; assume there is a~$K$-homotopy~$\psi$ between~$\psi_1$ and~$\psi_2$, then~$h:= (J \otimes \id) \circ \psi$ defines a~$K$-homotopy operator between~$\psi_1$ and~$\psi_2$. Here~$J$ is the same ``integration operator" as in \ref{eq:J}.

We will say that~$K$-differential space morphisms~$\psi: E \to E'$ and~$\varphi: E' \to E$ are $K$-homotopy inverses if~$\psi \circ \varphi$ is homotopic to~$\id_{E'}$ and~$\varphi \circ \psi$ is homotopic to~$\id_E$. If a~$K$-differential space morphism admits a~$K$-homotopy inverse it is called a~$K$-homotopy equivalence.

Let us introduce a differential~$\dd_{\otimes}$ on~$\bbC[t,dt] \otimes \tW(\fg)$ by~$\dd_{\otimes}(t) = dt$,
$\dd_{\otimes}(dt) = 0$, $\dd_{\otimes}(\mu) =  \bar{\mu}$.
In this way $\bbC[t,dt] \otimes \tW(\fg)$ becomes a~$K$-differential algebra with the trivial action of $K$
and $\iota_x$ on~$\bbC[t,dt]$.

\begin{lem}\label{lem:tWrelHomotopy}
    Let~$\pr_{\fk}: \tW(\fg) \to \tW(\fk)$ be as above. Then~$\pr_{\fk}$ is a $K$-homotopy equivalence, with inverse $j$. In more detail, $\pr_{\fk} \circ j = \id_{\tW(\fk)}$, while a $K$-homotopy between~$j \circ \pr_{\fk}$ and~$\id_{\tW(\fg)}$ is given by
    \[
    \psi_{\rel}: \tW(\fg) \to \bbC[t, dt] \otimes \tW(\fg),\qquad \psi_{\rel}(\mu) = (1-t) \otimes \mu + t \otimes j \circ \pr_{\fk}(\mu),\quad \mu \in \fg^{\ast},
    \]
    extended as a morphism of differential algebras using the universal property of~$\tW(\fg)$
    viewed as a~noncommutative Koszul algebra; see~\cite[Proposition~6.4 in~\S6.4]{MeinrenkenBook}.   
\end{lem}
\begin{proof}
  Obviously we have~$\pr_{\fk} \circ j = \id_{\tW(\fk)}$.
  Furthermore, it is clear that
  \[
    \id_{\tW(\frk)}= (\ev (0) \otimes \id) \circ \psi_\rel, \quad j \circ \pr_{\fk} = (\ev (1) \otimes \id) \circ \psi_\rel,
  \]
  so it suffices to show that $\psi_\rel$ is a $K$-differential algebra morphism.
  Note first that~$\pr_{\fk}$ is a morphism of~$K$-differential algebras, as is~$j$, so~$j \circ \pr_{\fk}$ is a~$K$-differential algebra morphism. 
  To see that $\psi_{\rel}$ is a $K$-differential algebra morphism, we first describe it on generators
  \begin{align}
  \psi_{\rel}(\mu) &{}= (1-t)\otimes \mu + t \otimes j \circ \pr_{\fk}(\mu) = 1 \otimes \mu - t \otimes \mu_{\fp},\notag\\
    \psi_{\rel}(\bar{\mu}) &{}= \psi_{\rel}(\dd \mu) = \dd_{\otimes}(\psi_{\rel}(\mu)) = \dd_{\otimes} (1 \otimes \mu - t \otimes \mu_{\fp}) \\
    &= 1\otimes \bar{\mu} - dt \otimes \mu_{\fp} - t \otimes \overline{\mu_{\fp}}\label{eq:psiRelBarMu}.
\end{align}
Since~$\psi_{\rel}$ is an algebra morphism, it is enough to check the conditions on generators.
Take~$x \in \fk$ and~$\mu \in \fg^{\ast}$ and note that
\[
  L_xt = L_x dt = 0,\qquad \iota_xt = \iota_x dt = 0.
\]
So we get
  \begin{align*}
    \iota_x \psi_{\rel}(\mu) = {}
    &\iota_x (1 \otimes \mu - t \otimes \mu_{\fp}) = 1 \otimes \iota_x \mu - t \otimes \iota_x \mu_{\fp}
      = 1\otimes \iota_x \mu = \psi_{\rel}(\iota_x \mu), \\
    \intertext{(The last equality follows from the fact that $\iota_x\mu \in \bbC$.)}
    \iota_x \psi_{\rel} (\bar{\mu}) ={}
    & \iota_x(1 \otimes \bar{\mu} - dt \otimes \mu_{\fp} - t \otimes \overline{\mu_{\fp}}) = 1 \otimes \iota_x \bar{\mu} - t \otimes \iota_x \overline{\mu_{\fp}}  \\
    {}={} & 1 \otimes L_x \mu - t \otimes L_x \mu_{\fp} = \psi_{\rel}(L_x \mu) = \psi_{\rel}(\iota_x \bar{\mu}).
  \end{align*}
  It now follows from Cartan's magic formula that $\psi_{\rel}$ commutes also with Lie derivatives,
  which proves the claim.
\end{proof}

We denote the $K$-homotopy operator between $\id_{\tW}$ and $j\circ \pr_{\fk}$ corresponding to~$\psi_{\rel}$
by
\begin{equation}\label{eq:hrel}
  h_{\rel} = (J \otimes \id)\circ \psi_{\rel}.
\end{equation}

\begin{rem}

\label{rem:homotopyop}   
  Since $h_{\rel}$ is a $K$-homotopy operator, it commutes with the action of $K$ on~$\tW(\fg)$ and  $\iota_x$
  for all $x\in\fk$. Therefore, $h_{\rel}$ preserves the $K$-basic subalgebra~$\tW(\fg,K)$,
  hence $h_{\rel}(p)\in\tW(\fg,K)$. We use the $K$-homotopy operator~$h_{\rel}$ below
  to define the transgression map in~$\tW(\fg,K)$.

For $\mu \in \frg^\ast$, $\delta^2(\mu) \in T^3(\frg^\ast[-1])$ is skew symmetric, as is each of its $(\frk^*,\frk^*,\frk^*)$,
$(\frk^*,\frk^*,\frp^*)$, $(\frk^*,\frp^*,\frp^*)$ and $(\frp^*,\frp^*,\frp^*)$ parts. Here the $(X,Y,Z)$ part of $\delta^2(\mu)$ is the sum of monomials in $\delta^2(\mu)$ with a variable in each of $X$, $Y$ and $Z$ in any order. 
Therefore $\ul{c}$ of each of these parts is contained in $q (\twedge^3\frg^\ast)$ (see \cite[\S 2.2.6]{MeinrenkenBook}) and thus these parts cannot cancel each other, therefore each part is in the kernel of
$\ul{c}$.
Since $c$ is graded each part is also in the kernel of $c$.  The different parts of $\delta^2(\mu)$ are each
in the kernel and $\psi_\rel$ multiplies each part $\delta^2(\mu)_{(X,Y,Z)}$ by $(1-t)^{\#\frp}$ where ${\#
  \frp}$ is the number of $\frp^\ast$ variables in ${(X,Y,Z)}$. Thus $\psi_\rel(\delta^2(\mu)) \in (\ker \id
\otimes c) \cap \ker (\id \otimes \ul{c})$. Using the same argument as Lemma \ref{lem:derIdelta2} we find that $h_\rel(I_{\delta^2}) \subset \ker c \cap \ker \ul{c}$. It follows that $h_\rel(\dd_{\tW}(\what{p}))$ is in $\ker c \cap \ker \ul{c}$.

\end{rem}

\begin{lem}\label{lem:PsiRel}
    For~$\mu \in \fk^{\ast}$ we have
    \[
    \psi_{\rel}(\what{\mu}) = 1 \otimes \what{\mu} + (2t - t^2) \otimes \delta_{\fp}(\mu)
    \]
    and for~$\mu \in \fp^{\ast}$ we have
    \[
    \psi_{\rel}(\what{\mu}) = (1-t) \otimes \what{\mu} - dt \otimes \mu.
    \]
\end{lem}
\begin{proof}

To calculate~$\psi_{\rel}(\what{\mu})$ we need~$\psi_{\rel}(\delta_{\fg}(\mu))$. Assume that our basis~$f_i$ of $\frg^*$ is split into~$\fk^\ast$ and~$\fp^\ast$ parts, meaning that every element belongs to exactly one of them.
\begin{align*}
    \psi_{\rel}&(\delta_{\fg}(f_i))  = \psi_{\rel}(- \frac{1}{2}\sum_{a,b}^{(\fg)} c_{ab}^{i} f_a \otimes f_{b}) = - \frac{1}{2} \sum_{a,b}^{(\fg)}c_{ab}^i \psi_{\rel}(f_a) \otimes \psi_{\rel}(f_b) \\
    &= - \frac{1}{2} \sum_{a,b}^{(\fg)} c_{a,b}^{i} \left(1 \otimes f_a - t \otimes (f_a)_{\fp}\right) \cdot \left( 1 \otimes f_b - t \otimes (f_b)_{\fp} \right) \\
    &= - \frac{1}{2}\sum_{a,b}^{(\fg)}c_{ab}^i \left( 1 \otimes f_a \otimes f_b - t \otimes f_a \otimes (f_b)_{\fp} - t \otimes (f_a)_{\fp} \otimes f_b + t^2 \otimes (f_a)_{\fp} \otimes (f_b)_{\fp} \right) \\
    &{}= 1 \otimes \delta_{\fg}(f_i) + \frac{t}{2} \otimes \left(\sum_{a}^{(\fg)}\sum_{b}^{(\fp)} c_{ab}^i f_a\otimes f_b + \sum_{a}^{(\fp)} \sum_{b}^{(\fg)} c_{ab}^i f_a\otimes f_b \right) - \frac{t^2}{2} \otimes \sum_{a,b}^{(\fp)}c_{ab}^i f_a\otimes f_b.
\end{align*}
If~$f_i \in \fk^{\ast}$, then both sums are over basis of~$\fp^{\ast}$ because~$c_{ab}^i = 0$ if~$f_a \in \fk^{\ast}, f_b \in \fp^{\ast}$ or the other way around. Also,~$\sum_{a,b}^{(\fp)}c_{ab}^i f_a \otimes f_b = - 2\delta_{\fp}(f_i)$ so we get
\begin{equation}\label{eq:psiRelDeltaFik}
\psi_{\rel}(\delta_{\fg}(f_i)) = 1 \otimes \delta_{\fg}(f_i) - 2t \otimes \delta_{\fp}(f_i) + t^2 \otimes \delta_{\fp}(f_i),\quad\text{if $f_i\in\fk^\ast$}.
\end{equation}

On the other hand, if~$f_i \in \fp^{\ast}$, then
\[
\sum_{a}^{(\fg)} \sum_{b}^{(\fp)} c_{ab}^i f_a \otimes f_b + \sum_{a}^{(\fp)} \sum_{b}^{(\fg)} c_{ab}^i f_a \otimes f_b = \sum_{a}^{(\fk)} \sum_{b}^{(\fp)} c_{ab}^i f_a \otimes f_b + \sum_{a}^{(\fp)} \sum_{b}^{(\fk)} c_{ab}^i f_a \otimes f_b = -2\delta_{\fg}(f_i),
\]
and the last summand in the expression above vanishes because~$c_{ab}^i = 0$ if~$f_i, f_a, f_b \in \fp^{\ast}$. 
(Here we are using that~$\frk$ is symmetric in $\frg$.)
We obtain
\begin{equation}\label{eq:psiRelDeltaFip}
  \psi_{\rel}(\delta_{\fg}(f_i)) = (1-t) \otimes \delta_{\fg}(f_i),\quad\text{if $f_i\in\frp^*$}.
\end{equation}

To summarize, if~$\mu \in \fk^{\ast}$, then, using the formulas~\eqref{eq:psiRelBarMu},
\eqref{eq:psiRelDeltaFik} and~\eqref{eq:psiRelDeltaFip},
\begin{align*}
  \psi_{\rel}(\what{\mu}) = {}& \psi_{\rel}(\bar\mu - \delta_\fg(\mu))
  = 1 \otimes \what{\mu} + (2t-t^2) \otimes \delta_{\fp}(\mu), &{}&\text{if~$\mu \in \fp^{\ast}$},\\
  \psi_{\rel}(\what{\mu}) ={}
   &\psi_{\rel}(\bar\mu - \delta_\fg(\mu)) = (1-t) \otimes \what{\mu} - dt \otimes \mu,
   &{}&\text{if~$\mu \in \fk^{\ast}$}.
\end{align*}
Which concludes the proof.
\end{proof}

\subsection{Transgression in the relative noncommutative Weil algebra}
Let $\cA$ be a~flat locally free $\fg$-differential algebra such that the action of $\fk\subset\fg$ integrates to
an~action of~$K$.
It follows that $\cA$ is a $K$-differential algebra compatible with the
$\fg$-differential algebra structure.
Set $\cB = \cA_{K-\bas}$,
the $K$-basic subalgebra of $\cA$, and let $\cA_{K-\hor}$ be the $K$-horizontal subalgebra of~$\cA$;
see Subsection~\ref{ssec:GdiffAlg}.
It is easy to see that $\Phi_{\fp}(\tW(\fg)) \subset \tW(\fg)_{K-\hor}$, so $\wtilde{\pi}_{\cA}\circ \Phi_{\fp}(\tW(\fg))\subset\cA_{K-\hor}$.
(Here $\Phi_\fp$ is defined by~\eqref{eq:PhiP}.)
Define the algebra homomorphism
\begin{equation}\label{eq:PhiAhor}
  \Phi_{\cA}^{\hor} := \wtilde{\pi}_{\cA}\circ \Phi_{\fp} \colon\tW(\fg)\to\cA_{K-\hor}.
\end{equation}
Let $\wtilde{\pi}_{\cB}\colon \tW(\fg,K)\to \cB$ be the restriction of~$\wtilde{\pi}_{\cA}$ to $K$-basic subalgebras.
Note that $\wtilde{\pi}_{\cB}$ is a~morphism of differential algebras.

\begin{Def}
  Analogously to the $\tW(\fg)$ case, we define the \emph{relative transgression space} of~$\tW(\fg,K)$ by
  \[
    \wtilde{\mathfrak{T}}(\fg,K) = \left\{
      \what{\sym}(p) \mid
      \text{$p \in S^+(\fg^\ast)^K$ such that $\pr_{\fk}(\what{\sym}(p)) = 0$}
    \right\} \subset \tW^{+,0}(\fg,K).
  \]

\end{Def}

\begin{lem}\label{lem:prkhrelp}
  For $p\in\wtilde{\frT}(\fg,K)$ we have that $\pr_\fk h_\rel(p) = 0$.
\end{lem}
\begin{proof}We have that
  \begin{align*}
    \pr_\fk h_\rel(p) =
    {}& \pr_\fk\circ (J\otimes\id)\circ \psi_\rel (p)
        = (J\otimes\id)\circ(\id\otimes\pr_\fk)\circ \psi_\rel (p).
  \end{align*}
  Since $\psi_\rel$ and $\pr_\fk$ are morphisms of $K$-differential algebras
  and every monomial in~$p$ must contain at least one factor $\what{\mu}$ for $\mu\in\fp^\ast$,
  to show that $(\id\otimes\pr_\fk)\circ \psi_\rel (p)=0$,
  it is enough to check this on (subset of) generators $\what{\mu}$, $\mu\in\fp^\ast$.
  Using the universal  property of~$\tW(\fg)$ and Lemma~\ref{lem:PsiRel}, we see:
  \[
    (\id\otimes\pr_\fk)\circ \psi_\rel (\what{\mu}) =
    (1-t) \otimes\pr_\fk\what{\mu} - dt \otimes\pr_\fk\mu = 0.
  \]
  Which proves the claim.
\end{proof}

\begin{Def}\label{def:relTtW}
The \emph{relative transgression map} in~$\tW(\fg,K)$ with values in~$\cB$ is the~linear map
\[
  \wtilde{\mathsf{t}}_{\cB} = \wtilde{\pi}_{\cB}\circ h_{\rel} \colon \wtilde{\mathfrak{T}}(\fg,K) \to \cB.
\]
\end{Def}

\begin{prop}\label{prop:reltrformula}
  Let $p\in \wtilde{\mathfrak{T}}^{m+1}(\fg,K)$, then
  \eq\label{tilde t 1}
    \wtilde{\mathsf{t}}_{\cB}(p)  = - \frac{4^m \cdot (m!)^2}{(2m+1)!} \Phi_{\cA}^{\hor}(g(p)),
  \eeq
  where $g$ is defined in Table~\eqref{tableder} and $\Phi_{\cA}^{\hor}$ is defined by~\eqref{eq:PhiAhor}.
  Explicitly, let
  \eq\label{dec p}
    p = p' + \sum\what{\mu}_{i_1} \otimes \dots \otimes \what{\mu}_{i_{k-1}} \otimes \what{\xi}_{i_k} \otimes  \what{\mu}_{i_{k+1}} \otimes \dots \otimes \what{\mu}_{i_{m+1}} \in \wtilde{\mathfrak{T}}^{m+1}(\fg,K)
  \eeq
  such that~$\mu_{i_j} \in \fk^\ast$, $\xi_{i_j} \in \fp^\ast$ and $p'$ does not contain any monomials of this form.
  Then
  \eq\label{tilde t 2}
    \wtilde{\mathsf{t}}_{\cB}(p) =- \frac{4^m \cdot (m!)^2}{(2m+1)!} \sum  \Phi_{\cA}^{\hor}(\what{\mu}_{i_1} \otimes \dots \otimes\what{\mu}_{i_{k-1}}) \otimes  \wtilde{\pi}_{\cA}({\xi}_{i_k}) \otimes \Phi_{\cA}^{\hor}(\what{\mu}_{i_{k+1}}  \otimes \dots \otimes\what{\mu}_{i_{m+1}}).
  \eeq
\end{prop}
\begin{proof}
  By the definition~\eqref{eq:hrel} of~$h_{\rel}$,
  on~$\wtilde{\mathfrak{T}}(\fg,K)$ we have
  \[
    \wtilde{\mathsf{t}}_{\cB} = \wtilde{\pi}_{\cB} \circ h_{\rel} = \wtilde{\pi}_{\cB} \circ (J \otimes \id) \circ \psi_{\rel}.
  \]  
  Note that $\wtilde{\pi}_{\cA}\circ h_{\rel}$ is a well-defined map on~$\tW(\fg)$ which restricts
  to~$\wtilde{\pi}_{\cB}\circ h_{\rel}$ on~$\tW(\fg,K)$.
  Therefore we can work with monomials
  \eq\label{monos nonzero}
  \what{\mu}_1 \otimes \dots \otimes \what{\mu}_{m+1} \in \tW^{m+1,0}(\fg),
  \eeq
  even though they are not in the relative transgression space $\wtilde{\mathfrak{T}}(\fg,K)$. Moreover, since every $\what\mu$ is the sum of its $\frk$-part and its $\frp$-part, we can assume that each $\mu_i$ is either in $\frk^\ast$ or in $\frp^\ast$.
  
  Let us see which of the monomials \eqref{monos nonzero} will
  survive the application of the map~$\wtilde{\pi}_{\cA}\circ h_{\mathrm{rel}}$.
  
  Recall from Lemma \ref{lem:PsiRel} that
\[
    \psi_{\rel}(\what{\mu}) = 
    \begin{cases}  1 \otimes \what{\mu} + (2t - t^2) \otimes \delta_{\fp}(\mu), \quad  \mu\in\frk^\ast;\\
     (1-t) \otimes \what{\mu} - dt \otimes \mu, \qquad\quad\, \mu \in \fp^{\ast}.
    \end{cases}
\]
Now recall that by Lemma~\ref{lem:Aflat},  $\wtilde{\pi}_{\cA}$  annihilates~$\what{\mu}$, so the surviving terms can only come from $(2t - t^2) \otimes \delta_{\fp}(\mu), \   \mu\in\frk^\ast$ and from $- dt \otimes \mu,\ \mu \in \fp^{\ast}$. Moreover, the result will be 0 if there is more than one factor containing $dt$ (since $(dt)^2=0$), and also if there are no terms containing $dt$ (since then the result of applying $J$ will be 0).

To summarize, to survive the application of $\wtilde{\pi}_{\cA}\circ h_{\mathrm{rel}}$, the monomial \eqref{monos nonzero} must contain exactly one factor $\what\mu_j$ from $\frp^\ast$, while all the other factors $\what\mu_i,\,i\neq j,$ must be from $\frk^\ast$.  

For such a~monomial, since by \eqref{eq:PhiAhor} and \eqref{eq:PhiP} we have
\[
\wtilde{\pi}_{\cA}\circ \delta_{\fp}(\mu_i) = \Phi_{\cA}^{\hor}(\what\mu_i),\quad \mu_i\in\frk^*,
\]
we get
\begin{align*}
  {}\wtilde{\pi}_{\cA}&{}\circ h_{\rel}(\what{\mu}_1 \otimes \dots \otimes \what{\mu}_j \otimes \dots \otimes
    \what{\mu}_{m+1}) = \\
  &{}= (J \otimes \id)\Bigl(
    \bigl( (2t-t^2) \otimes \Phi_{\cA}^{\hor}(\what\mu_1)\bigr)
    \otimes \cdots \otimes \bigl( - dt \otimes \wtilde{\pi}_{\cA}(\mu_j) \bigr) \otimes \cdots\\
  {}&\phantom{{}={}}{}\cdots \otimes
      \bigl( (2t -t^2) \otimes \Phi_{\cA}^{\hor} (\what\mu_{m+1})\bigr)
      \Bigr) \\
  &{}= - \int_{0}^{1}(2t - t^2)^m dt \otimes \Phi_{\cA}^{\hor}(\what\mu_1) \otimes \dots \otimes \wtilde{\pi}_{\cA}(\mu_j)  \otimes \dots \otimes \Phi_{\cA}^{\hor}(\what\mu_{m+1}) \\
  &{}= - \frac{4^m \cdot (m!)^2}{(2m+1)!} \;  \Phi_{\cA}^{\hor}(\what\mu_1) \otimes \dots \otimes \wtilde{\pi}_{\cA}(\mu_j)  \otimes \dots \otimes \Phi_{\cA}^{\hor}(\what\mu_{m+1}).
\end{align*}
Now note that any $p \in \wtilde{\mathfrak{T}}^{m+1}(\fg,K)$ can be written as in \eqref{dec p}. The above discussion then implies \eqref{tilde t 2}. The  formula \eqref{tilde t 1} now follows directly from the definition of the derivation $g$ (which is in Table \eqref{tableder}).
\end{proof}

\begin{cor}\label{cor:abs_rel_tr}
  For all $\wtilde{p}\in \wtilde{\mathfrak{T}}^{m+1}(\fg,K)$ we have that
  \[
    \wtilde{\mathsf{t}}_{(\wedge\fp^\ast)^K}(\wtilde{p}) = -4^m\Res^\fg_\fp \circ\wtilde{\mathsf{t}}_{\wedge\fg^\ast}(\wtilde{p}),
    \qquad
    \wtilde{\mathsf{t}}_{\Cl(\fp)^K}(\wtilde{p}) = -4^m q\circ\Res^\fg_\fp\circ q^{-1}\circ \wtilde{\mathsf{t}}_{\Cl(\fg)}(\wtilde{p}).
  \]
\end{cor}

\begin{proof}
  For $\cA = \twedge\fg$ or $\cA = \Cl(\fg)$ and $\wtilde{p}\in\wtilde{\frT}^{m+1}(\fg,K)\subseteq\wtilde{\frT}^{m+1}(\fg)$,
  using Proposition~\ref{prop:trformula}, we have
  \begin{align*}
    \Res^\fg_\fp \circ\wtilde{\mathsf{t}}_{\cA}(\wtilde{p}) = {}
    & \frac{(m!)^2}{(2m+1)!}\Res^\fg_\fp\circ \Phi_{\cA}\circ g (\wtilde{p}) \\
    \intertext{(Since $\Phi_\cA=\wtilde\pi_\cA\circ\Phi$ and $\Res^\fg_\fp\circ\wtilde{\pi}_{\cA} = \wtilde{\pi}_{\cA}\circ\Res^\fg_\fp$)}
    {}={}& \frac{(m!)^2}{(2m+1)!}\wtilde{\pi}_{\cA}\circ \Res^\fg_\fp\circ\Phi\circ g(\wtilde{p}) \\
    \intertext{(Using that $\Phi_{\cA}^{\hor} = \wtilde{\pi}_{\cA}\circ \Phi_{\fp} = \wtilde{\pi}_{\cA}\circ \Res^\fg_\fp\circ \Phi$.)}
    {}={}& \frac{(m!)^2}{(2m+1)!}\Phi_{\cA}^{\hor}\circ g (\wtilde{p}) \\
    \intertext{(Using Proposition~\ref{prop:reltrformula}.)}
    {}={}& - \frac{1}{4^m}\wtilde{\mathsf{t}}_{\cB}(\wtilde{p}),
  \end{align*}
  which proves the claim.
\end{proof}

\begin{prop}\label{prop:contr_trans_formula}
  Let $\cA=\twedge\fg^\ast$ or $\cA=\Cl(\fg)$, so $\cB = \left(\twedge \fp^{\ast}\right)^K$, respectively $\cB = \Cl(\fp)^K$.
  For $x\in\fp$ and  $\wtilde{p}\in \wtilde{\mathfrak{T}}^{m+1}(\fg,K)$ we have
  \[
   \iota_x\wtilde{\mathsf{t}}_{\cB}(\wtilde{p}) = -\frac{4^m(m!)^2}{(2m)!}\Phi_{\cA}^{\hor}(\iota_{\what{x}}\wtilde{p})
  \]
\end{prop}
\begin{proof}
  By Theorem~\ref{thm:iotaTtilde}, for~$\wtilde{p} \in \wtilde{\mathfrak{T}}^{m+1}(\fg,K)$ and~$x \in \fp$ we have
  \[
    \iota_x \wtilde{\mathsf{t}}_{\cA}(\wtilde{p}) = \frac{(m!)^2}{(2m)!} \Phi_{\cA}(\iota_{\what{x}} \wtilde{p}).
  \]
  Now, using Corollary~\ref{cor:abs_rel_tr} and noting that the maps~$\iota_x$, $x \in \fp$ commute with
  $\Res_{\fp}^{\fg}$ on~$\cA$, for $\cA=\twedge\fg^\ast$ and $\cB=(\twedge\fp^\ast)^K$ we have
  \begin{align*}
  \iota_x \wtilde{\mathsf{t}}_{\cB}(\wtilde{p}) &{}= \iota_x \left( -4^m \Res_{\fp}^{\fg} \circ \wtilde{\mathsf{t}}_{\cA}(\wtilde{p}) \right) 
    = -4^m \Res_{\fp}^{\fg} \circ  \iota_x\circ \wtilde{\mathsf{t}}_{\cA}(\wtilde{p})  \\
    &{}= -4^m \Res_{\fp}^{\fg} \left( \frac{(m!)^2}{(2m)!} \Phi_{\cA}(\iota_{\what{x}}\wtilde{p}) \right) \\
    \intertext{(Since for $\cA = \twedge\fg^\ast$ we have that
    $\Res^{\fg}_{\fp} \circ \wtilde{\pi}_{\cA} = \wtilde{\pi}_{\cA} \circ \Res^{\fg}_{\fp}$,
    using~\eqref{eq:PhiPeqResPhi} and the definition~\eqref{eq:PhiAhor} of~$\Phi_{\cA}^{\hor}$, we get
    that $\Phi_{\cA}^{\hor}=\Res^\frg_\frp\circ\Phi_{\cA}$)}
    &{}= - \frac{4^m(m!)^2}{(2m)!} \Phi_{\cA}^{\hor}(\iota_{\what{x}}\wtilde{p}).
  \end{align*}
  For the other case $\cA=\Cl(\fg)$, $\cB = \Cl(\fp)^K$ the claim is proved analogously.
\end{proof}

\subsection{Transgression in the relative commutative  Weil algebra} 
Define the relative commutative Weil algebra
\[
  W(\fg,K) = W(\fg)_{K-\bas} = \left( S(\fg^\ast) \otimes \twedge \fp^\ast \right)^K.
\]

In a similar manner to the case of~$\tW(\fg,K)$ one can define a $K$-homotopy operator for~$W(\fg,K)$ between the identity map $\id_{W(\fg,K)}$ and $j\circ \pr_{\fk}$, which we will denote
by abuse of notation by~$h_\rel$.

The \emph{relative transgression space} in~$W(\fg,K)$ is
\begin{align*}
  \mathfrak{T}(\fg,K) = \left\{ p\in W^{+,0}(\fg,K) \mid \exists C_p\in W(\fg,K): p = \dd_{W}C_p
  \right\}.
\end{align*}
For $p\in\frT(\fg,K)$ an~element~$C_p\in W(\fg,K)$ such that $\dd_W C_p = p$ and $\pr_\fk C_p = 0$ is called
a~\emph{relative cochain of transgression} for~$p$ in~$W(\fg,K)$. 
In Lemma~\ref{lem:TgKprops} below we show that for any element of the relative transgression
space~$\mathfrak{T}(\fg,K)$ there is a~relative cochain of transgression.
In other words,
\begin{align*}
  \mathfrak{T}(\fg,K) = \left\{ p\in W^{+,0}(\fg,K) \mid \exists C_p\in W(\fg,K):
  \text{$p = \dd_{W}C_p$ and $\pr_\fk C_p = 0$}
  \right\}.
\end{align*}

Let $\cA = \twedge\fg^\ast$  and set $\cB = \cA_{K-\bas} = (\twedge\fp^\ast)^K$.
Then the relative transgression map $\mathsf{t}_{\cB}\colon \mathfrak{T}(\fg,K)\to \cB$ for the relative
commutative Weil algebra is defined as follows.
\begin{Def}
  For $p\in \mathfrak{T}(\fg,K)$ let $C_p\in W(\fg,K)$ be a~relative cochain of transgression for~$p$.
  The linear map $\mathsf{t}_{\cB}\colon \mathfrak{T}(\fg,K)\to \cB$ defined by $\mathsf{t}_{\cB}(p) = \pi_{\cB}(C_p)$ is
  called the \emph{relative transgression map} in $W(\fg,K)$ with values in~$\cB$. 
\end{Def}

The following lemma is a commutative analogue of Lemma~\ref{lem:tWrelHomotopy} and Remark~\ref{rem:homotopyop}.
By abuse of notation, we use the same symbols for analogues of the maps $j$, $\pr_{\fk}$ and $h_{\rel}$.
We define the $K$-differential algebra morphism $\pr_{\fk}: W(\fg) \to W(\fk)$ analogously to $\pr_{\fk}:\tW(\frg) \to \tW(\frk)$,
see~\eqref{eq:DefPrK}.
Note that $W(\fg)$ is a~locally free $K$-differential algebra with the connection given by $\mu\mapsto\mu$ for
$\mu\in\fk^\ast$.
Let $j\colon W(\fk) \to W(\fg)$ be the corresponding characteristic homormorphism of
$K$-differential algebras defined analogously to~\eqref{eq:j}.

\begin{lem} \label{lem:CommWrelHomotopy}
    Let~$\pr_{\fk}: W(\fg) \to W(\fk)$ be as above. Then~$\pr_{\fk}$ is a $K$-homotopy equivalence, with inverse $j$. In more detail, $\pr_{\fk} \circ j = \id_{W(\fk)}$, while a $K$-homotopy between~$j \circ \pr_{\fk}$ and~$\id_{W(\fg)}$ is given by
    \[
    \psi_{\rel}: W(\fg) \to \bbC[t, dt] \otimes W(\fg),\qquad \psi_{\rel}(\mu) = (1-t) \otimes \mu + t \otimes j \circ \pr_{\fk}(\mu),\quad \mu \in \fg^{\ast},
    \]
    extended as a morphism of differential algebras using the universal property of~$\tW(\fg)$
    viewed as a~commutative Koszul algebra; see~\cite[Proposition~6.3 in~\S6.4]{MeinrenkenBook}.   
    Define \begin{equation}\label{eq:commhrel}
  h_{\rel} = (J \otimes \id)\circ \psi_{\rel}.
\end{equation}  then 
\begin{equation}   \label{eq:commhreld}
  \dd_{W} \circ h_{\rel} + h_{\rel} \circ \dd_{W} = \id_{W(\frg)} - j \circ \pr_{\fk}.
  \end{equation}  
\end{lem}

\begin{proof}
  The statement of the lemma and the proof are exactly the same as the proof of Lemma~\ref{lem:tWrelHomotopy} and Remark~\ref{rem:homotopyop}. 
\end{proof}
  
\begin{lem}\label{lem:TgKprops}
  We have that
  \begin{enumerate}
  \item for any element of~$\frT(\fg,K)$ there is a~relative cochain of transgression\label{item:TgKrelcohainExists},
  \item $\left\{ w \in W(\fg,K)\cap\ker\dd_W \mid \pr_\fk w = 0\right\} \subseteq \dd_W(W(\fg,K))$,\label{item:relKerDwcapKerPr}
  \item $\frT(\fg,K) = \frT(\fg)\cap \ker\pr_\fk$
    (which is equal to $S^+(\fg^\ast)^\fg \cap \ker\pr_\fk$).\label{item:TgKinTgKerPr}
  \end{enumerate}
\end{lem}
\begin{proof}
  (\ref{item:TgKrelcohainExists})
  We will argue that if there is a cochain of transgression $C \in W(\frg)$ for $p \in \frT(\fg,K)$,
  then we can find a relative cochain of transgression $C_p \in W(\frg,K)$ for~$p$.
  Since $h_{\rel}$ is a~$K$-homotopy operator, the $K$-differential algebra
  map~$\pr_\fk\colon W(\fg)\to  W(\fk)$ induces an isomorphism in $K$-basic cohomology:
  \begin{equation}\label{eq:prkIsoCohom}
    \pr_\fk \colon H(W(\fg,K)) \to H(W(\fk)_{K-\bas}) = S(\fk^\ast)^K.
  \end{equation}
  Suppose now that $p\in \frT(\fg,K)$, then $[p] \in H(W(\fg,K))$ is zero.
  Thus since the differential on $W(\fk)_{K-\bas} = S(\fk^\ast)^K$ is zero, we have
  \[
    0 = \pr_\fk[p] = [\pr_\fk p] = \pr_\fk p
  \]
    in $H(W(\fk)_{K-\bas}) = S(\fk^\ast)^K$.
  Set $C_p = C - j\circ \pr_\fk(C) \in W(\fg,K)$, then we have that
  \begin{gather*}
    \dd_W C_p = \dd_W C  - j\circ \pr_\fk(\dd_W C) =  p - j\circ\pr_\fk(p) = p,\\
    \pr_\fk C_p = \pr_\fk C - \pr_\fk\circ j \circ \pr_\fk(C) = \pr_\fk(C) - \pr_\fk(C) = 0.
  \end{gather*}
  Here we used the fact that $\pr_\fk\circ j=\id_{W(\fk)}$, which implies the second line above. 
  We conclude that $C_p$ is a~relative cochain of transgression for~$p$.

  (\ref{item:relKerDwcapKerPr})
  Assume that
  \[
    w \in W(\fg,K)\cap\ker\dd_W\qquad\text{and}\qquad\pr_\fk w = 0.
  \]
  Using~\eqref{eq:prkIsoCohom}, for such~$w$ we get that $[w]\in H(W(\fg,K))$ corresponds to $[\pr_\fk w]=0$
  under the isomorphism $\pr_\fk\colon H(W(\fg,K))\to H(W(\fk)_{K-\bas})$. Thus $[w]=0$, i.e., $w\in \dd_W(W(\fg,K))$.  

  (\ref{item:TgKinTgKerPr})
  Note that for any $x\in\fg$ and $p\in {\frT}(\fg,K)$, since $\dd_{W}{p} = 0$ and
  since $p\in S^+(\fg^\ast)^K$ implies $\iota_xp=0$, we have
  \[
    L_x {p} = \dd_{W}\circ\iota_x(p) + \iota_x\circ \dd_{W}({p}) = 0.
  \]
  Hence, $p\in W^{+,0}(\fg)^\fg$. Thus $\frT(\fg,K) \subseteq \frT(\fg) = S^+(\fg^\ast)^\fg$.
  It follows from~(\ref{item:TgKrelcohainExists}) that for $p\in\frT(\fg,K)$ there is a $C_p\in W(\fg,K)$ such that
  $p = \dd_W(C_p)$ and $\pr_\fk C_p =0$. Since $\pr_\fk$ is a~morphism of differential algebras, we get that
  $\pr_\fk p = 0$.
  So $\frT(\fg,K) \subseteq S^+(\fg^\ast)^\fg \cap \ker \pr_\fk$.

  Now assume that $p\in S^+(\fg^\ast)^\fg\cap\ker\pr_\fk$, then, using~\eqref{eq:commhreld}, we get that
  \[
    \dd_W\circ h_\rel(p) = p.
  \]
  Therefore, since $h_\rel(p)\in W(\fg,K)$, we get that $p\in \frT(\fg,K)$.
\end{proof}

\begin{lem}\label{lem:trelWellDef}
  The relative transgression map  $\mathsf{t}_{\cB}\colon \mathfrak{T}(\fg,K)\to \cB$ is well-defined, i.e.,
  does not depend on a~choice of a~relative cochain of transgression.  
\end{lem}
\begin{proof}
  Since $(\fg,\fk)$ is a~symmetric pair, the restriction of the differential from~$\cA =
  \twedge\fg^\ast$ to~$\cB = (\twedge\fp^\ast)^K$ is zero, see~\cite[\S9.11]{Onishchik1994} and~\cite[Remark~7.7]{MeinrenkenBook}.
  Therefore, $\cB \cap \im \dd_{\cB} = \{0\}$.
  Hence, $\pi_{\cB}$ vanishes on the space~$\dd_W(W(\fg,K))$.
  Let $C,C^\prime\in W(\fg,K)$ be two relative cochains of transgression for~$p$,
  then $\dd_W(C - C^\prime) = p - p = 0$ and $\pr_\fk(C-C^\prime)=0$. 
  Therefore, by Lemma~\ref{lem:TgKprops}(\ref{item:relKerDwcapKerPr}) we get that $C-C^\prime \in \dd_W(W(\fg,K))$.
  Hence, $\pi_{\cB}(C - C^\prime) = 0$. 
\end{proof}

\begin{lem}\label{lem:tResCommutative}
  For $p\in \mathfrak{T}^{m+1}(\fg,K)$ we have that
  \[
    \mathsf{t}_{(\wedge\fp^\ast)^K}(p) = - 4^m\Res^\fg_\fp\mathsf{t}_{\wedge\fg^\ast}(p)
  \]
\end{lem}

\begin{proof}

  Since $W(\fg)$ admits a~connection, we can use the universal property of~$\tW(\fg)$ to define the
  characteristic homomorphism $c\colon \tW(\fg) \to W(\fg)$ which is the quotient map in this
  case, see~\cite[\S6.11]{MeinrenkenBook}. Then we have that 
  \[    
    \dd_W(c(h_\rel( \what \sym(p)))) = c(\dd_{\tW}(h_\rel(\what{\sym}(p)))) =  c(\what{\sym}(p) -h_\rel(\dd_{\tW}(\what{\sym}(p))))= p,
  \]
  since $c(\what{\sym}(p)) = c( \sym_W(p))  = p$, $c$ is a~morphism of $\fg$-differential spaces, and by Remark \ref{rem:homotopyop} $c(h_\rel(\dd_{\tW}(\what{\sym}(p))))=0$.
  
  Using Lemma~\ref{lem:prkhrelp}, we see that
  \[
    \pr_\fk(c(h_\rel(\what{\sym}(p))))= c(\pr_\fk(h_\rel(\what{\sym}W(p)))) = 0.
  \]
  Hence $c(h_\rel(\what{\sym}(p)))$ is a~relative cochain of transgression for~$p$.

  It is easy to check on generators of $\tW(\fg)$ that $\wtilde{\pi}_{\wedge\fg^\ast} = \pi_{\wedge\fg^\ast}\circ c$.
  By passing to $K$-basic subalgebras, we get that $\wtilde{\pi}_{(\wedge\fp^\ast)^K} = \pi_{(\wedge\fp^\ast)^K}\circ c$.
  It follows that
  \[
    \mathsf{t}_{(\wedge\fp^\ast)^K}(p) = \pi_{(\wedge\fp^\ast)^K}(c(h_{\rel}(\what{\sym}(p)))) 
    = \wtilde{\pi}_{(\wedge\fp^\ast)^K}(h_{\rel}(\what{\sym}(p)))  = \wtilde{\mathsf{t}}_{(\wedge\fp^\ast)^K}(\what{\sym}(p)).
  \]
  and similarly
  \[
    \mathsf{t}_{\wedge\fg^\ast}(p) = \wtilde{\mathsf{t}}_{\wedge\fg^\ast}(\what{\sym}(p)).
  \]
  The claim now follows from Corollary~\ref{cor:abs_rel_tr}.
\end{proof}

\subsection{Transgression in the  relative quantum  Weil algebra} 
Define the relative quantum Weil algebra by
\[
  \cW(\fg,K) = \cW(\fg)_{K-\bas} = (U(\fg)\otimes\Cl(\fp))^K.
\]
From~\cite[\S7.4]{MeinrenkenBook} we have that $H(W(\fg,K)) \simeq S(\fk)^K$
and $H(\cW(\fg,K)) \simeq S(\fk)^K$.

\begin{Def}
The \emph{relative transgression space} in~$\cW(\fg,K)$ is 
\[
\ul{\mathfrak{T}}(\fg,K) = \left\{ \ul{p}\in \cW^{+,0}(\fg,K) \mid \exists \ul{C}_{\ul{p}}\in \cW(\fg,K): \ul{p} = \dd_{\cW}\ul{C}_{\ul{p}}\right\}.
\]
For $\ul{p}\in\ul{\frT}(\fg,K)$ an~element~$\ul{C}_{\ul{p}}\in\cW(\fg,K)$ such that $\ul{p} = \dd_{\cW}\ul{C}_{\ul{p}}$ 
and $\pr_\fk(q_{W}^{-1}( \ul{C}_{\ul{p}})) = 0$ is called a~\emph{relative cochain of transgression} for~$\ul{p}$.
\end{Def}

Note that since $q_W\colon W(\fg)\to \cW(\fg)$ is an~isomorphism of $\fg$-differential spaces, it can be restricted
to $K$-basic subalgebras. The restriction $q_W\colon W(\fg,K)\to\cW(\fg,K)$ is an isomorphism of differential
spaces.
Therefore, we have that
\begin{equation}\label{eq:resCandT}
  \ul{\mathfrak{T}}(\fg,K)  =q_W(\mathfrak{T}(\fg,K)).
\end{equation}
Therefore, there is a~relative cochain of transgression for every element of $\ul{\mathfrak{T}}(\fg,K)$.

\begin{Def}
For $\ul{p}\in \ul{\mathfrak{T}}(\fg,K)$ let $\ul{C}_{\ul{p}}\in\cW(\fg,K)$ be a~corresponding relative cochain of
transgression.
The linear map $\ul{\mathsf{t}}_{\Cl(\fp)^K}(\ul{p})\colon \ul{\mathfrak{T}}(\fg,K) \to \Cl(\fp)^K$ defined by
$\ul{\mathsf{t}}_{\Cl(\fp)^K}(\ul{p}) = \ul{\pi}_{\Cl(\fp)^K}(\ul{C}_{\ul{p}})$ is called the \emph{relative
transgression map} in~$\cW(\fg,K)$.
\end{Def}
The relative transgression map in~$\cW(\fg,K)$ is well-defined by similar arguments as in
Subsection~\ref{ssec:cWg}, using Lemma~\ref{lem:trelWellDef}.

The main result of this subsection is the following theorem:

\begin{thm}\label{thm:IotaTClP}  
  For any $x\in\fp$ and $p\in\ul{\mathfrak{T}}(\fg,K)$,  $\iota_x\ul{\mathsf{t}}_{\Cl(\fp)^K}(p) \in \im \alpha_\fp$.  
\end{thm}

Using the universal property of~$\tW(\fg)$ we show that the images of the relative transgressions in $\tW(\fg,K)$ and
$\cW(\fg,K)$ with values in~$\Cl(\fp)^K$ coincide (Lemma~\ref{lem:imTsameRel}),
and so do the images of relative quantum moment maps $\Phi^{\hor}_{\Cl(\fg)}$ and~$\alpha_\fp$ (Lemma~\ref{lem:imalphaisimphiRel}).
Thus Theorem~\ref{thm:IotaTClP} is a~direct consequence of the
universal property of~$\tW(\fg)$ and~Proposition~\ref{prop:contr_trans_formula}.

The proof of the following lemma is postponed to after the proof of the theorem:

\begin{lem}\label{lem:imTsameRel}
  For $p\in \frT^{m+1}(\fg,K)$,
  \[   \wtilde{\mathsf{t}}_{\Cl(\fp)^K}(\what{\sym}(p)) = \ul{\mathsf{t}}_{\Cl(\fp)^K}(\ul{c}(\sym_W(p))).
  \]
\end{lem}
\begin{proof}[Proof of Theorem~\ref{thm:IotaTClP}]
Using Lemma~\ref{lem:imTsameRel} and 
the same argument as in the proof of Theorem~\ref{thm:iotaxImAlphaG}, we conclude that $\im \underline{\mathsf{t}}_{\Cl(\fp)^K} = \im\wtilde{\mathsf{t}}_{\Cl(\fp)^K}$.

We now see from Proposition~\ref{prop:contr_trans_formula} that $\iota_x (\im \wtilde{\mathsf{t}}_{\Cl(\fp)^K}) \subset \Phi_{\Cl(\fg)}^{\hor}\left( \tW^{\bullet,0}(\fk) \right)$, so $\iota_x (\im \ul{\mathsf{t}}_{\Cl(\fp)}) \subset \Phi_{\Cl(\fg)}^{\hor}\left( \tW^{\bullet,0}(\fk) \right)$.

This is equal to~$\im \alpha_{\fp}$ by the following lemma.
\end{proof}

\begin{lem}\label{lem:imalphaisimphiRel}
  For the algebra morphism $\Phi_{\Cl(\fg)}^{\hor} \colon \tW(\fg)\to \Cl(\fp)$ defined in~\eqref{eq:PhiAhor}, we have~$\Phi_{\Cl(\fg)}^{\hor}\left( \tW^{\bullet,0}(\fk) \right) = \alpha_\fp( U(\fk))$.
\end{lem}
\begin{proof}
By abuse of notation let~$\ul{c}$ be the characteristic homomorphism~$\tW(\fk) \to \cW(\fk)$.
Since both~$\Phi_{\Cl(\fg)}^{\hor}$ and~$\alpha_{\fp}$ are algebra homomorphisms
and the map $\ul{c}$ is surjective,
  it suffices to show that~$\Phi_{\Cl(\fg)}^{\hor}(\what{\mu})$ and $\alpha_{\fp}(\ul{c}(\what{\mu}))$ are equal for~$\mu \in \fk^{\ast}$.
  Take a generator~$\what{f_i}$ for~$f_i \in \fk^*$, we have  
\begin{align*}
  \Phi_{\Cl(\fg)}^{\hor}(\what{f}_i)&{}= \wtilde{\pi}_{\Cl(\fp)}(\delta_{\fp}(f_i)) = \wtilde{\pi}_{\Cl(\fp)} \left( - \frac{1}{2} \sum_{a,b}^{(\fp)} c_{a,b}^i \, f_a \otimes f_b \right) = - \frac{1}{2}\sum_{a,b}^{(\fp)} c_{a,b}^i \, f_a \, f_b, \\
\intertext{where we used the formula~\eqref{eq:deltaP} for~$\delta_\fp$ and
  formulas~\eqref{eq:DeltaGasDeltaKP} and~\eqref{eq:PhiAhor}. 
  From the $\fp$-part of~\eqref{eq:alphafi} we get that}
  \alpha_{\fp}(\ul{c}(\what{f}_i))
 &{}= \alpha_{\fp}(2\what f_i)= - \frac{1}{2} \sum_{a,b}^{(\fp)} c_{b,a}^i\,  f_b\, f_a.
\end{align*}
This proves the claim. 
\end{proof}

\noindent{\it Proof of Lemma \ref{lem:imTsameRel}.} 
  First note that it is easy to check that $\wtilde{\pi}_{\Cl(\fg)} = \ul{\pi}_{\Cl(\fg)}\circ \ul{c}$. 
  Hence $\wtilde{\pi}_{\Cl(\fp)^K} = \ul{\pi}_{\Cl(\fp)^K}\circ \ul{c}$.
  By identical arguments to the proof of Lemma \ref{lem:tResCommutative},
  $\ul{c}(h_{\rel}(\what{\sym}(p))$
  is a~relative cochain of transgression for~$\ul{c}(\sym_{W}(p))$. We have
  \begin{align*}
  \wtilde{\mathsf{t}}_{\Cl(\fp)^K}(\what{\sym}(p)) ={}&  \wtilde{\pi}_{\Cl(\fp)^K}(h_{\rel}(\what{\sym}(p)))
    = \ul{\pi}_{\Cl(\fp)^K}(\ul{c}(h_{\rel}(\what{\sym}(p))\\
    {}={}&  \ul{\mathsf{t}}_{\Cl(\fp)^K}(\ul{c}(\sym_W(p))),
  \end{align*}
  which proves the claim. 
\qed

\smallskip

For later use, we also prove

\begin{lem}\label{lem:quanttr_rel}
  For any $p\in \mathfrak{T}^{m+1}(\fg,K)$, $q\circ \mathsf{t}_{(\wedge\fp)^K}(p) = \ul{\mathsf{t}}_{\Cl(\fp)^K}(q_{W}(p))$. 
\end{lem}
\begin{proof}
  Let $C_p\in W(\fg,K)$ be a~relative cochain of transgression for~$p$.
  Then by~\eqref{eq:resCandT} and the fact that $q_W$ is an~isomorphism of differential spaces,
  $q_W(C_p)\in\cW(\fg,K)$ is a~relative cochain of transgression for~$q_W(p)$. 
  Using~\eqref{q and pi} restricted to $(\twedge\fp^\ast)^K$ and $\Cl(\fp)^K$, we have that
  \[
    q\circ\pi_{(\wedge\fp^\ast)^K}(C_p) = \ul{\pi}_{\Cl(\fp)^K}(q_W(C_p)).
  \]
  Which proves the claim.
\end{proof}

To end this subsection, we note that Theorem \ref{thm:IotaTClP} also shows that the image of the transgression map~$\ul{\mathsf{t}}_{\Cl(\fp)^K}$ consists of
primitive 
elements in the sense of~\cite[\S5]{Rashevskii1969}.
This fact motivates our definition of the primitive invariants in~$\Cl(\fp)$, see Section~\ref{sec:Primitives}.

\section{Primitive invariants and relative transgression theorem for
  primary and almost primary types}
\label{sec:Primitives}
In this section we define the space of primitive $K$-invariants in~$\twedge\fp^\ast$ and prove the
relative analogue of the celebrated Transgression Theorem~\ref{thm:transgression}.

We specify $(\fg,\fk,K)$ to be in the following list:
\begin{equation}\label{eq:TheList} 
\begin{minipage}[l]{10cm}
\begin{enumerate}
\item $(\fsl(2n+1),\fo(2n+1),SO(2n+1))$,
\item $(\fsl(2n), \fo(2n), O(2n))$,
\item $(\fsl(2n), \fsp(2n), Sp(2n))$,
\item $(\fe(6), \ff(4), F(4))$.
\end{enumerate}
\end{minipage}
\end{equation}
\begin{lem}\label{lem:SplittingForPrk}
  For $(\fg,\fk,K)$ as in the list~\eqref{eq:TheList} we have that
  \[
    \pr_\fk\colon S(\fg^\ast)^\fg \to S(\fk^\ast)^K
  \]
  is surjective and admits a~splitting map
  \[
    s\colon S(\fk^\ast)^K\to S(\fg^\ast)^\fg
  \]
  which is an~algebra map.
\end{lem}

\begin{proof}
We first tackle items (a),(c),(d) on our list~\eqref{eq:TheList}, noting that in these settings $K$ is connected and $S(\frk^*)^K =S(\frk^*)^\frk$. We have the following commutative diagram of Harish-Chandra isomorphisms

 \[   
   \begin{tikzcd}
     S(\frg^\ast)^\frg \arrow{r}{\pr_\frk} \arrow{d}{\HC} & S(\frk^\ast)^\frk \arrow{d}{\HC} \\
     S(\frh^\ast)^{W_{\frg}} \arrow{r}{\Res^\frh_\frt} \arrow{d}{\cong} &  S(\frk^*)^{W_\frk} \arrow{d}{\cong}.\\
     S(P(W_\frg)) \arrow{r}{\Res^\frh_\frt} &S(P(W_\frk))
   \end{tikzcd}
 \]
 Here $P(W_\fg)$ and $P(W_\fk)$ are the Dynkin subspaces of $S(\fh^\ast)^{W_\fg}$, respectively
 $S(\ft^\ast)^{W_\fk}$, i.e., the orthogonal complements to the squares of the augmentation ideals,
 see~\cite[p. 264]{MeinrenkenBook}.
 
 In (a) (respectively (c)), let $\frh$ be traceless diagonal $(2n+1)\times (2n+1)$ (respectively $2n \times
 2n$) matrices, and let $x_i: \frh \to \bbC$ be the functional that picks out the $i$th entry. The space
 $P(W_\frg)$ can be taken to be the span of the degree $2$ to $2n+1$ (or $2n$) power sum symmetric polynomials
 in~$x_i$. Let $\frt$ be diagonal matrices which are antisymmetric under reflection along the antidiagonal,
 let $\{y_i \in \frt^*: i \in 1,\cdots,n\}$ pick out the $i$th entry.
 Then we can take $P(W_\frk)$ to be the span of
 the first $n$ power sum polynomials in $y_i^2$.
 For $i =1 ,\cdots, n$, $\Res^\frh_\frt(x_i) = y_i$ and
 $\Res^\frh_\frt(x_{n+i}) = -y_i$, hence the restriction to~$\ft$ of an odd power sum polynomial in $x_i$ is
 zero and the restriction of an even power sum polynomial is twice the power sum polynomial in $\{y_i^2\}$. Thus
 $\Res^\frh_\frt$ is surjective and we construct a~splitting by lifting the degree~$j$ power sum polynomial in~$y_i^2$ to degree
 $2j$ power sum polynomial in~$x_i$ for $j=1,\ldots,n$.
 
For (d) one could check by computer or appeal to \cite[Theorem~1 on p.~216]{Onishchik1994} which states that when $(\frg,\frk)$ is primary then $\pr_\frk: S(\frg^*)^\frg \to S(\frk^*)^\frk$ is surjective and a splitting can be constructed using any preimage of the generators for $S(\frh^*)^\frk$, this will be an algebra morphism since $S(\frh^*)^\frk$ is a free commutative algebra.

 When the triple $(\frg,\frk,K)$ is in the category (b) then $P(W_\frg)$ can again be taken to be the first
 $2n$ power sum polynomials in $x_i$, and again $\frt$ are diagonal matrices which are antisymmetric under
 reflection along the antidiagonal. Note that we can take $P(W_\frk)$ to be the span of the
 first $n-1$ power sum  polynomials in  $y_i^2$ and the Pfaffian $\prod_{i=1}^n y_i$. 
 In this case $\pr_\frk$ does not surject onto
 $S(\frh^*)^\frk$. However, the Pfaffian is not $K$ invariant; the algebra $S(\frk^*)^K$ is congruent under
 the Harish-Chandra isomorphism to $S(\frt^\ast)^{W_K}$, with $W_K \cong W_{B_n}$ (see~\cite[\S3.8]{CNP}).
 This algebra is equal to the polynomial algebra (of rank $\dim \frt$) generated by the first $n$ power sum polynomials in
 $y_i^2$. Using the reasoning for (c), this is clearly the image of $\Res^\frh_\frt:S(\frh^\ast)^{W_\frg} \to
 S(\frt^\ast)^{W_\frk}$. Again a splitting lifts the degree $j$ symmetric power sum polynomials in $y_i^2$ to
 degree $2j$ symmetric power sum polynomials in $x_i$.
\end{proof}

\begin{lem}\label{lem:ExactSeq}
  Consider $(\fg,\fk,K)$ as in the list~\eqref{eq:TheList}, then the sequence
  \[
    \begin{tikzcd}
      0\arrow[r] & {\frT}(\fg,K) \arrow[r, hookrightarrow] & S^+(\fg^\ast)^{\fg}\arrow[r,"\pr_\fk"] & S^{+}(\fk^\ast)^{K}\arrow[r] & 0
    \end{tikzcd}
  \]
  is exact and split.
\end{lem}
\begin{proof}
  By Lemma~\ref{lem:TgKprops}(\ref{item:TgKinTgKerPr}) we have that
  $\frT(\fg,K) = S^+(\fg^\ast)^{\fg} \cap \ker\fp_\fk$.
  It was shown in Lemma~\ref{lem:SplittingForPrk} that the map $\pr_{\fk}\colon S^+(\fg^\ast)^\fg\to S^+(\fk^\ast)^K$ is
  surjective and admits a~splitting map.
\end{proof}

\begin{Def}\label{def:PwedgeP}
The $\mathbb{Z}$-graded subspace $P_\wedge(\fp) := \Res^\fg_\fp P_\wedge(\fg)\subset(\twedge\fp^\ast)^\fk$ is called
the space of \emph{primitive invariants} in~$\twedge\fp^\ast$. 
\end{Def}

From the explicit form of primitives $P_\wedge(\fp)$ for $(\fsl(2n),\fso(2n))$ given in \cite{Dolce2016}, it
is easy to see that the primitives are in fact not only $\fk$-invariant but also $K$-invariant.
(They are given as traces, hence are invariant under conjugation by~$K$.)

It follows from~\cite[Theorem~12]{HochschildSerre1953}, see also the discussion in~\cite[\S4.2]{Dolce2015}
and~\cite{Dolce2016},
that for triples $(\fg,\fk,K)$ from the list~\eqref{eq:TheList} we have
\begin{equation}\label{eq:dimPwedgeP}
  \dim P_\wedge(\fp) = \rank\fg - \rank\fk = \dim\fa.
\end{equation}

We claim that the space $\mathfrak{T}(\fg,K)$ admits a~structure of a $S(\fg^\ast)^\fg$-module induced by the multiplication in
the symmetric algebra $S(\fg^\ast)$. Namely, let $p\in \mathfrak{T}(\fg,K)$ and $x\in S(\fg^\ast)^\fg$.
Let  $C_p\in W(\fg,K)$ be a~relative cochain of transgression for $p$.
Then $\dd_{W}(x C_p) = xp$. Since $x C_p\in W(\fg,K)$, $x C_p$ is a~relative cochain of transgression for $xp$.
Clearly, if $x\in S^+(\fg^\ast)^\fg$ we have that
\begin{equation}\label{eq:t2vanishes}
  \mathsf{t}_{\cB}(xp) = \pi_{\cA}(x C_p) = \pi_{\cA}(x)\pi_{\cA}(C_p) = 0
\end{equation}
since $\pi_{\cA}(x)$ vanishes by Lemma~\ref{lem:Aflat}\ref{item:Aflat1}.
It also follows that $\frT(\fg,K)$ is an ideal in~$S^+(\fg^\ast)^\fg$.

The following theorem is an~extended version of~\cite[Theorem~4.3.11]{GrizlyPhD}.

\begin{thm}[Relative transgression theorem]\label{thm:relTransgression}
  Let $(\fg,\fk,K)$ be as in the list~\eqref{eq:TheList},
  then the transgression map~$\mathsf{t}_{(\wedge\fp^\ast)^K}$ satisfies
  \[
    \ker \mathsf{t}_{(\wedge\fp^\ast)^K} = S^+(\fg^\ast)^\fg\cdot \mathfrak{T}(\fg,K),\qquad
    \im \mathsf{t}_{(\wedge\fp^\ast)^K} = P_\wedge(\fp).
  \]
  In particular, the map
  \begin{equation}\label{eq:trelIsIso}
    \mathsf{t}_{(\wedge\fp^\ast)^K}\colon\faktor{\mathfrak{T}(\fg,K)}{S^+({\fg^\ast})^{\fg}\cdot \mathfrak{T}(\fg,K)} \to P_{\wedge}(\fp).
  \end{equation}
  is an~isomorphism of linear spaces.
\end{thm}

\begin{proof}
  By~\eqref{eq:t2vanishes},
  \(
    S^+(\fg)^{\fg} \cdot \mathfrak{T}(\fg, K) \subseteq \ker \mathsf{t}_{(\wedge\fp^\ast)^K}
  \).
  
  By Lemma~\ref{lem:ExactSeq} the short exact sequence
  \begin{equation}\label{eq:SplitExactSeq}
    \begin{tikzcd}
      0\arrow[r] & \mathfrak{T}(\fg,K) \arrow[r, hookrightarrow] & S^+(\fg^\ast)^{\fg}\arrow[r, "\pr_\fk"] & S^+(\fk^\ast)^{K}\arrow[r] & 0
    \end{tikzcd}
  \end{equation}
  is split.  
  Let $s\colon S^+(\fk^\ast)^K\to S^+(\fg^\ast)^\fg$ be the
  splitting map constructed in Lemma~\ref{lem:SplittingForPrk}. We have
  \begin{equation}\label{eq:SplusSplit}
    S^+(\fg^\ast)^{\fg} = \mathfrak{T}(\fg,K) \oplus s(S^+(\fk^\ast)^{K}) \cong \mathfrak{T}(\fg,K) \oplus S^+(\fk^\ast)^{K}.
  \end{equation}
  In particular, noting that  $\mathfrak{T}(\fg,K) \cdot S^+(\fg^\ast)^{\fg}$ is an ideal in
  $\mathfrak{T}(\fg,K)$ and that $\mathfrak{T}(\fg,K)$ is an ideal in $S^+(\fg^\ast)^{\fg}$,
  we have that 
  \begin{equation}\label{eq:Splus2Split}
    \big(S^+(\fg^\ast)^{\fg}\big)^2 \cong \mathfrak{T}(\fg,K) \cdot S^+(\fg^\ast)^{\fg} \oplus \big(S^+(\fk^\ast)^{K}\big)^2.
  \end{equation}
  
  Since $S(\frg^\ast)^\frg\cong S(\fh^\ast)^{W_\fg}$ is a~polynomial algebra of rank $\dim \frh$ and
  $S(\frk^\ast)^K\cong S(\ft^\ast)^{W_K}$ is a~polynomial  algebra of rank $\dim \frt$,
  it follows  that
  \begin{displaymath}
    \dim \faktor{S^+(\fg)^{\fg}}{\big(S^+(\fg)^{\fg}\big)^2} = \dim \fh \quad \text{and} \quad \dim \faktor{S^+(\fk)^{K}}{\big(S^+(\fk)^{K}\big)^2} = \dim \ft.
  \end{displaymath}
  Therefore, using~\eqref{eq:SplusSplit} and~\eqref{eq:Splus2Split}, we see that
  \[
    \dim \faktor{\mathfrak{T}(\fg,K)}{S^+(\fg)^{\fg}\cdot \mathfrak{T}(\fg,K)} = \dim \fh - \dim \ft = \dim\fa.
  \]

  Let $P(W_\fg)$ and $P(W_K)$ be as in the proof of Lemma~\ref{lem:SplittingForPrk}.
  Set
  \begin{align*}
    P_S(\fg) ={}& \HC^{-1}(P(W_\fg))\cong \faktor{S^+(\fg)^{\fg}}{\big(S^+(\fg)^{\fg}\big)^2}, \\
    P_S(K) ={}& \HC^{-1}(P(W_K)) \cong \faktor{S^+(\fk)^{K}}{\big(S^+(\fk)^{K}\big)^2}.
  \end{align*}
  It follows from Lemma~\ref{lem:SplittingForPrk}
  that the image of~$P_S(\fg)$ under $\pr_\fk$ is $P_S(K)$ and the splitting map~$s$ sends $P_S(K)$ to~$P_S(\fg)$.

  Furthermore, let
  \[
    P_S(\fg,\fk) := P_S(\fg)\cap \frT(\fg,K).
  \]
  Using the split exact sequence~\eqref{eq:SplitExactSeq}, we get 
  \[
    P_S(\fg) = P_S(\fp) \oplus s(P_S(K)) \cong P_S(\fp) \oplus P_S(K)
  \]
  From the discussion above we get
  \[
    P_S(\fg,\fk) \cong \faktor{\mathfrak{T}(\fg,K)}{S^+(\fg)^{\fg}\cdot \mathfrak{T}(\fg,K)}.
  \]
  So the following diagram is a~split exact sequence.
  \begin{equation}\label{eq:PrimSeq}
    \begin{tikzcd}
      0\arrow[r] & P_S(\fg,\fk) \arrow[r, hookrightarrow]& P_S(\fg)\arrow[r, "\pr_\fk",two heads] & P_S(K)\arrow[r]\arrow[l,"s", bend left]& 0
    \end{tikzcd}
  \end{equation}
  
  Recall that  the involution $\theta$ fixes~$\fk$ and acts by $-1$ on~$\fp$.
  It is easy to check that the transgression
  map~$\mathsf{t}_{\wedge\fg^\ast}$ commutes with~$\theta$.

  We claim that $P_S(\fg)$ is stable under $\theta$. Since we can choose a~$\theta$-stable triangular
  decomposition $\fg = \fn^+\oplus\fh\oplus\fn^{-}$, the Harish--\/Chandra isomorphism can be assumed to be
  $\theta$-equivariant.
  So it is enough to show that the space $P(W_\fg)$ of Lemma~\ref{lem:SplittingForPrk} is $\theta$-stable.
  
  For $(\fg,\fk)=(\fsl(n),\fso(n))$ the action of $\theta$ on the space of diagonal matrices~$\fh\subset\fg$
  is given by minus the transpose with respect to the antidiagonal.
  (Here we are taking the realisation of~$\fk=\fo(n)$ to be matrices that are skew-symmetric with respect to the antidiagonal.)
  Hence the even power sums are fixed by~$\theta$ while the odd power sums are in
  the $(-1)$-eigenspace of~$\theta$. The claim follows.

  For $(\fg,\fk)=(\fsl(2n),\fsp(2n))$, the involution $\theta$ is $X\mapsto -JX^t J^{-1}=JX^tJ$, where  
  \(J = \left(\begin{smallmatrix} 0 & A_{n} \\  - A_{n} & 0 \end{smallmatrix}\right)\),
  with $A_{n}$ the $n\times n$ antidiagonal matrix with 1s on the antidiagonal.
  As in the $(\frsl(n),\frso(n))$ case, the action of $\theta$ on the diagonal matrices comprising $\frh$ is given by minus the transpose with respect to the antidiagonal. Hence the same proof works also in this case.

  Assume now that $(\fg,\fk)=(\fe(6),\ff(4))$. Let $\mu_1,\ldots,\mu_l\in\fh^\ast$ be the weights of one of
  $27$-dimensional irreducible $\fe(6)$-modules, then by~\cite[\S11 Theorem~3]{Onishchik1994} the polynomials
  \[
    p_k =  \sum_{j=1}^l \mu_{j}^{m_k+1},\qquad k =1,\ldots,6,
  \]
  where $m_i$ are the exponents of~$\fe(6)$, form a~basis of $P(W_\fg)$. The involution~$\theta$ which
  fixes~$\ff(4)$ in~$\fe(6)$ is induced by the outer automorphism.
  Therefore, $\theta(\mu_j) = - \mu_j$. Hence we get that
  \[
    \theta(p_k) = \begin{cases}
      p_k & \text{if $m_k$ is odd},\\
      - p_k & \text{if $m_k$ is even}.
    \end{cases}
  \]
  Which proves that $P_S(\fg)$ is $\theta$-stable.

  Set
  \[
    P_\wedge(\fg)^{\pm\theta} = \left\{ p\in P_\wedge(\fg) \mid \theta(p) = \pm p\right\},\qquad
    P_S(\fg)^{\pm\theta} = \left\{ p\in P_S(\fg) \mid \theta(p) = \pm p \right\}.
  \]
  We now consider the following diagram with the first row equal to~\eqref{eq:PrimSeq}.
  \begin{equation}\label{eq:PrimitiveDiagram}
    \begin{tikzcd}
      0\arrow[r] & P_S(\fg,\fk) \arrow[r, hookrightarrow]\arrow[d, "\sft_{(\wedge\fp^\ast)^K}"] & P_S(\fg)=P_S(\fg)^{\theta}\oplus P_S(\fg)^{-\theta}\arrow[r, "\pr_\fk"]\arrow[d, "\sft_{\wedge\fg^\ast}"] & P_S(K)\arrow[r]& 0\\
     {} & P_\wedge(\fp) & \arrow[l,"\Res^\fg_\fp"]P_\wedge(\fg)=P_\wedge(\fg)^{\theta}\oplus P_\wedge(\fg)^{-\theta} &
    \end{tikzcd}
  \end{equation}
  Clearly, $\theta$ acts by $-1$ on $\twedge^{\mathrm{odd}}\fp^\ast$.
  Since $P_\wedge(\fg)$ is concentrated in odd degrees, so is $P_\wedge(\fp)$.
  Thus $\theta$ acts by~$-1$ on~$P_\wedge(\fp)$,
  hence $\Res^\fg_\fp\sft_{\wedge\fg^\ast}(P_S(\fg)^\theta) = 0$ and $\Res^\fg_\fp\sft_{\wedge\fg^\ast}(P_S(\fg)^{-\theta}) = P_\wedge(\fp)$.
  Therefore, $\dim P_S(\fg)^{-\theta} \geq \dim P_\wedge(\fp) = \dim\fa$.
  Since $S(\fk)$ is $\theta$-invariant, $P_S(\fg)^{-\theta}\subset\ker\pr_{\fk} = P_S(\fg,\fk)$.
  We have seen that $\dim P_S(\fg,\fk) = \dim\fa$, it follows that $P_S(\fp) = P_S(\fg)^{-\theta}$.  
  Therefore,
  \[
    \Res^\fg_\fp\circ \mathsf{t}_{\wedge\fg^\ast}( P_S(\fg)) = \Res^\fg_\fp\circ \mathsf{t}_{\wedge\fg^\ast}( P_S(\fg)^{-\theta})
  \]
  Since the diagram~\eqref{eq:PrimitiveDiagram}
  commutes (up to multiplication by a non-zero scalar) by Lemma \ref{lem:tResCommutative},
  we have $\mathsf{t}_{(\wedge\fp^\ast)^K}(P_S(\fg,\fk)) = P_\wedge(\fp)$, hence~\eqref{eq:trelIsIso} is an~isomorphism.
\end{proof}

\section{Harish--\/Chandra projections}
In this section we construct relative analogues of the Harish--\/Chandra map $\Cl(\fg)^\fg\to\Cl(\fh)$.
They are attached to elements of~$W^1$ and denoted by~$\HC_w$.
We show that for each $w\in W^1$, $\HC_w\colon\Cl(\fp)^\fk\to\Cl(\fa)$ is a~surjective algebra homomorphim.
We stress that in general it is not a~filtered algebra homomorphism.
However, in Section~\ref{sec:KostantCliffordConjecture} 
we introduce a~filtration on~$\fa$ which makes~$\HC_w$
compatible with filtrations.

Fix a triangular decomposition of~$\fg$
\[
  \fg=  \fn^+\oplus \fh \oplus \fn^-.
\]
By the Poincar\'e-Birkhoff-Witt Theorem, this decomposition of $\frg$ induces a decomposition of~$\Cl(\fg)$
\[
  \Cl(\fg) = \Cl(\fh) \oplus \left( \fn^-\Cl(\fg) + \Cl(\fg)\fn^+\right),
\]
which in turn defines a~projection $\HC_\fg \colon \Cl(\fg)\to \Cl(\fh)$.
By results of Kostant and Bazlov~\cite[Theorem~4.1]{Bazlov2009HC}, see also~\cite[Theorem 11.8]{MeinrenkenBook},
the restriction of~$\HC_\fg$ to~$\Cl(\fg)^\fg$ is an~algebra isomorphism
\[
  \HC_\fg \colon \Cl(\fg)^\fg \to \Cl(\fh).
\]
As it was shown in~\cite{AlekseevMoreau2012}, 
the Harish-Chandra projection becomes a~filtered algebra isomorphism when~$\frh$ is given a~special
filtration, described in Section~\ref{sec:KostantCliffordConjecture}.

Assume now that $\frt=\frh\cap\frk$ is a Cartan subalgebra of $\frk$ and fix a choice~$\Delta^+(\fk,\ft)$ of positive $\ft$-roots in~$\fk$ and a choice of $\Delta^+(\fg,\ft)\supset\Delta^+(\fk,\ft)$. This determines the subset $W^1$ of the Weyl group $W_\frg$ (see the introduction).

Each $w\in W^1$ defines a~choice $\Delta^+_w(\frg,\frt) = w(\Delta^+(\frg,\frt))$
of positive $\ft$-roots in~$\fg$ containing $\Delta^+(\fk,\ft)$.
Consider the corresponding triangular decomposition
\[
  \fg = \fn^+_w \oplus \fh \oplus \fn^-_w.
\]
Set
\[
  \fp^+_w := \fp\cap \fn^+_w,\qquad  \fp^-_w := \fp\cap \fn^-_w,\qquad
\fa:=\fp\cap\fh,
\]
so that
\[  
  \fp = \fp^+_w \oplus \fa \oplus \fp^-_w.
\]
Note that $\fp^+_w $ and $\fp^-_w$ are isotropic and in duality under $B$.

By the Poincar\'e-Birkhoff-Witt Theorem, there is a basis of $\Cl(\frp)$ consisting of monomials
\eq\label{PBW}
f_1^{i_1}\dots f_p^{i_p}h_1^{j_1}\dots h_a^{j_a} e_1^{k_1}\dots e_p^{k_p},
\eeq
with $f_1,\dots,f_p$ a basis of $\frp_w^-$, $h_1,\dots,h_a$ a basis of $\fra$, $e_1,\dots,e_p$ a basis of $\frp^+$, and the exponents $i_r,j_s,k_t$ all equal to 0 or 1.
This induces a decomposition
\[
  \Cl(\fp) = \Cl(\fa) \oplus
  \left(\fp^-_w\Cl(\fp) + \Cl(\fp)\fp^+_w\right)
\]
and consequently a Harish-Chandra projection $\HC_w\colon \Cl(\fp) \to \Cl(\fa)$. 

\begin{lem}\label{hc alg hom}
The restriction of the linear map $\HC_w$ to the $\frt$-invariants $C(\frp)^\frt$ is an algebra homomorphism.  Consequently, $\HC_w\colon \Cl(\fp)^\frk \to \Cl(\fa)$ is also an algebra homomorphism.
\end{lem}
\pf
 If $x\in\Cl(\frp)^\frt$, then $x$ can be written as sum of monomials \eqref{PBW} of $\frt$-weight 0. For the weight to be 0, each such monomial contains a factor in $\frp^+_w$ if and only if it contains a factor in $\frp_w^-$. It follows that
\[
\Cl(\frp)^\frt\cap\left(\fp^-_w\Cl(\fp) + \Cl(\fp)\fp^+_w\right) = \Cl(\frp)^\frt\cap\left(\fp^-_w\Cl(\fp)\right)=\Cl(\frp)^\frt\cap\left(\Cl(\fp)\fp^+_w\right),
\]
and this is clearly an ideal in $\Cl(\frp)^\frt$. Consequently, the projection
$\HC_w:\Cl(\frp)^\frt\to \Cl(\fra)$ is an algebra homomorphism.
\epf

\begin{lem}\label{lem:HCpImAlpha} 
  The Harish-Chandra projection $\HC_w: \Cl(\fp) \to \Cl(\fa)$ maps $\alpha_\fp(U(\frk))$ to $\bbC$.  
\end{lem}
\begin{proof}
In writing the basis \eqref{PBW}, we may assume that $f_i$ and $e_j$ are dual bases (with respect to $B$) of $\frp_w^-$ respectively $\frp_w^+$, such that each $e_i$ is a root vector for the positive $(\frg,\frt)$-root $\beta_i$, while $f_i$  is a root vector for the $(\frg,\frt)$-root $-\beta_i$. We can also assume that the $h_i$ form an orthonormal basis for $\fra$.

Denote by $\rho_\frg$, respectively $\rho_\frk$, the half sum of roots in $\Delta^+(\frg,\frt)$, respectively $\Delta^+(\frk,\frt)$. The half sum of roots in $\Delta^+_w(\frg,\frt)$ is then $w\rho_\frg$.  We claim that
  \eq\label{eq:HConT}
    \HC_w(\alpha_{\fp}(H)) = (w\rho_\frg- \rho_{\fk})(H), \; \quad H \in \ft.
  \eeq
  To see this, we compute
  \begin{align*}
    \alpha_{\fp}(H) &= \frac{1}{4} \Big( \sum_{i=1}^{p} [H,e_i]f_i + \sum_{i=1}^{a} [H,h_i]h_i + \sum_{i=1}^{p} [H,f_i]e_i \Big) \\
    & \text{(since $H \in \ft$ we have $[H,h_i]=0$)} \\
    &= \frac{1}{4} \sum_{i=1}^{p} \beta_i(H)e_i f_i - \frac{1}{4} \sum_{i=1}^{p} \beta_i(H) f_i e_i \\
    &= \frac{1}{4} \sum_{i=1}^{p} \beta_i(H) (e_i f_i - f_i e_i) = \frac{1}{4} \sum_{i=1}^{p} \beta_i(H)(2-2f_i e_i) \\
    &= (w\rho_\frg-\rho_\frk)(H) - \frac{1}{2} \sum_{i=1}^{p} \beta_i(H)f_i e_i.
  \end{align*}
  The last sum vanishes under the map~$\HC_w$ by its definition and~$(w\rho_\frg-\rho_\frk)(H)$ stays unchanged by the same map because it is a constant. This proves \eqref{eq:HConT}.
  
  Next, we claim that~$\HC_w \circ \alpha_\fp$ annihilates elements of $\fn_w^+ \cap  \fk$.
  By the definition of $\HC_w$ it is enough to show that for any root vector $x \in \fn_w^+ \cap \fk$, $\alpha_\fp(x) \in \Cl(\fp)\frp_w^+$.
  The summands in the expression for~$\alpha_\fp(x)$ are of the form
  \[
  [x,e_i]f_i\quad\text{ or}\quad [x,h_i]h_i  \quad\text{ or}\quad [x,f_i]e_i.
  \]
  Of these summands, $[x,f_i]e_i$ is clearly in $\Cl(\fp)\frp_w^+$. Furthermore, $[x,h_i]h_i$ is in
  $\Cl(\fp)\frp_w^+$, since $[x,h_i]\in\fp_w^+$, and it anticommutes with $h_i$. Finally,
  $[x,e_i]\in\frp_w^+$, and it anticommutes with $f_i$ since it is not proportional to $e_i$.
  Hence $[x,e_i]f_i$ is also in $\Cl(\fp)\frp_w^+$.

  One shows similarly that if $y\in \fn_w^- \cap \fk$, then $\alpha_\fp(y)\in\frp_w^-\Cl(\frp)$
  (so $\HC_w(\alpha_\fp(y))=0$).  

  \smallskip
  
  Now let $M$ be a~monomial in~$U(\fk)$ of the form
  \[
    M = F_1 \cdots F_m H_1 \cdots H_k E_1 \cdots E_n,
  \]
  where~$F_1, \ldots, F_m \in \fn_w^- \cap \fk,\, H_1, \ldots, H_k \in \ft,\, E_1, \ldots E_n \in \fn_w^+ \cap \fk$. By the Poincar\'e-Birkhoff-Witt Theorem, such monomials span $U(\frk)$. We claim that $\HC_w(\alpha_\fp(M))$ is a constant; this will finish the proof.
  
  If $n>0$, then we can write~$M= M' E_n$, thus~$\alpha_\fp(M) = \alpha_\fp(M') \alpha_\fp(E_n) \in \Cl(\fp)\frp_w^+$, and $\HC_w(\alpha_\fp(M))=0$. 
  Similarly, if~$m>0$, then~$M=F_1 M''$, hence~$\alpha_\fp(M) \in \frp_w^- \Cl(\fp)$ and again~$\HC_w(\alpha_\fp(M)) = 0$.
  
  Finally, if $M=H_1 \cdots H_k$, then $\alpha_\fp(M)=\alpha_\fp(H_1)\dots\alpha_\fp(H_k)$ is $\frt$-invariant, by Lemma \ref{hc alg hom} and \eqref{eq:HConT}, 
  \[
    \HC_w(\alpha_\fp(M))=\HC_w(\alpha_\fp(H_1))\dots\HC_w(\alpha_\fp(H_k))=(w\rho_\frg-\rho_\frk)(H_1)\dots (w\rho_\frg-\rho_\frk)(H_k)\in \bbC.
\]
\end{proof}

Let $e_i$ and $f_i$ be dual bases for $\frp^+_w$, $\frp^-_w$ as before, let $p = \dim \frp^+_w$.  Define  
\[
P_w = \frac{1}{2^p}  e_1 \cdots  e_{p}  f_p  \cdots  f_1     \in \Cl(\frp).
\]
It is easy to verify that $P_{w}$ is a projection, i.e., $P_{w}^2 = P_{w}$; as we shall see, $P_w$ is closely related to the Harish-Chandra map $\HC_w$.

\begin{lem}\label{proj p+}
Let $\fra^+$ be a~fixed maximal isotropic subspace of~$\fra$ with respect to~$B$. If $\frp$ is odd dimensional
let $Z_1,Z_2 = \frac12 (1 \pm Z_{\mathrm{top}})$, where $Z_{\mathrm{top}} \in q(\twedge^{\mathrm{top}}\fp)$
such that $Z_{\mathrm{top}}^2=1$, be idempotents in the centre of $\Cl(\frp)$.
Furthermore,  let $\ul{ Z_1},\ul{Z_2}$ be analogous idempotents in the centre of $\Cl(\fra)$. 

We construct $\frp^+$ to be $\frp^+_w \oplus \fra^+$ and model $S$ by $\twedge \frp^+$
(respectively $S_1$, $S_2$ by~$\twedge \frp^+Z_1$, $\twedge \frp^+Z_2$).  The projection $P_w$ when considered in $\End(S)$ (respectively $\End(S_1) \oplus \End(S_2)$) is the projection of $S$ (respectively $S_1,S_2)$ to the subspace 
\[
e_{\mathrm{top}}\wedge\twedge\fra^+= e_1\wedge\dots\wedge e_p \wedge \twedge\fra^+ \quad (\text{resp. } e_{\mathrm{top}}\wedge\twedge\fra^+\ul{Z_1}, \quad e_{\mathrm{top}}\wedge\twedge\fra^+\ul{Z_2})
\]
along the span of monomials not containing $e_{\mathrm{top}}=e_1\wedge\dots\wedge e_p$.
\end{lem}
\pf
It is clear that $P_w$ fixes every monomial in $ e_{\mathrm{top}}\wedge \twedge\fra^+$, and kills every other monomial in $\twedge \frp^+$. The claim follows. 
\epf

Define the linear map $p_{w}:\Cl(\frp) \to \Cl(\frp)$ by $p_{w}(x) = P_wx P_w$.
Since $P_w$ commutes with $\Cl(\fra)$, $p_{w}|_{\Cl(\fra)}$ is a nonzero $\mathbb{Z}_2$-algebra morphism. The only $\mathbb{Z}_2$-ideals in $\Cl(\fra)$ are $\Cl(\fra)$ and $\{0\}$, thus $p_{w}|_{\Cl(\fra)}$ is injective, and its image is clearly the subalgebra  $P_w\Cl(\fra)P_w\subseteq\Cl(\frp)$.
Let $\phi_w: P_w\Cl(\fra)P_w \to \Cl(\fra)$ be the inverse algebra isomorphism 
\[\phi_w (P_wa P_w) = a, \quad a \in \Cl(\fra).\]

\begin{lem}\label{lem:HCisprojtowedgeaplus}
  The Harish-Chandra projection $\HC_w$ is equal to $ \phi_w \circ p_{w}$.    
\end{lem}

\begin{proof}
    Since $\frp^+_w P_w = 0$ and $P_w\frp^-_w = 0$, the kernel of $\HC_w$ is contained in the kernel of $p_{w}$. The image of $p_{w}$ is $P_w\Cl(\fra)P_w$, which has dimension $2^{\dim(\fra)}$, just as $\Cl(\fra)=\im\HC_w$. It follows that $\ker p_{w} = \ker \HC_w$. On the other hand, $\phi_w \circ p_{w}$ and $\HC_w$ are both equal to the identity on $\Cl(\fra)$, so they must be the same.
\end{proof}

Clearly, we can think of $p_{w}$ as the projection of $\End(S)$ to $\End( e_{\mathrm{top}} \wedge \twedge\fra^+)\cong\End(\twedge \fra^+)$ (respectively $p_{w}:\End(S_1) \oplus \End(S_2) \to \End(e_{\mathrm{top}}\wedge \twedge\fra^+\ul{Z_1}) \oplus \End(e_{\mathrm{top}}\wedge \twedge\fra^+\ul{Z_2})$).

\begin{prop}\label{prop:HCsurjective}
  For each $w\in W^1$ the corresponding Harish--\/Chandra projection $\HC_w$ restricts to a surjective algebra homomorphism
  \[
    \HC_w\colon \Cl(\fp)^\fk \to \Cl(\fa).
  \]
\end{prop}
\begin{proof} We already know (Lemma~\ref{hc alg hom}) that $\HC_w$ is an algebra homomorphism on
  $\Cl(\frp)^\frk$, so it remains to prove surjectivity.
  Using Lemma~\ref{lem:HCisprojtowedgeaplus} and the fact that $\phi_w$ is an~isomorphism, we see that it is
  enough to  show that $p_{w}$ is surjective.  

Assume first that $\Cl(\frp)=\End(S)$ (i.e., $\dim\frp$ is even). We want to show that any endomorphism of $V=e_{\mathrm{top}} \wedge \twedge\fra^+ $ extends to a $\frk$-morphism of $S$, i.e., to an element of $\Cl(\frp)^\frk$. Let $v_1,\dots,v_m$ be a basis of $V$, i.e., a basis of 
the space of highest weight vectors for the $\frk$-isotypic component of $S$ corresponding to $w$. Let $\psi_{ij}$ be the linear map on $V$ sending $v_i$ to $v_j$ and $v_k$, $k\neq i$, to 0; it is enough to lift such maps since they span $\End(V)$. The highest weight vectors $v_i$ and $v_j$ generate isomorphic copies of the $\frk$-module $E_{w\rho_\frg-\rho_\frk}$ with highest weight $w\rho_\frg-\rho_\frk$, and there is a unique $\frk$-isomorphism of these two copies of $E_{w\rho_\frg-\rho_\frk}$ sending $v_i$ to $v_j$. Now we extend this isomorphism by 0 to the copies of $E_{w\rho_\frg-\rho_\frk}$ generated by $v_k$, $k\neq i$, and to the other $\frk$-isotypic copies of $S$. In this way we get a $\frk$-endomorphism of $S$ which is sent to $\psi_{ij}$ by $p_{w}$.

The argument for $\dim\frp$ odd is similar.  Using the same reasoning as above, 
\[ P_w \Cl(\fra) P_w =\End(e_{\mathrm{top}} \wedge \twedge\fra^+\ul{Z_1}) \oplus \End( e_{\mathrm{top}}  \wedge \twedge\fra^+\ul{Z_2}).\]
Now let $V_1 = e_{\mathrm{top}}\wedge\twedge \fra^+\ul{Z_1}  $ and $V_2 =  e_{\mathrm{top}}  \wedge \twedge \fra^+\ul{Z_2}$.
In an identical way to the even dimension case, given a linear map $\psi_{ij}$ of either $V_1$ or $V_2$ we can lift this map to a $\frk$-module endomorphism of the Spin modules $S_1 = \twedge \frp^+ Z_1$ and $S_2 = \twedge \frp^+ \wedge Z_2$.
\end{proof}

\begin{rem}\label{rem::HCofproj}
  By Lemma~\ref{lem:HCpImAlpha}, the Harish-Chandra homomorphism $\HC_{w}$ sends
  $\Pr(S)=\alpha_\fp(Z(\frk))$, see Definition~\ref{def:PrSigma}, to constants.
  It is easy to see that in terms of projections, $\HC_w$ sends $\pr_w$ to 1, and other $\pr_\sigma$ to 0. 
    
  Indeed, since all elements of $ e_{\mathrm{top}}\wedge \twedge\fra^+$ are of weight $w\rho_\frg-\rho_\frk$,
  this space is contained in the $\frk$-isotypic component of $S$ corresponding to $w$. Thus $\pr_w$ is the
  identity on this space, so $p_{w}(\pr_w)=P_w$ and therefore $\HC_w(\pr_w)=\phi_w(P_w)=1$. On the other hand,
  for $\sigma\neq w$, $\pr_\sigma$ kills the space $ e_{\mathrm{top}} \wedge \twedge\fra^+$, so
  $p_{w}(\pr_\sigma)=0$ and thus also $\HC_w(\pr_\sigma)=0$.
\end{rem}

\section{Clifford algebra analogue of Cartan's theorem for primary and almost primary cases}
\label{sec:Cliff Alg analogue}
In this section we prove an~analogue of Cartan's theorem for $\Cl(\fp)^K$ for primary and almost
primary cases. Throughout this section $(\fg,\fk,K)$ is a~triple from the list~\eqref{eq:TheList}.

\begin{lem}
\label{lem:hc almprim}
  Assume that $(\fg,\fk,K) = (\fsl(2n),\fo(2n),O(2n))$. In this case $|W^1|=2$.
  Let~$\HC_1$ and~$\HC_2$ be the corresponding Harish--\/Chandra homomorphisms
  and~$\pr_1$, $\pr_2$ be the corresponding projections.
  Then the restrictions of~$\HC_1$ and~$\HC_2$ to~$\Cl(\fp)^K$ coincide and define an~isomorphism
  $\HC\colon\Cl(\fp)^K\to \Cl(\fa)$.
\end{lem}
\begin{proof}
  Let $k = \diag(1,\ldots,1,-1)$, so $k$ is a~representative of the disconnected component of~$K$.
  Let $\tilde{k}$ be a~lift of~$k$ in the pin double cover~$\tilde{K}$ of~$K$.
  The two $\fk$-types in~$S$ differ only by the sign of the last coordinate of the highest weight.
  The action of~$\tilde{k}$ switches these two highest weights
  because the action of $k$ on~$\fh$ changes the sign of the last coordinate.
  This implies that $\Ad_k \pr_1 = \pr_2$  and therefore  $\HC_1 = \HC_2\circ\Ad_k$.
  Thus for $\omega \in \Cl(\fp)^K$, $\HC_1(\omega) = \HC_2(\omega)$ since $\Ad_k\omega = \omega$.
  
  Suppose that $\HC(\omega) = 0$. Note that $\omega = \omega \pr_1 + \omega \pr_2$,
  $\HC_1(\omega \pr_2)=0$ and $\HC_2(\omega \pr_1) = 0$, so it follows that
  $\HC_i(\omega \pr_i) = 0$. Since $\HC_i$ is injective on~$\Cl(\fp)^\fk \pr_i$,
  it follows that $\omega \pr_i=0$.
  So $\omega =  \omega \pr_1 + \omega \pr_2 = 0$. This proves injectivity
  of~$\HC_i\colon\Cl(\fp)^K\to\Cl(\fa)$ and surjectivity follows by dimension counting
  (this requires a~version of Lemma~\ref{lem:ClifDecomp}).
\end{proof}

\begin{prop}\label{prop:HCisIsoInPrimary}
  For any $w\in W^1$ the restriction of $\HC_w$ to $\Cl(\fp)^K$ is an algebra isomorphism.
  Moreover, the restrictions of all $\HC_w$ to $\Cl(\fp)^K$ coincide.
\end{prop}
\begin{proof}
  For the primary cases we have that $|W^1|=1$ and
  \(
  \dim \Cl(\fa) = \dim \Cl(\fp)^K = 2^{\dim\fa}
  \),
  so the claim follows from Proposition~\ref{prop:HCsurjective} which states $\HC_w$ is surjective.
  The claim for the remaining almost primary cases $(\fsl(2n),\fo(2n),O(2n))$
  follows from Lemma~\ref{lem:hc almprim}.
\end{proof}

\begin{Def}
  The space of primitive invariants~$P_{\Cl}(\fp)$ in~$\Cl(\fp)^K$ is defined as $q(P_\wedge(\fp))$, where
  $P_\wedge(\fp)$ is as in Definition~\ref{def:PwedgeP}.
\end{Def}

It follows from the relative transgression Theorem~\ref{thm:relTransgression} and Lemma~\ref{lem:quanttr_rel} that
$P_{\Cl}(\fp) = \im\ul{\mathsf{t}}_{\Cl(\fp)^K}$.

\begin{thm}\label{thm:hcofprims}      
  For each $w\in W^1$ the Harish Chandra projection $\HC_w:\Cl(\fp) \to \Cl(\fa)$,
  takes primitives $\phi\in P_{\Cl}(\fp)$ to linear elements of~$\Cl(\fa)$, i.e., to~$q(\twedge^1\fa)$.
  Moreover, $\HC_w( P_{\Cl}(\fp)) = \fa$.
\end{thm}
\begin{proof}
  Lemma~\ref{lem:HCpImAlpha} shows that $\HC_w: \im \alpha_\frp \to \bbC$.
  Theorem~\ref{thm:IotaTClP} states that we have  $\iota_x \phi \in \im \alpha_\frp$ for all $x\in\fp$, hence $\HC_w \iota_x \phi \in \bbC$.
  The linear map~$\HC_w$ commutes with contractions by elements from~$\fra$ hence  
  \[ \iota_x \HC_w (\phi) = \HC_w(\iota_x\phi)  \in \bbC \quad \forall x \in \fra.\] 
  Since~$\phi$ is odd, this implies that $\HC_w(\phi)$ is linear.

  The fact that $\HC_w(P_{\Cl}(\fp)) = \fa$ follows from two facts:
  
  1) $\dim P_{\Cl}(\fp) = \dim \fa$ by~\eqref{eq:dimPwedgeP},

  2) $\HC_w\colon \Cl(\fp)^K\to \Cl(\fa)$ is an algebra isomorphism by Proposition~\ref{prop:HCisIsoInPrimary}.
\end{proof}

\begin{cor}\label{cor:primsanticomm}
  For any two elements $\phi,\psi\in P_{\Cl}(\fp)$,
  the commutator $[\phi,\psi]$ is in~$\bbC$.
\end{cor}
\begin{proof}
  Theorem~\ref{thm:hcofprims} shows that for every $w\in W^1$ the Harish Chandra map~$\HC_w$ takes
  elements of~$P_{\Cl}(\fp)$ to linear terms.
  For any $a_1,a_2\in\fa\subset\Cl(\fa)$, 
  $[a_1,a_2] \in \bbC$, hence the primitives $P_\wedge(\fp)$ have the same property in~$\Cl(\fp)^K$,
  since the Harish Chandra map is an~isomorphism of algebras $\HC_w: \Cl(\fp)^K \to \Cl(\fa)$.
\end{proof}

\subsection{The form on primitives in~$\Cl(\fp)^K$}\label{ss:formonprims}

For $X\in\Cl(\fp)$, let $X^L$ be the operator of left-multiplication by~$X$, and $X^R$ the operator of
($\bbZ_2$-graded) right-multiplication:
\[
  X^L(Y) = XY,\qquad X^R(Y) = (-1)^{p(X)p(Y)}YX.
\]
for $X,Y\in\Cl(\fp)$ of parity $p(X)$ and $p(Y)$, and $X^L - X^R = [X, - ]$.
If $x\in\fp$, we have
\[
  x^L = q\circ(\eps_x + \iota_x) \circ q^{-1},\qquad
  x^R = q\circ(\eps_x - \iota_x) \circ q^{-1}, 
\]
where $\eps_x$ is the operator of left-multiplication by~$x$ in~$\twedge\fp$.

\begin{lem}\label{lem:cliffdeg}
  Let $X\in\twedge^m \fp$, $Y\in\twedge^n \fp$, then
  \[
    q^{-1}([q(X),q(Y)]) \in  \bigoplus_{k= |m-n|}^{m+n}\twedge^k \fp.
  \]
  Moreover, if $m\leq n$ then $ q^{-1}([q(X),q(Y)])_{[n-m]} = (1-(-1)^m) \iota_XY$.
  Here $(-)_{[n-m]}$ denotes the $(n-m)$th graded piece, see Definition~\ref{def:biforms}.
\end{lem}
\begin{proof}
  Recall that $[q(X),q(Y)] = (q(X)^L - q(X)^R)q(Y)$.
 
  We can assume that $X = x_1\wedge \ldots \wedge x_m$, where $x_1,\ldots,x_m\in\fp$ are distinct members of
  an~orthonormal basis. Since in this case
  \[
    q(X) = q(x_1\wedge \ldots \wedge x_m) = q(x_1) \ldots q(x_m),
  \]
  we get that
  \[
    q(X)^L = q(x_1)^L\ldots q(x_m)^L,\qquad
    q(X)^R = q(x_1)^R\ldots q(x_m)^R.
  \]
  Therefore,
  \begin{align*}
    q(X)^L(q(Y)) = {}& q \circ (\eps_{x_1} + \iota_{x_1})\circ \ldots \circ (\eps_{x_m} + \iota_{x_m})(Y),\\
    q(X)^R(q(Y)) = {}& q \circ (\eps_{x_1} - \iota_{x_1})\circ \ldots \circ (\eps_{x_m} - \iota_{x_m})(Y).
  \end{align*}
  When expanding out the brackets $(\eps_{x_i} + \iota_{x_i})$ we find that the resulting operator $q(X)^L$ is a sum 
  of monomials with $k$ multiplication operators and $m-k$ contraction operators, for $k \in \{0,\cdots, m\}$.  Since $\deg \eps_{x_i} = 1$ and $\deg\iota_{x_i} = -1$, then $q(X)^L(q(Y)) \in  \bigoplus_{k= n-m}^{m+n}\twedge^k \fp$. A  similar argument shows that $q(X)^R(q(Y)) \in  \bigoplus_{k= n-m}^{m+n}\twedge^k \fp$, thus 
    \[
    q^{-1}([q(X),q(Y)]) \in  \bigoplus_{k=n-m}^{m+n}\twedge^k(\fp).
  \]

On the other hand, we have that
  \[
    [q(X),q(Y)] = -(-1)^{mn}[q(Y),q(X)]   = -(-1)^{mn} (q(Y)^L-q(Y)^R)q(X).
  \]
  Thus by an identical argument considering operators $q(Y)^L$ and $q(Y)^R$ 
  we find
  \[
    q^{-1}([q(X),q(Y)]) \in  \bigoplus_{k=m-n}^{m+n}\twedge^k(\fp).
  \]

  Thus $q^{-1} ([q(X),q(Y)])$ is concentrated between degrees $|n-m|$ and $n+m$. When $m \leq n$, the $(n-m)$th degree term of $[q(X),q(Y)]$ is equal to 
  \[ \left(q(X)^Lq(Y)\right)_{[n-m]} - \left(q(X)^R q(Y)\right)_{[n-m]} = \iota_X (Y) - (-1)^m \iota_X(Y).\]
  Thus $[q(X),q(Y)]_{[n-m]} = (1 - (-1)^m) q(\iota_X(Y))$.
\end{proof}

\begin{lem}\label{lem:linearterms} 
  Suppose two odd homogeneous elements $X,Y \in \twedge\fp$ are such that  $[q(X),q(Y)] \in \bbC$.
  Then $[q(X),q(Y)] = 2 (\iota_XY)_{[0]}$.
\end{lem}

\begin{proof}
  Suppose that $X\in\twedge^m \fp$, $Y\in\twedge^n\fp$ and $n\neq m$. Then $(\iota_XY)_{[0]} = 0$. On the other hand, by Lemma~\ref{lem:cliffdeg},
  \[
    q^{-1}([q(X),q(Y)]) \in \bigoplus_{k=|m-n|}^{m+n}\twedge^k \fp.
  \]
  In particular, since $|m-n|>0$ we have that $q^{-1}([q(X),q(Y)])_{[0]} = 0$.
  The assumption in the Lemma is that $[q(X),q(Y)] \in \bbC$ hence
  \[
    [q(X),q(Y)] = q^{-1}([q(X),q(Y)])_{[0]} =  0 = 2(\iota_XY)_{[0]}.
  \]
  
  Now suppose that $X,Y \in \twedge^m \fp$ have the same degree, then
  \[
    q^{-1}(q(X)^Lq(Y))_{[0]} = \iota_XY,\qquad
    q^{-1}(q(X)^Rq(Y))_{[0]} = (-1)^m\iota_XY.
  \]
  Hence we get that
  \begin{multline*}
    q^{-1}([q(X),q(Y)])_{[0]} = q^{-1}(q(X)^Lq(Y))_{[0]} - q^{-1}(q(X)^Rq(Y))_{[0]}\\
    = \iota_XY - (-1)^m\iota_XY
    = (1-(-1)^{m})\iota_XY  = (1-(-1)^{m})(\iota_XY)_{[0]}.
  \end{multline*}
  Since $m$ was assumed to be odd, we conclude that $[q(X),q(Y)] = 2(\iota_XY)_{[0]}$.
\end{proof}

Recall that $\tilde{B}:\Cl(\frp) \times \Cl(\frp) \to \bbC$ is the form defined on elements $a,b \in
\Cl(\frp)$ by
\[ \tilde{B}(a,b) = \iota_{q^{-1}(a)}(b)_{[0]}.\]
By abuse of notation, we  denote the restriction of the form~$\tilde{B}$ to~$P_{\Cl}(\fp)$
again by~$\tilde{B}$.
Applying Lemma~\ref{lem:linearterms} and Corollary~\ref{cor:primsanticomm} to the primitives $P_{\Cl}(\fp)$ we find:
\begin{cor}
  Let $\phi,\psi\in P_\wedge(\fp)$ be homogeneous elements. Then
  \[
    [q(\phi), q(\psi)] = 2(\iota_{\phi}\psi)_{[0]}= 2\tilde{B}(\phi,\psi).
  \]
\end{cor}

\subsection{Proof of Theorem \ref{thm:main} for primary and almost primary symmetric pairs. }

To prove statement (a) of Theorem \ref{thm:main}, note first that the elements of
$P_{\Cl}(\fp)\subset\Cl(\fp)^K$ satisfy the relations from Corollary~\ref{cor:primsanticomm}. Thus it follows from the universal property of~$\Cl(P_{\Cl}(\fp),\tilde{B})$ that the inclusion $P_\wedge(\fp) \hookrightarrow \Cl(\fp)^K$ extends to
an~algebra homomorphism \[f\colon \Cl(P_{\Cl}(\fp),\tilde{B})\to \Cl(\fp)^K.\]  
Furthermore, Theorem~\ref{thm:hcofprims}(a) shows that for every $w\in W^1$ \[\HC_w( P_{\Cl}(\fp)) \subseteq \fa \subset\Cl(\fa).\]
By Proposition~\ref{prop:HCisIsoInPrimary} 
we have that $\HC_w|_{\Cl(\fp)^K}$ is an~isomorphism and $\HC_w(P_{\Cl}(\fp)) = \fa \subset \Cl(\fa)$.
We conclude $\HC_w(P_{\Cl}(\fp))$ generates $\Cl(\fa)$ and $P_{\Cl}(\fp)$ generates $\Cl(\fp)^K$. 
In particular,  \[\HC_w\circ f\colon \Cl(P_{\Cl}(\fp),\tilde{B})\to\Cl(\fa)\] is an~algebra homomorphism and $\HC_w(f(P_{\Cl}(\fp)))=\fa$.
The restriction of the form~$\tilde{B}$ to $P_{\Cl}(\fp)$ must therefore be non-degenerate. Indeed, suppose that a~nonzero $\phi\in P_{\Cl}(\fp)$ is orthogonal to
$P_{\Cl}(\fp)$, then $\phi$ anticommutes in~$\Cl(P_{\Cl}(\fp),\tilde{B})$ with all of~$P_{\Cl}(\fp)$ (including itself), so
the nonzero element $\HC_w(f(\phi))\in\fa$ anticommutes with~$\fa$ (including itself).
But there is no such element in~$\fa$. 

(b) When $(\frg,\frk)$ is primary, $\Cl(P_{\Cl}(\frp)) \cong \Cl(\frp)^K$ and $\im \alpha_Z = \Pr(S)$ is
$\bbC$ (see Lemma~\ref{lem:ClifDecomp}). So there is nothing left to prove. In the almost primary cases, $\im\alpha_Z = \Pr(S)$ is two dimensional spanned by $\pr_1$ and $\pr_2$. It is clear that $\im \alpha_Z$ is in the centre of $\Cl(\frp)^\frk$ thus there is a  filtered algebra homomorphism from  $\Cl(P_{\Cl}(\frp))\otimes\Pr(S)$ to $\Cl(\frp)^\frk$. We are left to show that this map is injective.

Let $a_1,a_2 \in \Cl(P_{\Cl}(\frp))$ be such that $a_1 \otimes \pr_1 + a_2 \otimes \pr_2$ is in the kernel, then taking the Harish Chandra projection $hc_{i}$ shows that $a_i=0$, thus this algebra homormorphism is an injection. Dimension counting shows it is a filtered algebra isomorphism.

(c) This is proved in \cite{CNP} and \cite{LRpaper} for any symmetric pair. 

Statement (d) follows from the fact that $\gr \Cl(\frp)^\frk \cong  \left(\twedge\frp\right)^\frk$,
$\gr(\Cl(P_\wedge(\frp)) \cong \twedge P_\wedge(\frp)$, $\gr \bbC[\frt^*]^{W_\fk}/I_\rho \cong
\bbC[\frt^*]^{W_\fk}/I_+$ and $\gr( A \otimes B) \cong \gr A \otimes \gr B$
for any two filtered algebras~$A$ and~$B$.
\qed

\begin{rem}
  In a similar manner, starting from Theorem~\ref{thm:iotaxImAlphaG} and the Harish--\/Chandra
  projection $\HC_\fg\colon \Cl(\fg)\to\Cl(\fh)$, one can prove Kostant's result that $\Cl(\fg)^\fg$ equals
  $\Cl(q(P_\wedge(\fg)),\tilde{B})$ without using Hodge theory for~$\twedge\fg$.  
\end{rem}

\subsection{Contractions by primitive elements}
Let $x,a,b \in \Cl(\fp)$ then 
\begin{align*}
  [x,ab] = {} & xab - (-1)^{|x||a|+ |x||b|}ab x\\
  {}={}& xab - (-1)^{|x||a|}axb + (-1)^{|x||a|}axb - (-1)^{|x||a|+ |x||b|}ab x\\
  {}={}& [x,a]b + (-1)^{|x||a|} a[x,b].
\end{align*}
That is, $[x,-]$ is a derivation of $\Cl(\fp)$.

\begin{cor}\label{cor:cliffordder}
  Let $(\fg,\fk,K)$ be a~triple from the list~\eqref{eq:TheList}.
Let $\phi\in P_\wedge(\fp)$, then the restriction of~$\iota_\phi$ to~$\Cl(\fp)^K$ is equal to
$\frac12[q(\phi), - ]$ and hence is a~derivation of~$\Cl(\fp)^K$.
\end{cor}

\begin{proof}
  Let $\phi_i$ be a basis of $P_\wedge(\fp)$, orthonormal with respect to~$\tilde{B}$, 
 and such that $\phi_i\in( \twedge^{m_i}\fp)^K$.
  For $1\leq i_1<i_2<\ldots < i_k \leq\dim P_\wedge(\fp)$, set $\psi = \phi_{i_1}\wedge \ldots \wedge \phi_{i_k}\in (\twedge^M\fp)^K$,
  where $M = m_{i_1} + \ldots + m_{i_k}$.
  Since the $\phi_i$ are orthonormal, Theorem~\ref{thm:main}(a) implies that $q(\psi) = q(\phi_{i_1})\ldots q(\phi_{i_k})$.

  Using the fact that $[q(\phi_j),-]$ is a derivation of~$\Cl(\fp)$, we have that
  \begin{align*}
    [q(\phi_j),&q(\psi)]  =
     \sum_a(-1)^{a-1}q(\phi_{i_1})\ldots [q(\phi_j),q(\phi_{i_a})]\ldots q(\phi_{i_k})\\
    \intertext{(Since $[q(\phi_a),q(\phi_b)]=2\delta_{a,b}$ and all $i_a$ are distinct.)}
    {}={}& \begin{cases}
      2(-1)^{a-1}q(\phi_{i_1})\ldots\widehat{q(\phi_{i_a)})}\ldots q(\phi_{i_k})
      & \text{if there exists $a$ such that $i_a=j$}, \\
      0 & \text{otherwise}.
    \end{cases}           
  \end{align*}
  Here $\widehat{q(\phi_{i_a})}$ denotes the omission of $q(\phi_{i_a})$.
  We have that $q^{-1}([q(\phi_j), q(\psi)]) \in \left(\twedge^{M-m_j}\fp\right)^K$.
  By Lemma~\ref{lem:cliffdeg} we get that
  \[
    q^{-1}([q(\phi_j),q(\psi)]) = q^{-1}([q(\phi_j),q(\psi)])_{[M-m_i]} = 2\iota_{\phi_j}\psi.
  \]
  Thus $\frac{1}{2}[q(\phi_j),-]$ and $\iota_{\phi_j}$ agree on every monomial of the form
  $q(\phi_{i_1})\ldots q(\phi_{i_k})$ which form a~basis in $\Cl(\fp)^K$,
  hence they agree on~$\Cl(\fp)^K$.
\end{proof}

\begin{cor}\label{cor:primsareder}
  For $\phi\in P_\wedge(\fp)$ the restriction of $\iota_\phi$ to $\left( \twedge\fp \right)^K$ is a derivation of~$\left(\twedge\fp \right)^K$.
\end{cor}

\begin{proof}
  By Corollary~\ref{cor:cliffordder},
  for $\phi\in P_\wedge(\fp)$ the operator $\iota_\phi$ is a~derivation of~$\Cl(\fp)^K$. 
  The associated graded map of $\iota_\phi\in\End(\Cl(\fp))$ is $\iota_\phi\in \End(\twedge\fp)$ when we pass
  from~$\Cl(\fp)$ to the associated graded $\twedge\fp$.  
  Since $\iota_\phi$ is a~derivation of~$\Cl(\fp)^K$ it preserves $\Cl(\fp)^K$, hence its graded version $\iota_\phi$ preserves $(\twedge\fp)^K$.
  Moreover, the associated graded algebra to~$\Cl(\fp)^K$ is~$(\twedge\fp)^K$.  
  The associated graded map of a~derivation is  a~derivation of  the associated graded algebra, thus we get that
  $\iota_\phi$ is a~derivation of~$(\twedge\fp)^K$.
\end{proof}

\begin{rem}
  Let us note that the analogue of Corollary~\ref{cor:primsareder} is a~crucial starting point in Kostant's
  proof of the main result of~\cite{KostantRho}, i.e., $\Cl(\fg)^\fg = \Cl(P_{\Cl}(\fg),\tilde{B})$.
  In our approach it is just a~byproduct of our main result.
\end{rem}

\section{An analogue of Kostant's Clifford algebra conjecture.}\label{sec:KostantCliffordConjecture}

In this section we define a~filtration on~$\fa$ which makes the Harish--\/Chandra
projections~$\HC_w\colon\Cl(\fp)^\fk\to\Cl(\fa)$ filtered algebra morphisms.
This generalises Kostant's Clifford algebra conjecture to the relative case.

We first review some background material about principal $\fsl_2$-triples, for example, see~\cite{Dynkin1952,Kostant1959sl2,Vogan2001SL2}.
Let $\fg$ be a~semisimple Lie algebra. Up to conjugacy, there is a finite number of Lie algebra embeddings
$i\colon \fsl_2\to\fg$. An  $\fsl_2$-embedding is called \emph{principal} if $\fg$ splits into the
smallest number of $i(\fsl_2)$-submodules under the adjoint action:
$\fg \cong \bigoplus_{i=1}^{\rank \fg} L_{2m_i}$, where $m_i$ are the exponents of~$\fg$ and 
$L_{2m_i}$ is the irreducible $\fsl(2)$-module with highest weight~$2m_i$.
Recall that by the Hopf\/--\/Koszul\/--\/Samelson Theorem $P_\wedge(\fg)$ is a~graded vector space
with generators in degrees $2m_i+1 (=\dim L_{2m_i})$.

Let $\check{\fg}$ be the Lie algebra defined by the dual root system and $\check{\fh}\subset\check{\fg}$ be the
corresponding Cartan subalgebra. Kostant observed that $\rho_\fg\in\fh^\ast$ viewed as an~element
of~$\check{\fh}$ coincides with the regular element~$\check{h}$ of the principal $\fsl_2$-embedding
to~$\check{\fg}$ given by $(\check{e},\check{h},\check{f})\subset\check{\fg}$.
Since~$\check{\fh}^\ast$ and~$\fh$ are canonically isomorphic, the coadjoint action of~$\check{e}$ induces
a~filtration on~$\fh$:
\[
  \cF^{(m)} \fh = \left\{ x\in\fh \mid (\ad_{\check{e}}^\ast)^{m+1} x = 0 \right\}.
\]
The dimension of the vector space~$\cF^{(m)} \fh$ jumps at values $m=m_1,\ldots,m_r$ which follows from
results of~\cite{Kostant1959sl2}.

Since $\HC_\fg(q(P_\wedge(\fg))) = \fh$, we can define a~different filtration on~$\fh$:
\[
  F^{(m)} \fh = \HC_\fg(q(P_\wedge^{(2m+1)}(\fg))),
\]
where $P^{(k)}_\wedge(\fg)$ denotes the filtration on~$P_\wedge(\fg)$ induced by the $\bbZ$-grading.
Kostant suggested the following conjecture which was proved for type~$A$ in~\cite{BazlovPhD} and for any reductive Lie algebra
in~\cite{Joseph2012KostantClifford,AlekseevMoreau2012}.

\begin{thm}[Kostant's Clifford algebra conjecture] \label{thm:KostCliffAlgConj}
  The above two filtrations on~$\fh$ coincide, namely, for any $m\in\mathbb{N}$, we have $\cF^{(m)}\fh = F^{(m)}\fh$.
\end{thm}

Similarly, for any $w\in W^1$ define two filtrations on~$\fa$ by
\[
  \mathcal{F}^{(m)}\fa = \left\{x\in\fa \mid (\ad_{\check{e}}^\ast)^{m+1}x = 0 \right\},\quad
  F^{(m)} \fa = \HC_w( P^{(2m+1)}_\wedge(\fp)),
\]
where $P^{(k)}_\wedge(\fp)$ denotes the filtration on~$P_\wedge(\fp)$ induced by the $\bbZ$-grading.
Note that the filtration $F^{(m)}\fa$ does not depend on~$w$ since by Proposition~\ref{prop:HCisIsoInPrimary} we have that
the image of~$P_\wedge(\fp)$ under~$\HC_w$ does not depend on~$w$.

\begin{thm}[Relative Kostant's Clifford algebra conjecture]\label{thm:relKostCliffAlgConj}
  Let $(\fg,\fk)$ be a~symmetric pair of primary or almost primary type, then the two filtrations on~$\fa$
  coincide, namely, for any $m\in\mathbb{N}$, we have $\cF^{(m)}\fa = F^{(m)}\fa$. 
\end{thm}

\begin{proof}
  By definition we have that $P_\wedge(\frp) = \Res^\frg_\frp(P_\wedge(\frg))$,
  so the filtration $F^{(m)}\fra$ is just the restriction of~$F^{(m)}\frh$ to~$\fra$.
  Furthermore, since $\fra \subset \frh$,
  $\cal F^{(m)} \fra$ is the restriction of~$\cal F^{(m)} \frh$ to~$\fra$.
  Thus the theorem follows from restriction to $\fra$ of the results from Theorem~\ref{thm:KostCliffAlgConj};
  the following diagram commutes
 \[   
   \begin{tikzcd}
     \cal F^{(m)} \frh \arrow{r}{\cong} \arrow{d}{\Res^\fg_\fp} & F^{(m)}\frh \arrow{d}{\Res^\fg_\fp} \\
     \cal F^{(m)} \fra \arrow{r}{\cong}  &  F^{(m)} \fra.
   \end{tikzcd}
 \]
 Which concludes the proof.
\end{proof}

\providecommand{\bysame}{\leavevmode\hbox to3em{\hrulefill}\thinspace}
\providecommand{\MR}{\relax\ifhmode\unskip\space\fi MR }
\providecommand{\MRhref}[2]{%
  \href{http://www.ams.org/mathscinet-getitem?mr=#1}{#2}
}
\providecommand{\href}[2]{#2}

\end{document}